\documentclass[10pt]{amsart}
\usepackage{amssymb,amsmath,amsthm,amscd,latexsym}
\usepackage{epsfig, comment,enumerate,epic,bbm}
\usepackage[all]{xy}

\newtheorem{thm}{Theorem}
\newtheorem{cor}[thm]{Corollary}
\newtheorem{lem}[thm]{Lemma}
\newtheorem{prop}[thm]{Proposition}
\newtheorem{defn}[thm]{Definition}

\newtheorem{theorem-question}[thm]{Theorem-Question}

\theoremstyle{definition}
\newtheorem{example}[thm]{Example}
\newenvironment{exafont}{\begin{bf}}{\end{bf}}

\newenvironment{remark}{\vspace{0.3cm}\par\noindent\refstepcounter{thm}\begin{exafont}Remark \thethm.\end{exafont}\hspace{\labelsep}}{\vspace{0.3cm}\par}

\newcommand{\bfc}{{\boldsymbol{{c}}}} \newcommand{\bft}{{\boldsymbol{{t}}}}
\newcommand{\bfa}{\mathbf{a}}\newcommand{\bfb}{\mathbf{b}}\newcommand{\bfm}{\mathbf{m}}
\newcommand{\fc}{{\boldsymbol{{d}}}} \newcommand{\ft}{{\boldsymbol{{u}}}}
 \newcommand{\ut}{\underline{\bft}}
\newcommand{\ttz}{\mathtt{z}}\newcommand{\ttd}{\mathfrak{d}}
\newcommand{\ttu}{\mathtt{u}}
\newcommand{\ttv}{\mathtt{v}}
\newcommand{\ttw}{\mathtt{w}}
\newcommand{\uF}{\underline{F}}
\newcommand{\uN}{\underline{N}}
\newcommand{\um}{\underline{m}}
\newcommand{\uM}{\underline{M}}

\newcommand{\bbO}{\mathbb{O}}\newcommand{\fO}{\mathfrak{O}}
\newcommand{\bbP}{\mathbb{P}}

\newcommand{\bbH}{\mathbb{H}}

\newcommand{\cF}{\mathcal{K}}\newcommand{\bbF}{\mathbb{K}}

\newcommand{\ZZ}{\mathbb{Z}}

\newcommand{\bbT}{\mathbb{T}}

\newcommand{\x}{\mathbf{x}}
\newcommand{\y}{\mathbf{y}}

\DeclareMathOperator{\End}{End}
\DeclareMathOperator{\Hom}{Hom}

\DeclareMathOperator{\Ext}{Ext}

\DeclareMathOperator{\ml}{-mod}
\DeclareMathOperator{\mr}{mod-}
\DeclareMathOperator{\gr}{-gr}
\DeclareMathOperator{\bigra}{-dbigr}
\DeclareMathOperator{\trigra}{-dtrigr}

\DeclareMathOperator{\perf}{-perf}
\DeclareMathOperator{\rperf}{perf-}

\title{The Weyl extension algebra of $GL_2(\overline{\mathbb{F}}_p)$}

\author{Vanessa Miemietz and Will Turner}
\date{}

\begin{document}

\begin{abstract}
We compute the Yoneda extension algebra of the collection of Weyl
modules for $GL_2$ over an algebraically closed field of
characteristic $p>0$ by developing a theory of generalised Koszul duality for certain $2$-functors, one of which controls the rational representation theory of $GL_2$ over such a field.
\end{abstract}

\maketitle

\tableofcontents

%\newpage
\setlength{\parskip}{5pt}

\section{Intro.}

This paper belongs to a sequence of papers exploring structural features of the rational representation theory of the algebraic group $GL_2(F)$,
where $F$ is an algebraically closed field of characteristic $p>0$.
We denote by
$$\mathcal{W} = \{ Sym^\lambda(V) \otimes \omega^\mu \quad | \quad \lambda \in \mathbb{Z}_{\geq 0}, \mu \in \mathbb{Z} \}$$
the set of Weyl modules in the category $G \ml$ of rational representations of $G = GL_2(F)$,
where $V$ denotes the natural two dimensional representation of $G$ and $\omega$ the one dimensional determinant representation.
By definition, $Sym^\lambda(V)$ is the space of symmetric tensors of $V$ of degree $\lambda$,
which is to say the set of fixed points for the symmetric group $S_\lambda$ in its natural action on $V^{\otimes \lambda}$.
$\mathcal{W}$ is a complete set of standard objects in the highest weight category $G \ml$ of rational representations of $G$.
Our objective is to give an explicit description of the Yoneda extension algebra
$${\mathbf W} = \bigoplus_{\Delta, \Delta' \in \mathcal{W}, k \in \mathbb{Z}} \Ext_{G \ml}^k(\Delta,\Delta')$$
of the collection $\mathcal{W}$.
The most important previous discovery in this direction was an algorithm to
compute the dimension of $\Ext_{G \ml}^k(\Delta, \Delta')$ for
$\Delta, \Delta' \in \mathcal{W}$ and $k \geq 0$, written down by A. Parker (\cite{AP}, Theorem 5.1). 
Here we describe the algebra structure.
In previous papers, we have described certain $2$-functors which control the rational representation theory of $G$ and give a combinatorial description of blocks of rational representations.
We have developed a theory of Koszul duality for these operators,
allowing us to give an explicit description of the Yoneda extension algebra of the irreducible $G$-modules.
In this article, we extend these methods to more general homological dualities and a setting where gradings are not necessarily non-negative.

The category $G \ml$ has countably many blocks, all of which are equivalent.
Therefore, the algebra ${\mathbf W}$ is isomorphic to a direct sum of countably many copies of ${\mathbf w}$,
where ${\mathbf w}$ is the Yoneda extension algebra of the Weyl modules belonging to the principal block of $G$.
Our problem is to compute ${\mathbf w}$.

In the following section of the paper we outline a combinatorial description of ${\mathbf w}$ as $\lim_q \fO_F \fO_\Upsilon^q(F[z])$,
where $\fO_{\Upsilon}$ is a certain algebraic operator,
and where $\Upsilon$ is an algebra with monomial basis indexed by elements of a certain lattice polytope and state our main Theorem.
In the remaining sections we prove that this description and our main Theorem is correct.
A consequence of this description of ${\mathbf w}$
is that ${\mathbf w}$ has a monomial basis indexed by elements of an infinite dimensional polytope.
Our analysis proceeds via a more conceptual description of $\Upsilon$ 
involving the tensor algebra over a Koszul algebra $\fc$ of a certain bimodule ${}_{\fc}M_\fc$ (Theorem \ref{dgalgebraandhomology}). A guide through this analysis is provided in Section \ref{outline}, after giving an explicit example in Section \ref{example}, which gets taken up again in Examples \ref{c32}, \ref{exagain}, \ref{Mbasis} and \ref{yuck}.

We write the names of our algebraic operators in exotic scripts like $\mathbb{P}$, $\fO$, $\mathbb{K}$, and $\mathbb{H}$.
These operators are all $2$-functors on certain $2$-categories, a fact we do not prove to avoid frequent checking of axioms, but point out to the interested reader as we go along;
their dominant virtue is their natural behaviour, which means they give a conceptually simple way of encoding complicated homological information.

{\bf Acknowledgements.} The first author acknowledges support from ERC grant PERG07-GA-2010-268109. We would like to thank the referee for helpful suggestions on how to improve the exposition of the article.

\section{Combinatorial description of ${\mathbf w}$.} \label{monomial}

Suppose $\Gamma = \bigoplus_{i,j,k \in \mathbb{Z}} \Gamma^{ijk}$ is a
$\mathbb{Z}$-trigraded algebra. We have a combinatorial operator $\fO_\Gamma$
which acts on the collection of $\mathbb{Z}$-bigraded algebras
$\Sigma$ after the formula
$$\fO_\Gamma(\Sigma)^{ik} = \oplus_{j, k_1+k_2=k} \Gamma^{ijk_1} \otimes_F \Sigma^{jk_2},$$
where we take the super tensor product with respect to the $k$-grading.

There is a trigraded algebra $\Upsilon$ whose structure we partially describe via a polytopal monomial basis;
by a monomial basis we mean a basis in which the product of two elements is $\pm$ another basis element, or zero.
There is a reason for our idleness in giving only a partial definition:
in order to describe ${\mathbf w}$ we only need to know what the algebra structure of $\Upsilon$ looks like when projected onto the subspace
$\Upsilon^{\leq 1} = \oplus_{i, j, k \in \mathbb{Z}, i \leq 1} \Upsilon^{ijk}$.
The space
$\Upsilon^{\leq 1}$ has basis $\{ m_v \}_{v \in \mathcal{P}_{\Upsilon^{\leq 1}} }$ indexed by a polytope $\mathcal{P}_{\Upsilon^{\leq 1}}$
in $\mathbb{Z}^7$ and product (in its projection onto $\Upsilon^{\leq 1}$) given by
$$m_u m_{u'} =
\left\{ \begin{array}{cc}
(-1)^{a j'_0 + b j'_0 + b a'}m_{v} & \textrm{ if } v \in \mathcal{P}_{\Upsilon^{\leq 1}},
v_1 = u_1, u_7 = u_1', u_7'=v_7, \\
& \textrm{ and } v_l = u_l+u'_l \textrm{ for } 2 \leq l \leq 5. \\
0 & \textrm{otherwise.} \end{array} \right.$$
Here we write $u = (u_1,u_2,u_3,u_4,u_5,u_6,u_7) \in \mathbb{Z}^7$.
The $ijk$-degree of a basis element $m_u$ is $(u_2,u_3,u_4)$.

For details of how to define the polytope $\mathcal{P}_{\Upsilon^{\leq 1}}$ and the sign $(-1)^{a j'_0 + b j'_0 + b a'}$
we refer to Section \ref{polybasis}.
For a more elegant conceptual description of $\Upsilon$ we refer to Section \ref{explicit}.

Let us consider the field $F$ as a trigraded algebra concentrated in
degree $(0,0,0)$. We have a natural embedding of trigraded algebras $F
\rightarrow \Upsilon$, which sends $1$ to $m_{(1,0,0,0,0,0,1)}$. This embedding
lifts to a morphism of operators $\fO_F \rightarrow \fO_\Upsilon$.
We have $\fO_F^2 = \fO_F$. Putting these together, we obtain a
sequence of operators
$$\fO_F \rightarrow \fO_F \fO_\Upsilon \rightarrow \fO_F \fO_\Upsilon^2 \rightarrow ...$$
which, applied to the bigraded algebra $F[z,z^{-1}]$ with $z$ placed in $jk$-degree $(1,0)$, gives a sequence of algebra embeddings
$$\mu_1 \rightarrow \mu_2 \rightarrow \mu_3 \rightarrow ...,$$
where $\mu_q = \fO_F \fO_\Upsilon^q(F[z,z^{-1}])$. Taking the union of
the algebras in this sequence gives us an algebra ${\mathbf \mu}$. Our
main theorem, the proof of which will take the entirety of this article, and the structure of which is explained in Section \ref{outline}, is the following:

\begin{thm} \label{Yoneda}
The algebra ${\mathbf w}$ is isomorphic to ${\mathbf \mu}$.
\end{thm}

The algebras $\mu_q$ are isomorphic to Yoneda extension
algebras of collections of Weyl modules for certain Schur algebras, and are consequently finite dimensional.
The $k$-grading on ${\mathbf \mu}$ matches with the natural homological grading on ${\mathbf w}$.

It is straightforward to give a polytopal monomial basis for ${\mathbf \mu}$ as follows:
we define the \emph{height} of a monomial
$m_{v^1} \otimes ... \otimes m_{v^q} \otimes z^\alpha$ in $\Upsilon^{\otimes q} \otimes F[z,z^{-1}]$ to be
$$(v^2_2 - v^1_3,v^3_2 - v^2_3,....,v^q_2 - v^{q-1}_3, \alpha - v^q_3) \in \mathbb{Z}^{q+1}.$$
We have a map from the set of monomials in $\Upsilon^{\otimes q} \otimes F[z,z^{-1}]$
to the set of monomials in $\Upsilon^{\otimes q + 1} \otimes F[z,z^{-1}]$ that sends $a$ to $m_{(1,0,0,0,0,0,1)} \otimes a$ and preserves the last $q+1$ components of the height.
The union over $q$ of monomials in $\Upsilon^{\otimes q} \otimes F[z,z^{-1}]$ of height zero is a polytopal monomial basis for $\mu$.

\section{Example.}\label{example}

The algebra $\mu_q$ is isomorphic to the Yoneda extension
algebra of the Weyl modules for a block of a Schur algebra $S(2,r)$
with $p^q$ simple modules. To place our feet on the ground, let us
give an example of such an algebra.

Let $p=3$. We describe the Yoneda extension algebra of the
collection of Weyl modules in a block of a Schur algebra $S(2,r)$
containing $9$ simple modules.
The $\Ext^1$-quiver of this algebra is given by

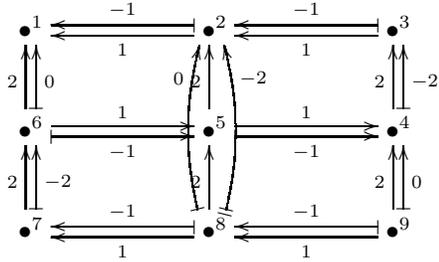
\xymatrix{
\bullet^1 && \bullet^2 \ar@<-2pt>@{|->}[ll]_{-1} \ar@<2pt>[ll]^{1} && \bullet^3 \ar@<-2pt>@{|->}[ll]_{-1} \ar@<2pt>[ll]^{1} \\
\bullet^6 \ar@<-2pt>@{|->}[rr]_{-1} \ar@<2pt>[rr]^{1} \ar@<-2pt>@{|->}[u]_0 \ar@<2pt>[u]^2  && \bullet^5  \ar@<-2pt>@{|->}[rr]_{-1} \ar@<2pt>[rr]^{1} \ar@<2pt>[u]^2 && \bullet^4 \ar@<-2pt>@{|->}[u]_{-2} \ar@<2pt>[u]^2 \\
\bullet^7 \ar@<-2pt>@{|->}[u]_{-2} \ar@<2pt>[u]^2  && \bullet^8 \ar@<-2pt>@{|->}[ll]_{-1} \ar@<2pt>[ll]^{1}  \ar@<2pt>[u]^2  \ar@<-2pt>@{||->}@/_/[uu]_(.75){-2} \ar@<4pt>@{|->}@/^/[uu]^(.75)0 && \bullet^9 \ar@<-2pt>@{|->}[ll]_{-1} \ar@<2pt>[ll]^{1} \ar@<-2pt>@{|->}[u]_0 \ar@<2pt>[u]^2
}

where the number of bars through the tail of the arrow marks the degree of the arrow in the $\Ext$-grading (aka $k$-grading)
and the number labelling the arrow denotes its $j$-grading.

The composition structure of the projectives, with the superscript denoting the $k$-degree and the subscript denoting the $j$-degree, is given by

{\footnotesize 
$1^0_0$ $\qquad$ $\qquad$
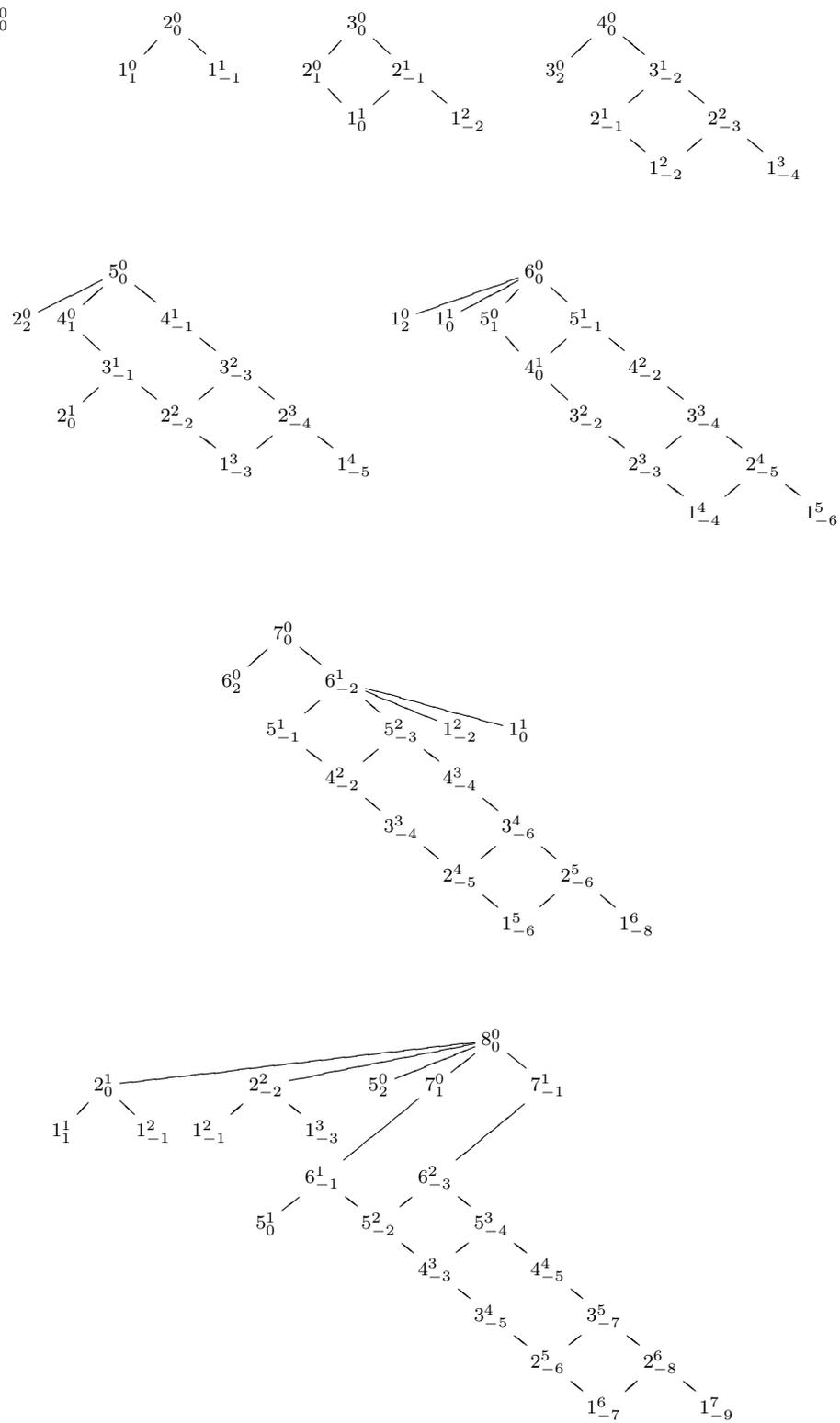
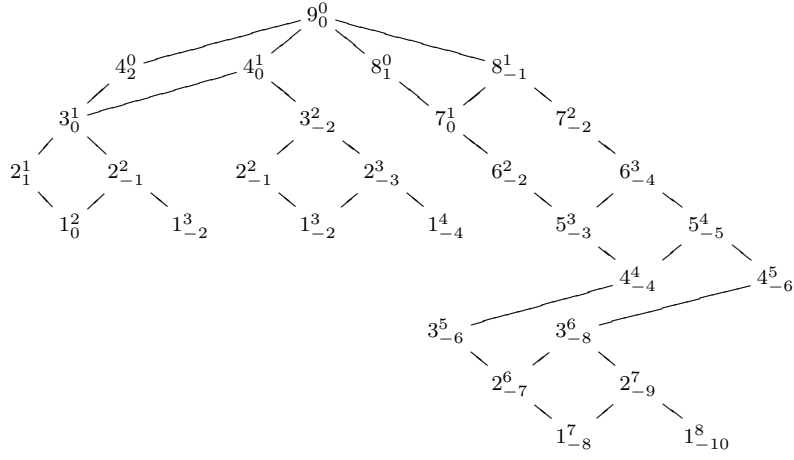
\xymatrix@=5pt{ & 2^0_0 \ar@{-}[dl]\ar@{-}[dr]& \\1^0_1 &&1^1_{-1} } $\qquad$ 
\xymatrix@=5pt{ & 3^0_0\ar@{-}[dl]\ar@{-}[dr] & & \\2^0_1 \ar@{-}[dr]&& 2^1_{-1}\ar@{-}[dl]\ar@{-}[dr] & \\&1^1_0 &&1^2_{-2} }
$\qquad$ 
\xymatrix@=5pt{ &4^0_0\ar@{-}[dl]\ar@{-}[dr]& &&\\ 3^0_{2}&&3^1_{-2} \ar@{-}[dl]\ar@{-}[dr]\\ &2^1_{-1} \ar@{-}[dr]&&2^2_{-3}\ar@{-}[dl]\ar@{-}[dr]& \\ && 1^2_{-2}&&1^3_{-4}  }

\medskip
$$\xymatrix@=5pt{ && 5^0_0 \ar@{-}[dll]\ar@{-}[dl]\ar@{-}[dr]&&&& \\
                2^0_2 & 4^0_1 \ar@{-}[dr]&& 4^1_{-1} \ar@{-}[dr]&&&\\
                && 3^1_{-1} \ar@{-}[dl]\ar@{-}[dr]&& 3^2_{-3} \ar@{-}[dl]\ar@{-}[dr]&&\\
                & 2^1_0 && 2^2_{-2} \ar@{-}[dr]&& 2^3_{-4} \ar@{-}[dl]\ar@{-}[dr]&\\
                 &&&& 1^3_{-3} && 1^4_{-5} } \textrm{   }
\xymatrix@=5pt{ &&& 6^0_0 \ar@{-}[dlll]\ar@{-}[dll]\ar@{-}[dl]\ar@{-}[dr]&&&&&\\
                1^0_2 & 1^1_0 & 5^0_1 \ar@{-}[dr]&& 5^1_{-1} \ar@{-}[dl]\ar@{-}[dr]&&&&\\
                &&& 4^1_{0} \ar@{-}[dr]&& 4^2_{-2} \ar@{-}[dr]&&&\\
                &&&& 3^2_{-2} \ar@{-}[dr]&& 3^3_{-4} \ar@{-}[dl]\ar@{-}[dr]&&\\
                &&&&& 2^3_{-3} \ar@{-}[dr]&& 2^4_{-5} \ar@{-}[dl]\ar@{-}[dr]& \\
                &&&&&& 1^4_{-4} && 1^5_{-6}}$$

\medskip
$$\xymatrix@=5pt{ && 7^0_0 \ar@{-}[dl]\ar@{-}[dr] &&&&&& \\
                & 6^0_2  && 6^1_{-2} \ar@{-}[dl]\ar@{-}[dr]\ar@{-}[drr]\ar@{-}[drrr]&&&&&\\
                && 5^1_{-1} \ar@{-}[dr] && 5^2_{-3} \ar@{-}[dl]\ar@{-}[dr]& 1^2_{-2} & 1^1_{0} &&   \\
                &&& 4^2_{-2} \ar@{-}[dr]&& 4^3_{-4} \ar@{-}[dr]&&&\\
                &&&& 3^3_{-4} \ar@{-}[dr]&& 3^4_{-6} \ar@{-}[dl]\ar@{-}[dr]&&\\
                &&&&& 2^4_{-5} \ar@{-}[dr]&& 2^5_{-6} \ar@{-}[dl]\ar@{-}[dr]&\\
                &&&&&& 1^5_{-6} && 1^6_{-8} }$$

\medskip
$$\xymatrix@=4pt{ &&&&&&&& 8^0_0 \ar@{-}[dlllllll]\ar@{-}[dllll]\ar@{-}[dll]\ar@{-}[dl]\ar@{-}[dr]&&&&&&& \\
                & 2^1_0 \ar@{-}[dl]\ar@{-}[dr]&&& 2^2_{-2} \ar@{-}[dl]\ar@{-}[dr]&& 5^0_2 & 7^0_1 \ar@{-}[ddll]&& 7^1_{-1} \ar@{-}[ddll]&&&&&& \\
                1^1_1 && 1^2_{-1} & 1^2_{-1} && 1^3_{-3} &&&&&&&&&& \\
                &&&&& 6^1_{-1} \ar@{-}[dl]\ar@{-}[dr]&& 6^2_{-3} \ar@{-}[dl]\ar@{-}[dr]&&&&&  \\
                &&&& 5^1_0 && 5^2_{-2} \ar@{-}[dr]&& 5^3_{-4}\ar@{-}[dl]\ar@{-}[dr] &&&& \\
                &&&&&&& 4^3_{-3} \ar@{-}[dr] && 4^4_{-5} \ar@{-}[dr]&&& \\
                &&&&&&&& 3^4_{-5} \ar@{-}[dr]&& 3^5_{-7} \ar@{-}[dl]\ar@{-}[dr]&& \\
                &&&&&&&&& 2^5_{-6} \ar@{-}[dr]&& 2^6_{-8} \ar@{-}[dl]\ar@{-}[dr]& \\
                &&&&&&&&&& 1^6_{-7} && 1^7_{-9}  }$$

\medskip
$$\xymatrix@=5pt{ &&&&&& 9^0_0 \ar@{-}[dlll] \ar@{-}[dl]\ar@{-}[dr]\ar@{-}[drrr]&&&&&&&&&& \\
                &&& 4^0_2 \ar@{-}[dl]&& 4^1_0\ar@{-}[dlll] \ar@{-}[dr] && 8^0_1\ar@{-}[dr] && 8^1_{-1} \ar@{-}[dl]\ar@{-}[dr]&&&&&&&\\
                && 3^1_0 \ar@{-}[dl]\ar@{-}[dr]&&&& 3^2_{-2}\ar@{-}[dl]\ar@{-}[dr] && 7^1_0 \ar@{-}[dr] && 7^2_{-2} \ar@{-}[dr] &&&&&& \\
                & 2^1_1 \ar@{-}[dr]&& 2^2_{-1} \ar@{-}[dl]\ar@{-}[dr]&& 2^2_{-1} \ar@{-}[dr]&& 2^3_{-3} \ar@{-}[dl]\ar@{-}[dr]&& 6^2_{-2} \ar@{-}[dr]&& 6^3_{-4} \ar@{-}[dl]\ar@{-}[dr]&&&&&\\
                && 1^2_0 && 1^3_{-2} && 1^3_{-2} && 1^4_{-4} && 5^3_{-3} \ar@{-}[dr]&& 5^4_{-5} \ar@{-}[dl]\ar@{-}[dr]&&&&&\\
                &&&&&&&&&&& 4^4_{-4} \ar@{-}[dlll]&& 4^5_{-6} \ar@{-}[dlll]&&& \\
                &&&&&&&& 3^5_{-6} \ar@{-}[dr]&& 3^6_{-8} \ar@{-}[dl]\ar@{-}[dr]&& \\
                &&&&&&&&& 2^6_{-7} \ar@{-}[dr]&& 2^7_{-9}  \ar@{-}[dl]\ar@{-}[dr]&\\
                &&&&&&&&&& 1^7_{-8} && 1^8_{-10} }$$
}
\begin{remark}
We remark that the algorithm for the computation of 
$\Ext^k(\Delta, \Delta')$ described in \cite[Theorem 5.1]{AP} produces a direct sum of extension groups of lower Yoneda degrees, but not all summands actually describe extensions of modules in the same block (hence giving zero) and in general many summands will produce zero. 
Since our focus is on a fixed block, the algorithm cannot be deduced from our results, and due to the extra zero summands in her description, it does not seem possible to deduce a refined version (focussing on a single block) of her algorithm either. This is illustrated by the following:

Suppose we were just interested in the dimension of $\Ext^7(\Delta_{0},\Delta_{8p})$ where by $\Delta_{m}$ we denote the standard module attached to the regular weight $m$. It is well-known when such modules are in the same block (see e.g. \cite[Section 1]{AP}) and our two chosen weights are in the same block, and are the first and ninth weights in this block repectively, hence translate to $\Delta(1), \Delta(9)$ in our setup and from the above example we easily see that $\Ext^7(\Delta(1), \Delta(9))$ is one-dimensional (there is precisely one $1$ with superscript $7$ occurring in the ninth projective). The algorithm in \cite[Theorem 5.1]{AP} tells us that $\Ext^7(\Delta_{0},\Delta_{8p})\cong \Ext^6(\Delta_{0},\Delta_{6} ) \oplus \Ext^4(\Delta_{0},\Delta_{4} ) \oplus \Ext^2(\Delta_{0},\Delta_{2} ) \oplus \Ext^0(\Delta_{0},\Delta_{0} )$. Now $\Delta_{2}$ is in a different block, and $\Ext^4(\Delta_{0},\Delta_{4} ) = \Ext^4(\Delta(1),\Delta(2))=0$, as well as $\Ext^6(\Delta_{0},\Delta_{6} ) = \Ext^6(\Delta(1),\Delta(3)) =0$, and only  $\Ext^0(\Delta_{0},\Delta_{0} )$ is one-dimensional, giving our desired dimension.
\end{remark}

\section{Outline.}\label{outline}

Our general approach to our problem is to use our previous description of the category of rational representations of $GL_2$ as the module category of an iterative application of certain operators to the ground field (together with its regular bimodule taken twice). These operators include $\bbP_{a,\um}$, which depend on a (differential bigraded) algebra $a$ with pairs of bimodules $\um$ (in our case a very small quasi-hereditary algebra, (a projective resolution of) its tilting module and the adjoint thereof), and a variant $\mathfrak{O}$ of these which has nicer properties with respect to taking homology. In section \ref{theory} we describe some theory of these algebraic operators, recalling and generalising some previous definitions and results \cite{MT3}. In particular we show that under suitable conditions, dg-derived equivalences in the sense of Keller's theory \cite{Ke}, which we describe by an algebraic operator $\bbF$, behave nicely with respect to our algebraic operators, allowing us to pull the dg-derived equivalence past while acting on the index, i.e. $\bbF\bbP_{a,\um} = \bbP_{\bbF(a,\um)}\bbF$ (see Theorem \ref{chainquasi}). 

We then recall some facts about the representation theory of $GL_2$ in Section \ref{GL2} and introduce the main objects of interest, the algebras that we call $\bfc$ and their tilting modules $\bft$ (or rather a pair $\ut$ of projective bimodule resolutions of this and its adjoint). In particular, computing $\bbF\bbP^q_{\bfc,\ut}(F,(F,F))$ for any $q$ will become our goal.
In Section \ref{hypos}, we verify that in this case the conditions of Theorem \ref{chainquasi} are satisfied, so we can ideed apply our theory. We introduce the explicit little algebra $\fc$ (Section \ref{psi}), which is the algebra component of $\bbF(\bfc,\ut)$, and proceed to trace $\ut$ through this duality (producing a pair of dg bimodules which we call $\underline{\ft}$), which takes us Sections \ref{general} and \ref{tthroughK}, to finally give a description of $\bbF(\bfc,\ut)$ in Proposition \ref{daggermaltese}.

In Section \ref{computation}, we go through the reduction process, showing that applying Theorem \ref{chainquasi} recursively, we have $\bbF\bbP^q_{\bfc,\ut} (F,(F,F))=\bbP^q_{\fc,\underline{\ft}} (F,(F,F))$. Since we are interested in the homology of the resulting dg-algebra, we use results from Section \ref{comparison} to show that we can use our operators $\mathfrak{O}$ instead, with the homology of the tensor algebra $\mathbb{T}_\fc(\underline{\ft})$ as input. 

We then spend Section \ref{tensoralg} explicitly computing this homology $\mathbb{HT}_\fc(\underline\ft)$, which we call $\Upsilon$. In order to achieve this we examine a certain bimodule $M$, which is prominent in this homology algebra and has many intriguing properties, see Propositions \ref{Ms}, \ref{Ms2}. Finally in Theorem \ref{basisupsilon} we give the combinatorial description of $\Upsilon$ referred to in Section \ref{monomial}, which then also completes the proof of Theorem \ref{Yoneda}.
We have relegated certain generalities 
on signs and Koszul duality to appendices.

The reader mainly interested in actions on $2$-categories will be most interested in Section \ref{theory}, whereas the reader interested in looking for explicit combinatorial data about $GL_2$ might wish to skip this section, and only look up the definitions of the relevant operators when needed. This reader could start with Section \ref{GL2}, check out Section \ref{psi} and the end of Section \ref{tthroughK} for the definitions of $\fc$ and $\ft$, proceed to Section \ref{computation} to see what needs to be done, and then concentrate on the description of $\Upsilon$.

\section{Homological duality for algebraic operators.}\label{theory}

\subsection{Grading conventions.} \label{grs}

For our computation of the Yoneda extension algebras of simple modules for $GL_2$ in \cite{MT3}, we used trigraded structures. Here, we
also use quadragraded structures
$$S = \bigoplus_{d,i,j,k \in \mathbb{Z}} S^{dijk}.$$
It will be necessary to differentiate between the four gradings; we
will call them the $d$-grading, the $i$-grading, the $j$-grading,
and the $k$-grading.

As in our previous paper, the $i$ and $j$-gradings will be
algebraic. We denote by $\langle 1
\rangle$ a shift by $1$ in the $j$-grading, thus $(M \langle n
\rangle)_j = M_{j-n}$.

The $k$-grading will always be a homological grading, and differentials always have
$k$-degree $1$, and $(i,j)$-degree $(0,0)$.
When we speak of a differential (bi-, tri-, quadra-) graded algebra,
we mean (bi-, tri-, quadra-)graded algebra which is a differential
graded algebra with respect to the $k$-grading. We denote by
$\mathbb{H}$ the cohomology functor, which takes a differential
$k$-graded complex $C$ to the $k$-graded vector space $\mathbb{H} C
= H^{\bullet} C$. We denote by $[1]$ a shift by $1$ in the
$k$-grading.

The $d$-grading will be algebraic, and in our application will form
a $\Delta$-grading and will in all relevant cases be positive.

In practice, for the development of our theory, one of the gradings will in fact be obscured,
so we only have to consider trigraded structures $S^{djk}$ or $S^{ijk}$;
nevertheless it seems to be important to distinguish between the four types of grading.

There is a fifth grading. When we discuss Koszul duality,
we will refer to the grading of a quotient of a quiver algebra by path length, which we call it the $r$-grading.
In certain special cases, this $r$-grading coincides with our $j$-grading, but in other cases it does not;
this is why we denote it with a different letter.
But outside our discussion of Koszul duality, we will not be concerned with the $r$-grading, and we disregard it.
Likewise, the sixth and seventh gradings (the $f$- and $h$-gradings).

{\bf Given a mistake in Lemma \ref{bondedtensor} below, various statements in the rest of this section are incorrect. Their appropriate replacements are given in Appendix \ref{errata}. The reader is advised to ignore Section 5.}

\subsection{Bonded pairs of bimodules.}\label{bonded} Let $A$ be a finite dimensional algebra.
We say a pair $\uM= (M, M')$ of
$A$-$A$-bimodules are \emph{bonded} if we have homomorphisms $M
\otimes_A M' \rightarrow A$ and $M' \otimes_A M \rightarrow A$, such that
the resulting pair of maps
$$M \otimes_A M' \otimes_A M \rightarrow M$$
are equal, and the resulting pair of maps
$$M' \otimes_A M \otimes_A M' \rightarrow M'$$
are equal.

\begin{lem}\label{bondedtensor}
Given a bonded pair $\uM$ of $A$-$A$-bimodules, the space
$$\mathbb{T}_A(\uM)=(\bigoplus_{i \geq 1} M^{\otimes_A i}) \oplus A \oplus (\bigoplus_{i \leq -1} M'^{\otimes_A -i})$$
is a $\mathbb{Z}$-graded algebra, with product given by the natural bimodule homomorphisms
$$M^{\otimes_A i_1} \otimes M'^{\otimes_A -i_2} \rightarrow M^{\otimes_A i_1+i_2}, \quad i_1 \geq i_2,$$
$$M^{\otimes_A i_1} \otimes M'^{\otimes_A -i_2} \rightarrow M'^{\otimes_A -i_1-i_2}, \quad i_2 \geq i_1,$$
$$M'^{\otimes_A -i_1} \otimes M^{\otimes_A i_2} \rightarrow M'^{\otimes_A -i_1-i_2}, \quad i_1 \geq i_2,$$
$$M'^{\otimes_A -i_1} \otimes M^{\otimes_A i_2} \rightarrow M^{\otimes_A i_2+i_1}, \quad i_2 \geq i_1,$$
obtained by applying the maps $M \otimes_A M' \rightarrow A$, $M' \otimes_A M \rightarrow A$ a number of times;
here we write $M^{\otimes_A 0} = M'^{\otimes_A 0} = A$
\end{lem}
\proof We define $\mathcal{A}$ to be the tensor algebra of $M \oplus M'$ over $A$,
modulo relations implying the product of $M$ and $M'$ lies in $A$, and the product map on these bimodules is given by the maps
$M \otimes_A M' \rightarrow A$, $M' \otimes_A M \rightarrow A$.
The algebra $\mathcal{A}$ acts naturally on the space $\mathbb{T}_A(\uM)$ via the homomorphisms
$M \otimes_A M' \rightarrow A$, $M' \otimes_A M \rightarrow A$;
the fact this is an algebra action follows from the fact that $M$ and $M'$ are bonded.
The relations in $\mathcal{A}$ imply that $\mathcal{A}$ is a quotient of $\mathbb{T}_A(\uM)$ as an $A$-$A$-bimodule.
The subspace $A$ of $\mathbb{T}_A(\uM)$ generates $\mathbb{T}_A(\uM)$ as an $\mathcal{A}$-module,
which implies $\mathbb{T}_A(\uM)$ is a quotient of $\mathcal{A}$.
Consequently $\mathcal{A}$ acts freely on $\mathbb{T}_A(\uM)$, and we can identify $\mathcal{A}$ with $\mathbb{T}_A(\uM)$;
once we do this we obtain a product on $\mathbb{T}_A(\uM)$ as advertised.
\endproof

We have a ready supply of bonded pairs of bimodules:

\begin{lem}
Suppose $M$ is a differential graded $A$-$A$-bimodule which is projective on the left and right as an $A$-module.
Then $M$ and $\Hom_A(M, A)$ are a bonded pair of dg bimodules.
\end{lem}
\proof
We have a natural isomorphism of dg bimodules $$M \rightarrow \Hom_{A}(\Hom_A(M,A),A).$$
Our bonded structure is given by the pair of natural maps
$$M \otimes \Hom_A(M,A) \rightarrow A, \quad \Hom_A(M,A) \otimes \Hom_A(\Hom_A(M,A),A)\rightarrow A.$$
We check the natural composition maps
\begin{equation*}
\begin{split}
M \otimes_A \Hom_A(M,A) &\otimes_A \Hom_A(\Hom_A(M,A),A) \\&\rightarrow M \otimes_A A \rightarrow M \rightarrow \Hom_A(\Hom_A(M,A),A)
\end{split}
\end{equation*}
and
\begin{equation*}
\begin{split}
M \otimes_A& \Hom_A(M,A) \otimes_A \Hom_A(\Hom_A(M,A),A) \\&\rightarrow A \otimes_A \Hom_A(\Hom_A(M,A),A) \rightarrow \Hom_A(\Hom_A(M,A),A)
\end{split}
\end{equation*}
are equal. If we identify these with maps
\begin{equation*}
\begin{split}\Hom_A(A,M) \otimes_A \Hom_A(M,A) \otimes_A \Hom_A(A,M)\rightarrow \Hom_A(\Hom_A(M,A),\Hom_A(A,A))\end{split}
\end{equation*}
then both send $\alpha \otimes \beta \otimes \gamma \in \Hom_A(A,M) \otimes_A \Hom_A(M,A) \otimes_A \Hom_A(A,M)$
to the morphism sending $\eta \in \Hom_A(M,A)$ to $\alpha \beta \gamma \eta \in \Hom_A(A,A)$.
\endproof

\subsection{The collections $\mathcal{T}$ and $\mathcal{U}$.}\label{2cats}
Let $\mathcal{T}$ denote the
collection of dg algebras with a pair of bonded dg bimodules
$$\mathcal{T} = \left\{
\begin{array}{rr}
  (A,\uM) | & A = \bigoplus_k A^k \textrm{a dg algebra, } \uM = (\bigoplus_k M^k,\bigoplus_k M'^k) \\
& \textrm{bonded dg } A \textrm{-} A \textrm{-bimodules}.
\end{array} \right\}$$

Let $\mathcal{U}$ denote the
collection of differential bigraded algebras with a pair of bonded
differential bigraded bimodules
$$\mathcal{U} = \{ (A,\uM) | \quad A = \bigoplus_{k \in \ZZ, d \geq 0} A^{dk}, \uM=(M,M') =(\bigoplus_{k \in \ZZ, d \geq 0} M^{dk}, \bigoplus_{k \in \ZZ, d \geq 0} M'^{dk}) \}.$$
%where  $A$ is concentrated in positive $d$-degrees and we require isomorphisms $A^{0 \bullet} \otimes_A M \cong M \otimes_A A^{0 \bullet} \cong M^{0 \bullet}$ and similarly for $M'$.

For the reader interested in $2$-categories, we note that these collections indeed form $2$-categories as follows:

\begin{itemize}
\item Objects are given by the collections defined above;
\item $1$-morphisms between two objects $(A,\uM)$ and $(B,\uN)$ are given by a triple $(S, \phi_S,\phi_S')$, consisting of a dfferential (bi-)graded $A$-$B$-bimodule ${}_AS_B$ and quasi-isomorphisms
$$\phi_S: S \otimes_BN\to M\otimes_AS  \qquad \textrm{ and } \qquad  \phi'_S: S \otimes_BN'\to M'\otimes_AS $$
such that the diagrams 

$$\xymatrix@C=20pt{S\otimes_BN\otimes_BN'\ar^{\phi_S\otimes id }[d]\ar[r] & S \otimes_B B \ar^{\sim}[r]& S \ar^{id}[dd] & S\otimes_BN'\otimes_BN\ar[r]  \ar^{ \phi'_S\otimes id }[d]& S \otimes_B B \ar^{\sim}[r]& S\ar^{id}[dd] \\
M\otimes_AS\otimes_BN' \ar^{id \otimes \phi'_S}[d]&& & M'\otimes_AS\otimes_BN\ar^{id \otimes \phi_S}[d]\\
M\otimes_AM'\otimes_AS \ar[r]&A\otimes_AS\ar^{\sim}[r] & S & M'\otimes_AM\otimes_AS \ar[r]&A\otimes_AS\ar^{\sim}[r] & S \\
}$$
commute, where the first horizontal morphism in each row of each diagram is given by the corresponding bonding map and the second horizontal morphism in each row of each diagram is just the canonical isomorphism;
\item $2$-morphisms from $(S, \phi_S,\phi_S')$ to $(T, \phi_T,\phi_T')$ are given by homomorphisms of dfferential (bi-)graded $A$-$B$-bimodules $f: S \to T$ auch that the diagrams
$$\xymatrix{
S \otimes_BN\ar^{\phi_S}[r]\ar^{ f\otimes id }[d]&M \otimes_A S \ar^{id \otimes f}[d] &&S \otimes_BN'\ar^{ f\otimes id }[d]\ar^{\phi'_S}[r]& M' \otimes_A S \ar^{id \otimes f}[d]\\
T \otimes_BN\ar^{\phi_T}[r]                      &M \otimes_A ST                      &&T \otimes_BN'\ar^{\phi'_T}[r]& M' \otimes_A T 
}$$
commute.
\end{itemize}

\subsection{Dg equivalences and quasi-isomorphisms.}\label{dgeqs}

Our main interest in this article is using Keller's Morita theory for dg derived categories \cite[Theorem 3.10]{Ke}, in order to compute the dg algebra of which our desired extension algebra is the homology. Since the algebra to which we wish to apply Keller's equivalence is constructed via an iteration which drags bimodules around, we need to extend the well-known notion to this more cumbersome setting.

We define a \emph{Rickard object} of $\mathcal{U}$ (or $\mathcal{T}$) to be an object
$(A,\uM)$ of $\mathcal{U}$ (or $\mathcal{T}$), where ${}_AM_A \in A\perf \cap \rperf A$, the homomorphisms $M \otimes_A M' \rightarrow A$ and $M'
\otimes_A M \rightarrow A$ are quasi-isomorphisms and there is a quasi-isomorphism $A \to \mathbb{H}A$. In this case we will call $M$ invertible and often denote $M'$ by $M^{-1}$. The naming of these objects is motivated by the fact that in this case the bimodules $M$ and $M'$ induce mutually inverse self-equivalences of $D_{dg}(A)$ as studied by Rickard in the ungraded setting \cite{Rickard}.

\begin{example}
If $a$ is a quasi-hereditary Ringel self-dual algebra, and $m$ its tilting module, we can consider the object $(a,(m, \Hom_a(m,a)))$ in $\mathcal{T}$ which is a Rickard object.
\end{example}

Let $(A,\uM)$ and $(B,\uN)$ be objects of $\mathcal{T}$, where $\uM = (M,M')$ and $\uN = (N,N')$.

\begin{defn}\label{dgequiv} \cite[Definition 3]{MT3} A \emph{dg equivalence} between objects $(A,\uM)$ and
$(B,\uN)$ of $\mathcal{T}$ is

(i) a dg $A$-$B$-bimodule $X$ such that ${}_AX$ belongs to $A \perf$,
such that ${}_AX$ generates $D_{dg}(A)$, and the natural map
$$B \rightarrow \End_A(X)$$ is a quasi-isomorphism;

(ii) quasi-isomorphisms
$$X \otimes_B N \rightarrow M \otimes_A X,$$
$$X \otimes_B N' \rightarrow M' \otimes_A X,$$
such that the resulting diagrams of maps
$$\xymatrix{
X \otimes_B N \otimes_B N' \ar[r] \ar[d] & M \otimes_A \otimes M' \otimes_A X \ar[r] & A \otimes_A X \ar[d] \\
X \otimes_B B \ar[rr] & & X
}$$
$$\xymatrix{
X \otimes_B N' \otimes_B N \ar[r] \ar[d] & M' \otimes_A \otimes M \otimes_A X \ar[r] & A \otimes_A X \ar[d] \\
X \otimes_B B \ar[rr] & & X
}$$
commute. If there is a dg equivalence between $(A,\uM)$ and $(B,\uN)$,
we write $(A,\uM) \gtrdot (B,\uN)$. 
\end{defn} 

In other words, a dg equivalence between $(A,\uM)$ and
$(B,\uN)$ is a $1$-morphism between them, such that the bimodule given in the data of the $1$-morphism induces an equivalence between $D_{dg}(A)$ and $D_{dg}(B)$.
If we have a dg equivalence between $(A,\uM)$ and
$(B,\uN)$, we have a derived equivalence
$$\xymatrix{ D_{dg}(A) \ar@/^/[r]^{\Hom_A(X,-)} & \ar@/^/[l]^{X \otimes_B -} D_{dg}(B)},$$ and diagrams
which commute up to a natural isomorphism:
$$\xymatrix{
D_{dg}(B)\ar^{X \otimes_B -}[r]\ar^{N(') \otimes_B -}[d] & D_{dg}(A)\ar^{M(') \otimes_A-}[d]\\
D_{dg}(B)\ar^{X \otimes_B -}[r] & D_{dg}(A). }$$ We define a
\emph{quasi-isomorphism} from $(A,\uM)$ to $(B,\uN)$ to be a
quasi-isomorphism $A \rightarrow B$, along with compatible
quasi-isomorphisms $${}_AM_A \rightarrow {}_BN_B, \quad {}_AM'_A \rightarrow {}_BN'_B$$ such that the diagrams
$$\xymatrix{ M \otimes_A M' \ar[r]\ar[d]& N\otimes _BN'\ar[d]\\ A \ar[r] & B   } \textrm { and }\xymatrix{ M' \otimes_A M \ar[r]\ar[d]& N'\otimes _BN\ar[d]\\ A \ar[r] & B   }$$
commute, where the horizontal maps are given by above quasi-isomorphisms, and vertical maps by the bondings.

\subsection{The operator $\mathbb{P}$.}\label{bbP}
In order to define the algebras describing blocks of $GL_2$, we will need certain algebraic operators, which, when applied repeatedly to the object $(F,(F,F)) \in \mathcal{U}$, produce the desired algebras, for details please see Section \ref{GL2}.
In order to define these operators, we call an object $(a,\um)$ of $\mathcal{U}$ a \emph{$j$-graded
object}, if $a= \bigoplus a^{djk}$ is a differential trigraded algebra,
and $\um = (m,m') = (\bigoplus m^{djk},\bigoplus m'^{djk})$ a bonded pair of differential trigraded
$a$-$a$-bimodules. Given a $j$-graded object of $\mathcal{U}$, we
have an operator
$$\mathbb{P}_{a,\um} \circlearrowright \mathcal{U}$$
given by
$$\mathbb{P}_{a,\um}(A,\uM) = (\bigoplus a^{djk} \otimes_F M^{\otimes_A j}, (\bigoplus m^{djk} \otimes_F \uM^{\otimes_A j},\bigoplus m'^{djk} \otimes_F \uM^{\otimes_A j})),$$
where for $j>0$, $\uM^{\otimes_A j} = M^{\otimes_A j}$ and for $j<0$, $\uM^{\otimes_A j} = M'^{\otimes -j}$.
The algebra structure on $\bigoplus a^{djk} \otimes_F \uM^{\otimes_A j}$ is the restriction of the algebra structure on the tensor
product of algebras $a \otimes \bbT_{A}(\uM)$.
The $k$-grading and differential on the complex $\bigoplus a^{djk} \otimes \uM^{\otimes_A j}$ are defined to be the total $k$-grading and total differential on the tensor product of complexes.
The $d$-grading is defined to be the total $d$-grading, with the degree $d$ part given by
$$\bigoplus_{d_0 + d_1 +...+ d_j = d} a^{d_0j \bullet} \otimes \uM^{d_1 \bullet} \otimes ... \otimes \uM^{d_j \bullet}.$$
The bimodule structure, grading and differential on $\bigoplus m^{djk} \otimes \uM^{\otimes_A j}$ are defined likewise.
We sometimes write
$$\mathbb{P}_{a,\um}(A,\uM) = (a(A,\uM), \um(A,\uM)).$$
If $a$ and $b$ are differential trigraded algebras, and ${}_a x_b$ is a differential trigraded $(a,b)$-bimodule, then we have a
differential bigraded $(a(A,\uM),b(A,\uM))$-bimodule
$$x(A,\uM) := \bigoplus_{d,j,k} x^{djk} \otimes \uM^{\otimes_A j}.$$

\begin{example}\label{c32}
Let $\bfc_3$ be the algebra given by quiver and relations as 
$$
\xymatrix{ \ar@/^/[r]^{\eta} \overset{1}{\bullet} 
&\overset{2}{\bullet}\ar@/^/[r]^{\eta} \ar@/^/[l]^{\xi}
&\ar@/^/[l]^{\xi} \overset{3}{\bullet} 
 },$$
modulo relations $\xi^2 = \eta^2 = \xi \eta + \eta \xi = 0 , \eta\xi e_3=0$. This is a quasi-hereditary Ringel self-dual algebra with respect to the natural order and  composition structure of the left projectives is given by the large numbers in 
{\footnotesize$$
\xymatrix@=8pt{1_0^0\ar@{-}[d] &&&2_0^0\ar@{-}[dr]\ar@{-}[dl]&&&3_0^0\ar@{-}[d] \\
2_1^1\ar@{-}[d] &&1_1^0\ar@{-}[dr]&&3_1^1\ar@{-}[dl]&&2_1^0\\
1_2^1&&& 2_2^1&&&
 }$$}
and there is the obvious $j$-grading by path length, which we have indicated as an index. The $k$-grading is identically zero, and the $d$-grading is taken to be the $\Delta$-grading in the quasi-hereditary structure, indicated by the superscripts.
Its tilting module $t$ as a left module looks like
{\footnotesize$$
\xymatrix@=8pt{&&1\ar@{-}[d] &&&2\ar@{-}[dr]\ar@{-}[dl] \\
1&&2\ar@{-}[d] &&1\ar@{-}[dr]&&3\ar@{-}[dl]\\
&&1&&&2&
 }$$}
and we have an isomorphism of $t \otimes t \cong \bfc_3^*$, which thanks to a simple-preserving duality looks like $\bfc_3$ upside down as a left module. In computing the algebra $\bfc_3(\bfc_3, (t, \Hom_{\bfc_3}(t,\bfc_3))$, the adjoint to the tilting module does not come into play, since $\bfc_3$ is positively graded.
The algebra $\bfc_3(\bfc_3, (t, \Hom_{\bfc_3}(t,\bfc_3))$ has 9 simple modules $(i,j)$ where $1 \leq i,j\leq 3$. We compute the projective $(1,2)$. The  algorithm prescribes that we should get a quotient of this by taking the simple 
$1_0^0$ and tensoring it with $t^{\otimes 0}e_2 = \bfc_3e_2$, a subquotient by taking  the simple $2_1^1$ and tensoring it with $te_2$, and a submodule by tensoring $1_2^1$ with $\bfc_3^*e_2$, giving the three modules
 {\footnotesize
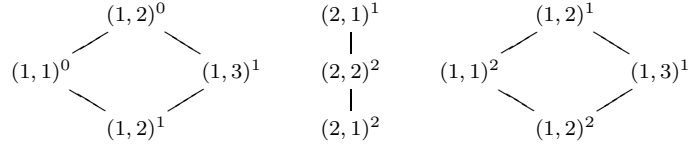
$$
\xymatrix@=8pt{&(1,2)^0\ar@{-}[dr]\ar@{-}[dl] \\
(1,1)^0\ar@{-}[dr]&&(1,3)^1\ar@{-}[dl]\\
&(1,2)^1&
 } \qquad
\xymatrix@=8pt{(2,1)^1\ar@{-}[d]  \\
(2,2)^2\ar@{-}[d]      \\
(2,1)^2
 }\qquad
\xymatrix@=8pt{&(1,2)^1\ar@{-}[dr]\ar@{-}[dl]& \\
(1,1)^2\ar@{-}[dr]&&(1,3)^1\ar@{-}[dl]\\
&(1,2)^2&
 }$$}
respectively. Identifying $(i,j)$ with $3(i-1)+j$ in a $3$-adic expansion, we obtain the second projective
{\footnotesize
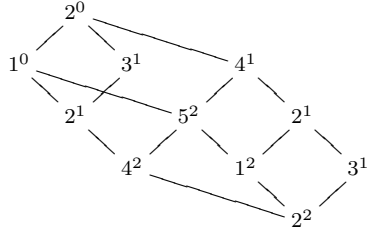
$$\xymatrix@=8pt{&2^0\ar@{-}[dr]\ar@{-}[drrr]\ar@{-}[dl] &&&&&\\
1^0\ar@{-}[drrr]\ar@{-}[dr]&&3^1\ar@{-}[dl]&&4^1\ar@{-}[dl]\ar@{-}[dr]&&\\
&2^1\ar@{-}[dr]&&5^2\ar@{-}[dl]\ar@{-}[dr]&&2^1\ar@{-}[dr]\ar@{-}[dl] &\\
&&4^2\ar@{-}[drrr]&&1^2\ar@{-}[dr]&&3^1\ar@{-}[dl]\\
&&&&&2^2&
}$$}
where the precise extensions require a more careful analysis (done in \cite[Section 1]{MT2}) and we have again indicated the $d$-degree by superscripts. The algebra $\bfc_3(\bfc_3, (t, \Hom_{\bfc_3}(t,\bfc_3))$ describes a block of polynomial representations of $GL_2$ in characteristic $3$ with $p^2=9$ simple modules, of which we have given the extension algebra of standard modules in Section \ref{example}. (Note that Ringel duality produces a twist on $t$, so that for computing the projective at $(1,1)$, we need the third summand of the tilting module and to compute the projective at $(1,3)$, we need the simple summand of $t$. Details on this can be found in \cite{MT2}.)
\end{example}

The operators $\bbP$ coincide with the operators $\bbO$ from our previous papers \cite{MT2,MT3} 
if we restrict them to positively graded objects of $\mathcal{U}$, concentrated in $d$-degree zero.
As the $d$-grading is just an extra grading dragged around, which does not interfere with the constructions,
the only place we need to take care in extending results from \cite{MT3} is the extension to the bonded setting.

However, using Lemma \ref{bondedtensor}, we find all proofs go through without problem. In particular, the operators $\bbP$ define $2$-endofunctors of $\mathcal{U}$ (see \cite[Lemma 9]{MT2} for the original $\bbO$), but we will not use this here.

The following lemmas are the bonded analogues of results we have established previously \cite{MT3} and are proved in exactly the same fashion.

\begin{lem} \label{tensor} \cite[Lemma 14]{MT3}
Let $a$ be a differential bigraded algebra, $x_a$ and ${}_ay$
differential bigraded modules, $(A, \uM)$ an object of $\mathcal{U}$.
Then
$$x(A,\uM) \otimes_{a(A,\uM)} y(A,\uM) \cong (x \otimes_a y)(A,\uM).$$
\end{lem}

\begin{lem} \label{homo}\cite[Lemma 15]{MT3}
Let $c$ be a differential trigraded algebra, and $(A,\uM)$ a Rickard
object of $\mathcal{U}$. If $x$ and $y$ are differential trigraded
$c$-modules, then we have a quasi-isomorphism
$$(\Hom_{c}(x,y))(A,\uM) \rightarrow \Hom_{c(A,\uM)}(x(A,\uM), y(A,\uM)).$$
\end{lem}

\begin{lem}\label{equivin}\cite[Lemma 22]{MT3}
Let $(A,\uM)$ and $(B,\uN)$ be objects of $\mathcal{U}$ such that $(A,\uM)
\gtrdot (B,\uN)$. Let $(a,\um)$ be a $j$-graded Rickard object of
$\mathcal{U}$. Then $\mathbb{P}_{a,\um}(A,\uM) \gtrdot
\mathbb{P}_{a,\um}(B,\uN)$.
\end{lem}

\begin{lem} \label{equivinagain}\cite[Lemma 23]{MT3}
 Let $(A,\uM)$ and $(B,\uN)$ be quasi-isomorphic objects of $\mathcal{U}$.
Let $(a,\um)$ be a $j$-graded Rickard object of $\mathcal{U}$. Then
$\mathbb{P}_{a,\um}(A,\uM)$ and $\mathbb{P}_{a,\um}(B,\uN)$ are
quasi-isomorphic objects of $\mathcal{U}$.
\end{lem}

\subsection{The operator $\fO$.}\label{SfO} We now recall the definition of the operator $\fO$
as well as the result which asserts that it behaves favourably with respect to taking homology \cite{MT3}.

Let $\Gamma = \bigoplus \Gamma^{ijk}$ be a
differential trigraded algebra. We have an operator
$$\fO_\Gamma \circlearrowright \{ \Sigma | \mbox{ $\Sigma = \bigoplus \Sigma^{jk}$ a differential bigraded algebra } \}$$
given by
$$\fO_\Gamma(\Sigma)^{ik} = \bigoplus_{j, k_1+k_2=k} \Gamma^{ijk_1} \otimes \Sigma^{jk_2}.$$
The algebra structure and differential are obtained by restricting
the algebra structure and differential from $\Gamma \otimes \Sigma$.
If we forget the differential and the $k$-grading, the operator
$\fO_\Gamma$ is identical to the operator $\fO_\Gamma$ defined in the
introduction.

\begin{lem} \label{commute} \cite[Lemma 18]{MT3} We have $$\mathbb{H} \fO_{\Gamma} \cong \mathbb{H} \fO_{\mathbb{H} \Gamma} \cong \fO_{\mathbb{H} \Gamma}
\mathbb{H},$$ for a differential trigraded algebra $\Gamma$.
\end{lem}

\subsection{Comparing $\mathbb{P}$ and $\fO$.}\label{comparison}
We have previously made a comparison of operators $\bbO$ and $\fO$ \cite{MT3}.
We could add a $d$-grading to $\fO$ to obtain a direct generalisation of this comparison result.
However, in our application,
we use the operator $\fO$ only in a setting where the $d$-grading is the negative of the $k$-grading,
and so to simplify we disregard it.
Let $\mathbb{D}$ denote the $2$-functor from $\mathcal{U}$ to $\mathcal{T}$ which disregards the $d$-grading.
Then, extending the operator $\bbO$ to the bonded setting, we have $\mathbb{D}\bbP = \bbO$.
Extending our comparison result in a bonded setting, again using Lemma \ref{bondedtensor},
we obtain the following:

\begin{lem}\label{compare} \cite[Corollary 21]{MT3}
Let $a$ be a differential bigraded algebra and $\um$ a bonded pair of $a$-$a$ dg bimodules. Then we have an
isomorphism of dg algebras
$$\mathbb{D}\bbP_{F,0} \bbP^n_{a,\um}(F,(F,F)) \cong \fO_F\fO_{\bbT_a(\um)}^n(F[z,z^{-1}]).$$
\end{lem}

\subsection{Keller's homological duality.} \label{kellerkoszul}\label{bbF}
Let $A$ be (quasi-isomorphic to) a finite dimensional algebra
with modules $S_1, \dots, S_f$ which generate the derived category
$D(A)$ of $A$. Let $P_l = \bigoplus_k P^k_l$ be a projective
resolution of $S_l$. Let $\cF (A)$ denote the dg algebra
$$\cF(A) = \bigoplus_{k,k'}\Hom_{A}(\bigoplus_{l=1}^fP^k_l,\bigoplus_{l=1}^f P^{k'}_l
).$$ Then $P = \bigoplus_{l=1}^f P_l$ is a differential graded
$A$-$\cF(A)$-bimodule. There are mutually inverse equivalences
$$\xymatrix{ D_{dg}(A) \ar@/^/[r]^{\Hom_A(P,-)} & \ar@/^/[l]^{P \otimes_{\cF (A)}-} D_{dg}(\cF (A))},$$
by a theorem of Keller \cite[Theorem 3.10]{Ke}. Since $P$ is
projective as an $A$-module, we have a natural isomorphism of
functors
$$\Hom_A(P,-) \cong \Hom_A(P,A) \otimes_A -.$$

 Let $(A,\uM)$ be an object of $\mathcal{U}$, where $A^d = 0$ for $d<0$.
We now define a homological duality operator $\bbF$ which sends $(A,\uM)$ to another object of $\mathcal{U}$.

We define $S_1,...,S_f$ to be the direct summands of the degree zero part of $A$ taken with respect to the $d$-grading.
Assume $P$ is $djk$-graded.
Given a differential graded $A$-$A$-bimodule $M$, denote by $\cF(M)$ the dg $\cF(A)$-$\cF(A)$-bimodule
$$\cF(M) = \Hom_A(P,A) \otimes_A M \otimes_A P.$$
We have $\cF(M) \cong \Hom_A(P,M) \otimes_A P$.
For a bonded pair $\uM = (M,M')$ of bimodules, denote by $\cF(\uM)$ the pair $(\cF(M),\cF(M') )$.
Let $\bbF(A,\uM):=(\cF(A), \cF(\uM))$, with bonded structure defined in the natural way as follows.
We have maps 

\begin{equation*} \begin{split}
       \cF(M) &\otimes_{\cF(A)}\cF(M') \\&=        \Hom_A(P,A) \otimes_A M \otimes_A P \otimes_{\cF(A)}\Hom_A(P,A) \otimes_A M' \otimes_A P \to \cF(A) 
\end{split}
\end{equation*}
 and
\begin{equation*}
\begin{split}
 \cF(M')& \otimes_{\cF(A)}\cF(M) \\&=\Hom_A(P,A) \otimes_A M' \otimes_A P \otimes_{\cF(A)}\Hom_A(P,A) \otimes_A M \otimes_A P \to \cF(A)
\end{split}
\end{equation*}
which are given by the composition of the natural map $P \otimes_{\cF(A)}\Hom_A(P,A) \to A$ and the maps $M \otimes_AM' \to A$ and $M' \otimes_AM \to A$ respectively, the latter maps being given by the bonding.
It follows immediately that these maps define a bonding, as the relevant morphisms
\begin{equation*}
\begin{split}
\Hom_A(P,A) \otimes_A M &\otimes_A P \otimes_{\cF(A)}\Hom_A(P,A) \otimes_A M' \otimes_A P  \otimes_{\cF(A)}
\Hom_A(P,A) \otimes_A M \otimes_A P \\&\rightarrow \Hom_A(P,A) \otimes_A M \otimes_A P\end{split}
\end{equation*}
factor through the morphisms
$$\Hom_A(P,A) \otimes_A M \otimes_A M' \otimes_A M \otimes_A P \rightarrow \Hom_A(P,A) \otimes M \otimes P,$$
where we can apply the bonding of $(M,M')$.

\begin{lem} \label{extsim} Suppose
$(A,\uM)$ is a Rickard object of $\mathcal{U}$, where $A^d = 0$ for $d<0$.
The bimodule $P$ induces a dg equivalence $(A,\uM) \gtrdot \bbF(A,\uM)$.
\end{lem}

\proof
The quasi-isomorphisms $$P \otimes_{\cF(A)}
\cF(M) \rightarrow M \otimes_{A} P$$ and $$P \otimes_{\cF(A)}
\cF(M) \rightarrow M' \otimes_{A} P$$ follow in the same fashion as in \cite[Lemma 9]{MT3}.
It remains to be checked that the diagrams
$$\xymatrix{
P \otimes_{\cF(A)} {\cF(M)} \otimes_{\cF(A)} {\cF(M')} \ar[r] \ar[d] & M \otimes_A \otimes M' \otimes_A P \ar[r] & A \otimes_A P \ar[d] \\
P \otimes_{\cF(A)} {\cF(A)} \ar[rr] & & P
}$$
$$\xymatrix{
P \otimes_{\cF(A)} {\cF(M')} \otimes_{\cF(A)} {\cF(M)} \ar[r] \ar[d] & M' \otimes_A \otimes M \otimes_A P \ar[r] & A \otimes_A P \ar[d] \\
P \otimes_{\cF(A)} {\cF(A)} \ar[rr] & & P
}$$
commute, which is straightforward.
\endproof

\begin{remark}\label{dagger}
Suppose that $(a,\um) = (a,(m, m^{-1}))$ is a $j$-graded Rickard object of $\mathcal{U}$, such that $a$ is concentrated in non-negative $d$-degrees and that a projective resolution $P$ of the degree zero part of $a$ is differential $djk$-trigraded.
Then $\cF(a)$ inherits a differential $djk$-trigrading.
In this case we call $(a,\um)$ a \emph{dagger object} of $\mathcal{U}$, and the dg equivalence in Lemma  \ref{extsim} is differential $djk$-trigraded.
\end{remark}

\subsection{Homological duality for operators $\bbP$.}
\label{Koszulops}

We now consider how the Keller duality operator $\bbF$ behaves with respect to the operators $\bbP$.

\begin{lem}\label{koszul}
Let $(A,\uM)$ be a Rickard object of $\mathcal{U}$.
Let $(a,\um)$ be a dagger object of $\mathcal{U}$, with homological dual $\bbF(a,\um)$. Assume further that $P(A,\uM)$ generates $D_{dg}(a(A,\uM))$.
Then we have a differential $dk$-bigraded equivalence $$\bbP_{a,\um}(A,\uM) \gtrdot \bbP_{\bbF(a,\um)}(A,\uM).$$
\end{lem}

\proof
Let $P$ be the differential $djk$-trigraded $a$-module inducing the duality between $a$ and $\cF(a)$ as in Subsection \ref{bbF}.
Notice that $P(A,\uM)$ is a differential $dk$-bigraded $a(A,\uM)$-$\cF(a)(A,\uM)$-bimodule.

Now $\cF(a)(A,\uM)$ is by definition $\Hom_a(P,P)(A,\uM)$, which by Lemma \ref{homo} is quasi-isomorphic to $\Hom_{a(A,\uM)}(P(A,\uM),P(A,\uM))$.

As $P(A,\uM)$ generates $D_{dg}(a(A,\uM))$ by assumption, by Keller's theory $P(A,\uM)$ induces a dg-equivalence between $D_{dg}(a(A,\uM)) $ and
$D_{dg}(\cF(a)(A,\uM)) $ \cite[Theorem 3.10]{Ke}.

To check the conditions of Definition \ref{dgequiv} on bimodules is an easy repeated application of Lemma \ref{tensor} to the diagrams for the equivalence in the proof of Lemma \ref{extsim}.
\endproof

\subsection{A quasi-isomorphism of operators.} \label{qiops}

Here we show that $\bbF \bbP_{a,\um}$ is quasi-isomorphic to
$\bbP_{\bbF(a,\um)} \bbF$ under suitable conditions. We first note that we always have a  $dk$-graded equivalence between the
objects $\bbF(\bbP_{a,\um}(A,\uM))$ and $\bbP_{\bbF(a,\um)}(\bbF(A,\uM))$ of $\mathcal{U}$.

\begin{thm}\label{chain} Let $(A,\uM)$ be a Rickard object of $\mathcal{U}$.
Let $(a,m)$ be a be a dagger object of $\mathcal{U}$. Assume further that $P(A,\uM)$ generates $D_{dg}(a(A,\uM))$. We have a
chain of differential $dk$-bigraded equivalences
$$\bbF(\bbP_{a,\um}(A,\uM))  \lessdot  \bbP_{a,\um} (A,\uM) \gtrdot \bbP_{\bbF(a,\um)} (A,M) \gtrdot \bbP_{\bbF(a,\um)}(\bbF(A,\uM)). $$
\end{thm}

\proof This follows from Lemmas \ref{extsim}, \ref{koszul} and
\ref{equivin}.
\endproof

Under special conditions, which in Proposition \ref{conds} will turn out to hold for our application to Weyl modules for $GL_2$, we can strengthen
this as follows:

\begin{thm}\label{chainquasi} Assume the following:
\begin{enumerate}[(i)]
\item\label{cond2} Let $(a,\um)$ be a dagger object of $\mathcal{U}$ such that $a$ is concentrated in non-negative $j$-degrees.
\item\label{cond1} Let $(A,\uM) = (A,(M,M^{-1}))$ be a Rickard object of $\mathcal{U}$, such that
\begin{enumerate}[(a)]
\item\label{cond11} both $A$ and $M$ are concentrated in non-negative $d$-degrees;
\item\label{cond12} $M^{\otimes j} \otimes_A A^{0 \bullet} \cong (M^{\otimes j})^{0 \bullet} $ for all $j$ such that $a^{0 j \bullet}\neq 0$.
\end{enumerate}

\item\label{cond3} The differential $dk$-bigraded $a(A,\uM)$-module $(a(A,\uM))^{0 \bullet}$ generates the derived category $D_{dg}(a(A,\uM))$.
\end{enumerate}
Then the chain of equivalences in Theorem \ref{chain} lifts to a %$dk$-graded
quasi-isomorphism from $\bbP_{\bbF(a,\um)}(\bbF(A,\uM))$ to $\bbF(\bbP_{a,\um}(A,\uM))$.
\end{thm}

\proof
We denote by $P_a$ the projective resolution of the $d$-degree zero part of $a$ and by $P_A$ the projective resolution of the $d$-degree zero part of $A$ and
note that we have an isomorphism
$$\cF(M)^{\otimes_{\cF(A)} r} \cong \Hom_A(P_A, M^{\otimes_A r} \otimes_A P_A)$$

and hence an isomorphism

\begin{equation}\label{weird}
 \cF(a)(\bbF(A, \uM)) \cong \Hom_A(P_A, \cF(a)(A,\uM) \otimes_A P_A).
\end{equation}

We have a quasi-isomorphism $P_a(A,\uM) \to a^{0 \bullet \diamond}(A,\uM)$ by \cite[Lemma 15]{MT2} and quasi-isomorphisms

\begin{equation*}\begin{split}M^{\otimes j} \otimes P_A &\to M^{\otimes j} \otimes_A A^{0 \clubsuit}\\
& \cong( M^{\otimes j})^{0 \clubsuit}
\end{split}
\end{equation*}

for all j such that $a^{0 j \bullet}\neq 0$ by \eqref{cond12}.
Thanks to \eqref{cond2} and \eqref{cond11} we also have
\begin{equation*}
\begin{split}
 (a(A,\uM))^{0 \bullet} & \cong \bigoplus_{j \geq 0, d \geq 0} a^{d j\diamond} \otimes_F(M^{\otimes j})^{-d \clubsuit}
\\&\cong \bigoplus a^{0 j\diamond} \otimes_F(M^{\otimes j})^{0 \clubsuit}
\end{split}
\end{equation*}
 and putting the latter two observations together, we obtain that $P_a(A,\uM) \otimes_A P_A$ is a projective resolution of $(a(A,\uM))^{0 \bullet}$.

Therefore \begin{equation}\label{resol}\cF(a(A,\uM)) \cong \End_{a(A,\uM)}(P_a(A,\uM) \otimes_A P_A).\end{equation}

To prove the theorem we need to show that $\cF(a(A,\uM))$ is quasi-isomorphic to $\cF(a)(\bbF(A, \uM))$. We have

\begin{equation*}
\begin{split}
\cF(a)(\bbF(A, \uM)) & = \Hom_A(P_A, \cF(a)(A,\uM) \otimes_A P_A)\\& \qquad\mbox{by \eqref{weird}}\\
& = \Hom_A(P_A, \Hom_a(P_a,P_a)(A,\uM) \otimes_A P_A)\\ &\qquad \mbox{by definition of $\cF(a)$}\\
& \cong \Hom_A(P_A, \Hom_{a(A,\uM)}(P_a(A,\uM),P_a(A,\uM))\otimes_A P_A )\\ &\qquad \mbox{by Lemma \ref{homo}}\\
& \cong \Hom_A(P_A, \Hom_{a(A,\uM)}(P_a(A,\uM),P_a(A,\uM)\otimes_A P_A ))\\ &\qquad \mbox{by projectivity ot $P_A$}\\
& \cong \Hom_{a(A,\uM)}(P_a(A,\uM)\otimes_A P_A,P_a(A,\uM)\otimes_A P_A) \\ &\qquad \mbox{by adjunction}\\
& \cong  \cF(a(A,\uM)) \\ &\qquad \mbox{by \eqref{resol}}
\end{split}
\end{equation*}

Similarly we have

\begin{equation*}
\begin{split}
& \cF(m)(\bbF(A, \uM)) \\
& = \Hom_A(P_A, \cF(m)(A,\uM) \otimes_A P_A)\\
& = \Hom_A(P_A, \Hom_a(P_a,m\otimes_a P_a)(A,\uM) \otimes_A P_A)\\
& \cong \Hom_A(P_A, \Hom_{a(A,\uM)}(P_a(A,\uM),m(A,\uM)\otimes_{a(A,\uM)}P_a(A,\uM))\otimes_A P_A )\\
& \cong \Hom_A(P_A, \Hom_{a(A,\uM)}(P_a(A,\uM),m(A,\uM)\otimes_{a(A,\uM)}P_a(A,\uM)\otimes_A P_A ))\\
& \cong \Hom_{a(A,\uM)}(P_a(A,\uM)\otimes_A P_A,m(A,\uM)\otimes_{a(A,\uM)}P_a(A,\uM)\otimes_A P_A) \\
& \cong  \cF(m(A,\uM))
\end{split}
\end{equation*}
and the analogous chain for $m^{-1}$.

We need to check compatibility of these quasi-isomorphisms with the bonding.

The bonding on $\cF(m)(\bbF(A,\uM))$ and $\cF(m^{-1})(\bbF(A,\uM))$ is given via
$$\cF(m)(\bbF(A,\uM))\otimes_{\cF(a)(\bbF(A,\uM))}\cF(m^{-1})(\bbF(A,\uM)) \cong (\cF(m)\otimes_{\cF(a)}\cF(m^{-1}))(\bbF(A,\uM))$$
from Lemma \ref{tensor} and the bonding $\cF(m)\otimes_{\cF(a)}\cF(m^{-1}) \to \cF(a)$ which is induced by the evaluation map  $P_a \otimes_{\cF(a)} \Hom_a(P_a, a) \overset{\sim}{\rightarrow} a$  and the bonding on $m\otimes_a m^{-1}$.

For the bonding on $\cF(m(A,\uM))$ and $\cF(m^{-1}(A,\uM))$ note that, since $(a(A,\uM))^{0 \bullet}$ generates $D_{dg}(a(A,\uM))$, its projective resolution $P:=P_a(A,\uM)\otimes_A P_A$ contains a progenerator of $a(A,\uM)$ so, the evaluation map $$P \otimes_{\cF(a(A,\uM))}\Hom_{a(A,\uM)}(P, a(A,\uM)) \to a(A,\uM)$$ is an isomorphism and, using projectivity of $P$ and, in the last step, Lemma \ref{homo} we obtain

\begin{equation*}
\begin{split}
\cF& (m(A,\uM)) \otimes_{\cF(a(A,\uM))}  \cF(m^{-1}(A,\uM)) \\
 = &\Hom_{a(A,\uM)}(P,m(A,\uM)\otimes_{a(A,\uM)}P) \\&\otimes_{\cF(a(A,\uM))}   \Hom_{a(A,\uM)}(P,m^{-1}(A,\uM)\otimes_{a(A,\uM)}P)\\
\cong & \Hom_{a(A,\uM)}(P, a(A,\uM)) \otimes_{a(A,\uM)}m(A,\uM) \\ &\otimes_{a(A,\uM)}P \otimes_{\cF(a(A,\uM))}  \Hom_{a(A,\uM)}(P,a(A,\uM)) \\ &\otimes_{a(A,\uM)}m^{-1}(A,\uM)\otimes_{a(A,\uM)}P\\
 \cong &\Hom_{a(A,\uM)}(P, a(A,\uM)) \otimes_{a(A,\uM)}m(A,\uM)\\&\otimes_{a(A,\uM)}m^{-1}(A,\uM)\otimes_{a(A,\uM)}P\\
\cong & \Hom_{a(A,\uM)}(P, a(A,\uM)) \otimes_{a(A,\uM)} (m\otimes_a m^{-1})(A,\uM) \otimes_{a(A,\uM)}P
\end{split}
\end{equation*}

which together with the bonding on $m\otimes_a m^{-1}$ induce the bonding. The analogous statements hold for reversed roles of $m$ and $m^{-1}$ and our quasi-isomorphisms are hence compatible with the bonding, completing the proof of the theorem.
\endproof

\section{Recollections on $GL_2$.}\label{zigzag}\label{GL2}

We say a $\mathbb{Z}$-grading $A = \oplus_{d \geq 0} A^d$ on a quasi-hereditary algebra is a $\Delta$-grading if $A^d$ is isomorphic to a direct sum of standard modules, for all $d$.

Let $Z$ denote the algebra given by the quiver
$$
\xymatrix{ \cdots & \overset{0}{\bullet} \ar@/^/[r]^{\eta}
& \ar@/^/[r]^{\eta} \overset{1}{\bullet} \ar@/^/[l]^{\xi}
&\overset{2}{\bullet}\ar@/^/[r]^{\eta} \ar@/^/[l]^{\xi}
&\ar@/^/[l]^{\xi} \overset{3}{\bullet} &\cdots
 },$$
modulo relations $\xi^2 = \eta^2 = \xi \eta - \eta \xi = 0$. We will call this the \emph{zigzag algebra} and it or its various truncations show up in many guises in representation theory. It is therefore a very well-studied little infinite dimensional algebra whose projective indecomposable modules have a Loewy structure given by {\scriptsize 
$$\xymatrix@=5pt{&l\ar@{-}[dl]\ar@{-}[dr]&\\l-1\ar@{-}[dr]&&l+1\ar@{-}[dl]\\&l&}
$$}
where $l$ denotes the simple module at vertex $l$, 
and is well-known to have
a number of interesting homological properties, all of which are easily checked by hand. For example, it is
quasi-hereditary, symmetric, and Koszul.
Furthermore, its derived category admits an action of the braid $2$-category, in which braids act faithfully \cite[Section 2]{KS}.
We are especially interested in the quasi-hereditary structure on $Z$,
in which $\Delta(l)$ has top $l$ and socle $l-1$ for $l \in \mathbb{Z}$.

For our application to $GL_2$, we will be interested in a finite
truncation of $Z$. Let $\bfc$ be the finite dimensional subquotient of
$Z$ generated by $$\xymatrix{ \overset{1}{\bullet}
\ar@/^/[r]^{\eta_1} &\overset{2}{\bullet}\ar@/^/[l]^{\xi_1}
\ar@/^/[r]^{\eta_2} &\ar@/^/[l]^{\xi_2} \overset{3}{\bullet} &\cdots
&\overset{p-1}{\bullet}\ar@/^/[r]^{\eta_{p-1}} &
\ar@/^/[l]^{\xi_{p-1}} \overset{p}{\bullet},}$$ modulo the ideal
$I=(\xi_{l+1} \xi_{l}, \eta_{l} \eta_{l+1}, \xi_l \eta_l +
\eta_{l+1} \xi_{l+1} , \eta_{p-1} \xi_{p-1}  \mid 1 \leq l \leq
p-2)$.

We will now recall some facts about the rational representation theory
of $G = GL_2(F)$. The category of polynomial representations of $G$
of degree $r$ is equivalent to the category $S(2,r) \ml$ of
representations of the Schur algebra $S(2,r)$ \cite{Green}. %All blocks of $S(2,r) \ml$ whose number of isomorphism classes of simple objects is $p^q$ are equivalent and 
%A block with $p$ simple module is described by the algebra $\bfc$ (\cite[Proposition 4.1]{Erd}). 
A block is Ringel self-dual if and only if it has $p^q$ simple modules \cite{EH}. In, \cite{MT1}.\cite{MT2} we developed a combinatorial way to describe these blocks, which we now describe.

We denote by $\sigma \nu$ the algebra involution of $Z$ which sends vertex
$i$ to vertex $p-i$ and exchanges $\xi$ and $\eta$. Let $e_l$ denote
the idempotent of $Z$ corresponding to vertex $l \in \mathbb{Z}$.
Let
$$\bft = \sum_{1 \leq l \leq p, 0 \leq m \leq p-1} e_l Z e_m.$$ Then
$\bft$ admits a natural left action by the subquotient $\bfc$ of $Z$. By
symmetry, $\bft$ admits a right action by $\bfc$, if we twist the regular
right action by $\sigma \nu$. In this way, $\bft$ is naturally a
$\bfc$-$\bfc$-bimodule (see \cite[Lemma 14]{MT1} in a more general case, this case is easily checked by hand).

The algebra $\bfc$ is a quasi-hereditary algebra, and the left
restriction $_\bfc \bft$ of $\bft$ is a full tilting module for $\bfc$. The
natural homomorphism $\bfc \rightarrow \Hom_\bfc({}_\bfc \bft,{}_\bfc \bft)$ defined by the
right action of $\bfc$ on $\bft$ is an isomorphism, implying that $\bfc$ is indeed
Ringel self-dual (again \cite[Theorem 19]{MT1} for a more general case, this case is easily checked by hand).

Let $\tilde{\bft}$ denote a projective
resolution of $\bft$ as a $\bfc$-$\bfc$ bimodule, then $\tilde{\bft}$ is a
two-sided tilting complex, and $\tilde{\bft} \otimes_\bfc -$ induces a
self-equivalence of the derived category $D^b(\bfc)$ of $\bfc$. We also have the adjoint complex $\tilde{\bft}^{-1} = \Hom_\bfc(\tilde{\bft}, \bfc)$.
We denote by $\ut$ the pair of bonded $\bfc$-$\bfc$ dg bimodules $(\tilde{\bft},\tilde{\bft}^{-1})$.

The operator $\mathbb{P}_{\bfc,\ut}$ (see Section \ref{bbP}) acts on the collection of algebras
with a bonded pair of bimodules, such as the pair $(F, \uF)$ whose algebra is $F$ and
whose bonded bimodules $\uF$ are just the pair of regular bimodules $(F,F)$. The operator
$\mathbb{P}_{F,0} \mathbb{P}_{\bfc,\ut}^q$ takes an algebra with a bonded pair of bimodules to an algebra, along with a pair of zero bimodules which we disregard.  We define $\mathbf{b}_q$ to be the category of modules over the
$\Delta$-graded algebra $\mathbb{P}_{F,0} \mathbb{P}_{\bfc,\ut}^q(F,\uF)$.

We have an algebra homomorphism $\bfc \rightarrow F$ which sends a path
in the quiver to $1 \in F$ if it is the path of length zero based at
$1$, and $0 \in F$ otherwise. This algebra homomorphism lifts to a
morphism of operators $\mathbb{P}_{\bfc,\ut} \rightarrow
\mathbb{P}_{F,0}$. We have $\mathbb{P}_{F,0}^2 = \mathbb{P}_{F,0}$;
We thus have a natural sequence of operators
$$\mathbb{P}_{F,0} \leftarrow \mathbb{P}_{F,0} \mathbb{P}_{\bfc,\ut} \leftarrow \mathbb{P}_{F,0} \mathbb{P}_{\bfc,\ut}^2 \leftarrow ...,$$
which, if we apply each term to $(F,\uF)$ and take representations,
gives us a sequence of embeddings of highest weight categories
$$\mathbf{b_1} \rightarrow \mathbf{b_2} \rightarrow \mathbf{b_3} \rightarrow ...$$
We denote by $\mathbf{b}$ the union of these highest weight categories. In a
previous paper, we have proved the following:

\begin{thm} \cite[Corollary 21,Corollary 27]{MT2} \label{GL2theorem}
Every block of $G \ml$ is $\Delta$-equivalent to $\mathbf{b}$. Every block of
$S(2,r) \ml$ whose number of isomorphism classes of simple objects
is $p^q$ is $\Delta$-equivalent to $\mathbf{b_q}$.
\end{thm}

Note that this was formulated only in terms of the operators $\mathbb{O}_{\bfc, \bft}$ (\cite[Corollary 21]{MT2}) and $\mathbb{O}_{\bfc, \tilde\bft}$ (\cite[Lemma 22]{MT2}). However, both $\bfc$ and $\bft$ are concentrated in non-negative $j$-degrees, and therefore
$\tilde{\bft}^{-1}$ is only present for formal purposes: it is irrelevant for computation and $\mathbb{P}_{\bfc,\ut}=\mathbb{O}_{\bfc, \tilde\bft}$.
The fact that the operator $\mathbb{P}_{\bfc,\ut}=\mathbb{O}_{\bfc, \tilde\bft}$ truly gives us the correct $\Delta$-grading
also follows from our previous work \cite[Theorem 18]{MT1}.

\section{The homological dual of $(\bfc, (\bft, \bft^{-1}))$.}\label{basestep}

In view of the previous section, if we take the $d$-grading to be that given by the $\Delta$-grading, our task is to compute the homology of the algebra component of $\bbF\mathbb{P}_{\bfc,\ut}^q(F,\uF)$. We would like to use Theorem \ref{chainquasi} to move the operator $\bbF$ to the right until we are down to the trivial task of applying  $\bbF$ to $(F,\uF)$. We therefore first check that the conditions given in this theorem are satisfied in our given situation.

\subsection{Hypotheses of Theorem \ref{chainquasi}.}\label{hypos}

We now gather the various facts that we need in order to ensure the hypotheses of Theorem \ref{chainquasi} are satisfied in our situation.

First note that, as the $d$-grading on $\bft$ does not feature in our algebraic constructions, we are free to choose it as it suits us.
We define $\bft^{0 \bullet \diamond}$ to be the natural quotient of $\bft$ which is the direct sum of costandard modules and $\bft^{1 \bullet \diamond}$ as the kernel of this surjection. This is obviously a grading on $\bft$ with respect to which $\bft$ is concentrated in degree $0$ and $1$.

\begin{prop}\label{check15} Set $(A_q,\uM_q)=\bbP_{\bfc,\ut}^q(F,F)$. Then $(A,\uM_q)$ is a Rickard object of $\mathcal{U}$ such that
\begin{enumerate}[(i)]
\item\label{22i} Both $A_q$ and $M_q$  are concentrated in non-negative $d$-degrees.
\item\label{22ii} We have an isomorphism $M_q \otimes_{A_q} A_q^{0 \bullet} \cong M_q^{0 \bullet}$.
\item\label{22iii} The differential $dk$-bigraded $\bfc(A_q,\uM_q)$-module $(\bfc(A_q,\uM_q))^{0 \bullet}$ generates
the derived category $D_{dg}(\bfc(A_q,\uM_q))$.
\end{enumerate}
\end{prop}

\proof 
$(A_q,\uM_q)$ is a Rickard object by \cite[Lemma 30]{MT3}.

Claim \eqref{22i} for $A_q$ follows from the fact that the $d$-grading on $\bbP_{F,0} \bbP_{\bbF(\bfc,\ut)}^q \bbF(F,F)$ is just given by the $\Delta$-grading in the quasi-hereditary structure. As $\bft$ is concentrated in $d$-degrees $0$ and $1$, $\bfc$ is also non-negatively $d$-graded, a projective $\bfc$-$\bfc$-bimodule resolution $\tilde{\bft}$ will also be concentrated in non-negative $d$-degrees. Then Claim \eqref{22i} for $M_q$ follows from this fact that $\bfc$ and $\tilde \bft$ are concentrated in positive $d$-degrees (as the base step in an induction) and the fact that $\bft$ is concentrated in positive $j$-degrees, to ensure in the inductive step that only $\bfc$ and $\bft$ are used in the iterative construction of $M_q$.

Claim \eqref{22ii} is also checked by induction on $q$.
 As $\bft \otimes_\bfc \Delta \cong \nabla$, we automatically obtain a quasi-isomorphism $\tilde{\bft} \otimes_\bfc \bfc^{0 \bullet \diamond} \cong \bft \otimes_\bfc \bfc^{0 \bullet \diamond} \cong \bft^{0 \bullet \diamond}$, which is the base step for Claim \eqref{22ii}.
For the inductive step we compute that
\begin{equation*}\begin{split}
    M_q \otimes_{A_q} A_q^{0 \bullet} & \cong \bft(A_{q-1}, \uM_{q-1} )    \otimes_{\bfc(A_{q-1}, \uM_{q-1})} (\bfc (A_{q-1}, \uM_{q-1}))^{0 \bullet}\\
&\cong    \bft(A_{q-1}, \uM_{q-1} )    \otimes_{\bfc(A_{q-1}, \uM_{q-1})} (\bigoplus_{j \geq 0, d \geq 0}\bfc^{dj \bullet}\otimes_F (M_{q-1}^{\otimes j})^{-d \bullet})  \\ &\qquad \mbox{as $c$ is concentrated in non-negative $d$- and $j$-degrees}\\
& \cong \bft(A_{q-1}, \uM_{q-1} )  \otimes_{\bfc(A_{q-1}, \uM_{q-1})} (\bigoplus_{j \geq 0}\bfc^{0j \bullet}\otimes_F (M_{q-1}^{\otimes j})^{0 \bullet})\\& \qquad \mbox{as $M_{q-1}$ is concentrated in non-negative $d$-degrees}\\
& \cong \bft(A_{q-1}, \uM_{q-1} )  \otimes_{\bfc(A_{q-1}, \uM_{q-1})} (\bigoplus_{j=0,1}\bfc^{0j \bullet}\otimes_F (M_{q-1}^{\otimes j})^{0 \bullet})\\& \qquad \mbox{as $c^{0 \diamond \bullet}$ is concentrated in $j$-degrees $0,1$}\\
&\cong \bft(A_{q-1}, \uM_{q-1} ) \otimes_{\bfc(A_{q-1}, \uM_{q-1})} \bfc^{0\diamond \bullet}(A_{q-1}, \uM_{q-1} ) \otimes_{A_{q-1}} A_{q-1}^{0 \clubsuit}\\ & \qquad \mbox{by the induction hypothesis}\\
& \cong (\bft \otimes_{\bfc} \bfc^{0\diamond \bullet})(A_{q-1}, \uM_{q-1} ) \otimes_{A_{q-1}} A_{q-1}^{0 \clubsuit} \\ & \qquad \mbox{by Lemma \ref{tensor}}\\
& \cong \bft^{0\diamond \bullet}(A_{q-1}, \uM_{q-1} ) \otimes_{A_{q-1}} A_{q-1}^{0 \clubsuit} \\
& \cong \bigoplus_{j = 0,1}\bft^{0j \bullet} \otimes_F\uM_{q-1}^{\otimes j}\otimes   A_{q-1}^{0 \clubsuit}\\& \qquad \mbox{as $t^{0\diamond \bullet}$ is concentrated in $j$-degrees 0,1}\\
& \cong \bigoplus_{j = 0,1}\bft^{0j \bullet} \otimes_F (\uM_{q-1}^{\otimes j})^{0 \clubsuit} \\& \qquad \mbox{by the induction hypothesis}\\
&\cong \bft(A_{q-1}, \uM_{q-1})^{0 \bullet} \\& \qquad \mbox{as both $\bft$ are concentrated in non-negative $d$-degrees}\\
& \cong M_q^{0 \bullet}.
\end{split}
\end{equation*}

For Claim \eqref{22iii} notice that by Theorem \ref{GL2theorem} $(\bfc (A_{q-1}, \uM_{q-1}))^{0 \bullet}$ is isomorphic to the sum over all the standard modules in the block $\mathbf{b_q}$ which is well known to be a generator of the derived category.
\endproof

In order to apply Theorem \ref{chainquasi}, we further need the following lemma.

\begin{lem}\label{cdagger}
The object $(\bfc, \ut) \in \mathcal{U}$ is a dagger object (cf. Remark \ref{dagger}) and $\bfc$ is concentrated in non-negative $d$- and $j$-degrees.
\end{lem}

\proof
We need to be careful about our gradings: $\bfc$ is $djk$-graded, with $\eta$ in degree $(1,1,0)$ and $\xi$ in degree $(0,1,0)$. Hence $\bfc$ is concentrated in non-negative $d$- and $j$-degrees.
In order to obtain $\cF(\bfc)$, we need to take a projective resolution of $\bfc^{0,\bullet, \diamond}$, which is given by $P=\bigoplus_{1 \leq l \leq p} P_l$ where $P_l$ is a linear projective resolution of $\Delta(l) = \bfc^{0,\bullet, \diamond}e_l$. A straightforward computation shows that this is given by
$$P_l = \bigoplus_{k=l}^p \bfc e_k \lceil k-l \rfloor \langle k-l \rangle [-(k-l)]$$
where $\lceil \cdot \rfloor$ denotes a shift in the $d$-grading.
As required, the differential, which is given by right multiplication by $\eta$, has $djk$-grading $(0,0,1)$, turning $P$ into a differential $djk$-trigraded left $\bfc$-module. This completes the proof of the lemma.
\endproof

We have now proved the following.

\begin{prop}\label{conds}
Setting $(a,\um)= (\bfc,\ut)$ and $(A,\uM)=\bbP_{\bfc,\ut}^q(F,F)$, the assumptions of Theorem \ref{chainquasi} are satisfied.
\end{prop}

\proof
This is a summary of Lemma \ref{cdagger} for Condition \eqref{cond2},
and Proposition \ref{check15} for Conditions \eqref{cond1} and \eqref{cond3}.
\endproof

We have now ensured that we can apply the theory developed in Section \ref{theory}.
Towards describing ${\mathbf w}$ combinatorially using our algebraic operators,
we next give a combinatorial description of $\bbF(\bfc,\ut)=(\cF(\bfc),\cF(\ut))$ (see Section \ref{bbF}).

\subsection{The algebra $\fc$.}\label{psi}

Note that in view of Remark \ref{dagger} and Lemma \ref{cdagger}, the dg-equivalence between $(\bfc, \ut)$ and $\bbF(\bfc, \ut)$ is $djk$-graded.

Using the notation introduced in the proof of Lemma \ref{dagger}, i.e.\ $P=\bigoplus_{1 \leq l \leq p} P_l$ being a projective resolution of the $d$-degree $0$ part of $\bfc$, given by $$P_l = \bigoplus_{k=l}^p \bfc e_k \lceil k-l \rfloor \langle k-l \rangle [-(k-l)],$$ we now describe the algebra $\cF(\bfc) = \End_\bfc(P)$.
Let $\fc$ be the algebra given by the quiver
$$\xymatrix{ \overset{1}{\bullet}  & \ar@/^/[l]^{\xi}
 \ar@/_/[l]_{x}
\overset{2}{\bullet} & \ar@/^/[l]^{\xi}
 \ar@/_/[l]_{x}
\overset{3}{\bullet} &  ... & \overset{p-1}{\bullet} \ar@/^/[l]^{\xi}
 \ar@/_/[l]_{x} &\ar@/^/[l]^{\xi}
 \ar@/_/[l]_{x}
 \overset{p}{\bullet}},$$ modulo relations $\xi x-x \xi = 0, \xi^2 = 0$, where $x$ is given $djk$-degree $(-1,-1,1)$ and $\xi$ is given $djk$-degree $(0,1,0)$.

The map $\fc \rightarrow \End_\bfc(P)$ given by

$$e_l \mapsto (id_{P_l}:P_l \to P_l)$$

\begin{equation*}\begin{split}
 e_{l-1}xe_l \mapsto\qquad \qquad \qquad  ( P_{l-1} &\to P_l ): \\ \alpha e_{l-1+k}\lceil k \rfloor \langle k\rangle [-k] &\mapsto \begin{cases}
0 & \text{if } k=0,\\
 \alpha e_{l-1+k}\lceil k-1 \rfloor \langle k-1\rangle [-(k+1)]& \text{if } k > 0.
\end{cases}
\end{split}
\end{equation*}

which has $djk$-degree $(-1,-1,1)$, and

$$e_{l-1}\xi e_l \mapsto (P_{l-1} \to P_l : \alpha e_{l-1+k}\lceil k \rfloor \langle k\rangle [-k]) \mapsto \alpha \xi e_{l+k}\lceil k \rfloor \langle k\rangle [-k]))$$

which has $djk$-degree $(0,1,0)$ is easily seen to give a right $\fc$-action on $P$,
turning $P$ into a differential $djk$-trigraded $\bfc$-$\fc$ bimodule.

The algebra $\fc$ is isomorphic to $\Ext^\bullet_\bfc(\Delta,\Delta)$ (cf. \cite[Example 5.1.1]{Mad}),
which is quasi-isomorphic to $\End_\bfc(P)$, implying we have a derived equivalence
$$\xymatrix{ D(\bfc \trigra_{djk}) \ar@/^/[r]^{\Hom_\bfc(P,-)} & \ar@/^/[l]^{P \otimes^L_{\fc}-} D(\fc \trigra_{djk})}.$$
where $D(\bfc \trigra_{djk})$ denotes the derived category of differential $djk$-trigraded modules.
Given that the $k$-grading is exactly the $\Ext$-grading on
$\Ext^\bullet_\bfc(\Delta,\Delta)$ and the $d$-grading on $\fc$ is the negative of the $k$-grading,
we will from now on ignore the $d$-grading on $\fc$.

\begin{remark} \label{koszulselfdualalgebra} The algebra $\fc$ is Koszul.
The algebra $\fc^!$ is the quadratic dual of $\fc$, and is generated by the quiver
$$
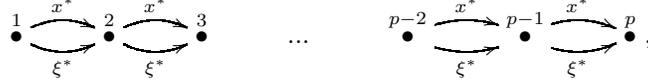
\xymatrix{ \overset{1}{\bullet}  \ar@/^/[r]^{x^*}
 \ar@/_/[r]_{\xi^*} &
\overset{2}{\bullet}  \ar@/^/[r]^{x^*}
 \ar@/_/[r]_{\xi^*} &
\overset{3}{\bullet} &  ... & \overset{p-2}{\bullet} \ar@/^/[r]^{x^*}
 \ar@/_/[r]_{\xi^*} & \overset{p-1}{\bullet} \ar@/^/[r]^{x^*}
 \ar@/_/[r]_{\xi^*} &
 \overset{p}{\bullet}
},$$ modulo relations $\xi^* x^*+ x^* \xi^* = 0, x^{*2} = 0$.
We have an algebra isomorphism $\fc \cong \fc^!$
which exchanges $e_s$ with $e_{p-s+1}$, exchanges $\xi^*$ with $x$, and exchanges $e_{s} x^* e_{s-1}$ with $(-1)^{s}e_{p-s+1} \xi e_{p-s+2}$.
\end{remark}

Our next task will be to trace $\ut$ through this duality, but first we need to consider a special case of homological duality.
Note that $\fc$ looks like a graded version of a truncation of the tensor product of a polynomial ring (generated by the arrows $x$) and its Koszul dual, an exterior algebra (generated by the arrows $\xi$). This will be made more precise in the proof of Proposition \ref{longlemma}, but motivates the analysis of the tensor product of an algebra and its Koszul dual in the following subsection.

\subsection{The tensor product of an algebra and its Koszul dual.}\label{general}

Throughout the rest of this section, we work in a monoidal $2$-category $\mathfrak{Dgalg}$ whose objects are differential $k$-graded algebras,
whose arrows $A \rightarrow B$ are objects in the derived category of dg $A$-$B$-bimodules with product given by derived tensor,
whose $2$-arrows are morphisms in the derived category of bimodules, and
whose monoidal structure is given by the super tensor product $\otimes$ over $F$.

Suppose $A$ is a Koszul algebra with Koszul dual $A^!$. For a general account of Koszul duality, please see Appendix 2 in Section \ref{kosz}.
Let $K_A = A \otimes_{A^0} A^{!*}$ denote the left Koszul complex for $A$ and $C_A = A \otimes_{A^0} A^!$ the adjoint of the Koszul complex for $A^!$.
Let $_AK = A^* \otimes_{A^0} A^!$ denote the associated right Koszul complex.
Given two $k$-graded algebras $A$ and $B$ we denote $\rho$ the algebra isomorphism $A \otimes B \rightarrow B \otimes A$ that sends
$a \otimes b$ to $(-1)^{|b||a|}b \otimes a$.

The tensor product $A \otimes A^! = A \otimes_F A^!$ is a Koszul self-dual algebra.
The following lemma establishes that Koszul self-duality for $A \otimes A^!$
is homologically dual to $\rho \circlearrowleft A^! \otimes A^!$ when $A^!$ is a symmetric algebra.

\begin{lem} Suppose $A^!$ is a symmetric algebra, thus $A \cong A^![s]\langle m \rangle$ for some $s$.
Then we have a commutative diagram
$$\xymatrix{
A \otimes A^! \ar[rr]^{C_A \otimes A^!} \ar[d]_{(K_A \otimes K_{A^!})^\rho}
& & A^! \otimes A^! \ar[d]^{\rho[s] \langle m \rangle} \\
A \otimes A^! \ar[rr]^{C_{A} \otimes A^!} & & A^! \otimes A^!  \\
}$$

in $\mathfrak{Dgalg}$.
\end{lem}
\proof
We have a diagram
$$\xymatrix{
A \otimes A^! \ar[rr]^{C_A \otimes A^!} \ar[d]_{K_A \otimes K_{A^!}} & & A^! \otimes A^! \ar[d]^{(_{A^!}K \otimes_A K_A) \otimes A^!} \\
A^! \otimes A \ar[rr]^{A^! \otimes C_A} \ar[d]_\rho & & A^! \otimes A^! \ar[d]^\rho \\
A \otimes A^! \ar[rr]^{C_A \otimes A^!} & & A^! \otimes A^!
}$$
which commutes because $C_A {_{A^!}K}$ is isomorphic to $A$ and $K_{A^!}C_{A}$ is isomorphic to $A^!$ (in the relevant derived categories of bimodules).
We have ${}_{A^!}K \otimes_{A} K_A \cong A^{!*} \otimes_{A^0} A \otimes_{A^0} A^{!*}$. This is isomorphic to
$A^{!*} \otimes_{A^!} C_{A^!} \otimes_A K_A$ which is isomorphic to $A^{!*}$ since $C_{A^!}$ and $K_A$
are adjoint equivalences. We have $A^{!*} \cong A[s]\langle m \rangle$ since $A$ is a symmetric algebra.
This implies our diagram is equivalent to
$$\xymatrix{
A \otimes A^! \ar[rr]^{C_A \otimes A^!} \ar[d]_{K_A \otimes K_{A^!}} & & A^! \otimes A^! \ar[d]^{[s]\langle m \rangle} \\
A^! \otimes A \ar[rr]^{A^! \otimes C_A} \ar[d]_\rho & & A^! \otimes A^! \ar[d]^\rho \\
A \otimes A^! \ar[rr]^{C_A \otimes A^!} & & A^! \otimes A^!
}$$
\endproof

\subsection{The homological dual of $\bft$.}\label{tthroughK}
The algebra $\fc$ is Koszul self-dual, and quasi-isomorphic to  $\cF(\bfc)$.
The differential $jk$-graded bimodule $\fc^\sigma \otimes_{\fc^0} \fc^*$
with differential given by internal multiplication by
$x \otimes \xi + \xi \otimes x$, i.e.
$$\ttd(a \otimes b) = (-1)^{|a|_k}(ax \otimes \xi b + a \xi \otimes xb),$$
is a twisted version of the Koszul complex for $\fc$.
Here we show this bimodule $\fc^\sigma \otimes_{\fc^0} \fc^* \langle 1 \rangle$ is dual to the $\bfc$-$\bfc$-bimodule $\bft$ under $\cF$, i.e. is quasi-isomorphic to $\cF(\bft)$.

\begin{prop} \label{longlemma}
We have a commutative diagram
$$\xymatrix{
\fc   & &  \ar[ll]^{P} \bfc  \\
\fc  \ar[u]^{\fc^\sigma \otimes_{\fc^0} \fc^*} & & \ar[ll]^{P} \bfc \ar[u]_{\bft \langle -1 \rangle} \\
}$$
in $\mathfrak{Dgalg}$.
\end{prop}
\proof
To begin with consider the Koszul algebra $A = S(x)$, with $x$ in $jk$-degree $(-1,1)$.
Its Koszul dual is $A^! = \bigwedge(\eta)$, with $\eta$ in $jk$-degree $(1,0)$. We also consider the Koszul algebra $ \bigwedge(\xi)$ with Koszul dual $S(y)$, where again $y$ lives in $jk$-degree $(-1,1)$ and $\xi$ in $jk$-degree $(1,0)$.
We have $\bigwedge(\eta)^* \cong \bigwedge(\eta)[0]\langle -1 \rangle$, and thus have a commutative diagram
$$\xymatrix{
S(x) \otimes \bigwedge(\xi) \ar[rr]^{C_x \otimes \bigwedge(\xi)} \ar[d]_{(K_x \otimes K_{\xi})^\rho}
& & \bigwedge(\eta) \otimes \bigwedge(\xi) \ar[d]^{\rho \langle -1 \rangle} \\
S(y) \otimes \bigwedge(\eta) \ar[rr]^{C_{y} \otimes \bigwedge(\eta)} & & \bigwedge(\xi) \otimes \bigwedge(\eta)  \\
}$$
where $K$ and $C$ denote suitable Koszul complexes and their adjoints.

We record that the differential on $K_x \otimes K_{\xi}$ is given by internal multiplication by $x \otimes \xi + \xi \otimes x$, and
thus maps a basis element of the form $x^n \otimes b \otimes c \otimes d$ for $b,c,d$ in $\bigwedge(\eta)^*, \bigwedge(\xi), S(y)^*$ respectively,
to $(-1)^n(x^{n+1} \otimes \eta b \otimes c \otimes d +  x^n \otimes b \otimes c \xi \otimes y d )$.
Here the factor $(-1)^n$ is the sign obtained from the $k$-degree of $x^n$.

We have a natural isomorphism
\begin{equation*}\begin{split}K_x \otimes K_{\xi} &= S(x) \otimes \bigwedge(\eta)^* \otimes \bigwedge(\xi) \otimes S(y)^* \\
& \cong S(x) \otimes \bigwedge(\xi) \otimes \bigwedge(\eta)^* \otimes S(y)^*;\end{split}\end{equation*}
denoting $S(x) \otimes \bigwedge(\xi)$ by $B$ and $S(y) \otimes \bigwedge(\eta)$ by $\tilde B$  we obtain the $B$-$\tilde B$-bimodule $B \otimes \tilde B^{*}$
with differential acting on the left, sending an element $a \otimes b$ to the element $(-1)^{|a|}(ax \otimes \eta b + a \xi \otimes yb)$, where again the degree of $a$ is the power of $x$ appearing  
(note the twist with $\rho$ gets swallowed into $\tilde B^* = (S(y) \otimes \bigwedge(\eta))^* \cong \bigwedge(\eta)^* \otimes S(y)^*$). We fix the canonical isomorphism $\vartheta: B \to \tilde B$ which takes $x$ to $y$ and $\xi$ to $\eta$.
Our commutative diagram now takes the form
$$\xymatrix{
B \ar[rr]^{C_x \otimes \bigwedge(\xi)} \ar[d]_{(B \otimes \tilde B^{*})^{ \vartheta}}
& & \bigwedge(\eta) \otimes \bigwedge(\xi) \ar[d]^{(\bigwedge(\eta) \otimes \bigwedge(\xi))^\rho\langle -1 \rangle} \\
B \ar[rr]^{ {}^\vartheta (C_{y} \otimes \bigwedge(\eta))} & & \bigwedge(\xi) \otimes \bigwedge(\eta).  \\
}$$

This is clearly equivalent to
$$\xymatrix{
B \ar[rr]^{C_x \otimes \bigwedge(\xi)} \ar[d]_{(B \otimes \tilde B^{*})^{\vartheta}}
& & \bigwedge(\eta) \otimes \bigwedge(\xi) \ar[d]^{(\bigwedge(\eta) \otimes \bigwedge(\xi))^\rho\langle -1 \rangle} \\
B \ar[rr]^{ (C_{x} \otimes \bigwedge(\xi))^{\vartheta'^{-1}} } & & \bigwedge(\xi) \otimes \bigwedge(\eta)  \\
}$$
where ${\vartheta'}^{-1}$ is the isomorphism $\bigwedge(\xi) \otimes \bigwedge(\eta) \to \bigwedge(\eta) \otimes \bigwedge(\xi)$ which takes $\xi \otimes 1$ to $\eta \otimes 1$ and $1 \otimes \eta$ to $1 \otimes \xi$. This is again equivalent to
$$\xymatrix{
B \ar[rr]^{C_x \otimes \bigwedge(\xi)} \ar[d]_{(B \otimes \tilde B^{*})^{ \vartheta}}
& & \bigwedge(\eta) \otimes \bigwedge(\xi) \ar[d]^{(\bigwedge(\eta) \otimes \bigwedge(\xi))^{\rho\vartheta' }\langle -1 \rangle} \\
B \ar[rr]^{ C_{x} \otimes \bigwedge(\xi) } & & \bigwedge(\eta) \otimes \bigwedge(\xi)  \\
}$$

To prove the lemma we string this commutative diagram along a line. In order to do so, we place another grading (called the \emph{$f$-grading}) on these algebras as follows:
$\xi$ and $x$ have degree $-1$, $\eta$ and $y$ have degree $1$. If we denote a shift by $1$ in the $f$-grading by
$\triangleleft 1 \triangleright$ then we have $\bigwedge(\eta)^* \cong \bigwedge(\eta) \triangleleft -1 \triangleright$.
It is easy to see that by identifying the $i$th shift of the unique simple module with the vertex $i$, the category of $f$-graded modules for $\bigwedge(\eta) \otimes \bigwedge(\xi)$
is isomorphic to a category of modules over the zigzag algebra $Z$, which is the infinite dimensional quasi-hereditary algebra described in Section \ref{zigzag}.
By the same argument, the category of $f$-graded modules for $B = S(x) \otimes \bigwedge(\xi)$ is isomorphic to a category of modules over the quasi-hereditary algebra $\fc_\infty$ with quiver
$$\xymatrix{ ... & \overset{1}{\bullet}  & \ar@/^/[l]^{\xi}
 \ar@/_/[l]_{x}
\overset{2}{\bullet} & \ar@/^/[l]^{\xi}
 \ar@/_/[l]_{x}
\overset{3}{\bullet} & \ar@/^/[l]^{\xi}
 \ar@/_/[l]_{x}
 \overset{4}{\bullet} & ...
}$$ and relations $\xi x-x \xi = 0, \xi^2 = 0$.

Similarly, the category of graded modules for $\tilde B =S(y) \otimes \bigwedge(\eta)$ is isomorphic to a category of modules over the quasi-hereditary algebra $\widetilde{\fc}_\infty$ with quiver
$$\xymatrix{ ...   &\overset{1}{\bullet}   \ar@/^/[r]^{\eta}
 \ar@/_/[r]_{y} &
\overset{2}{\bullet}  \ar@/^/[r]^{\eta}
 \ar@/_/[r]_{y} &
\overset{3}{\bullet}  \ar@/^/[r]^{\eta}
 \ar@/_/[r]_{y} &
 \overset{4}{\bullet}  & ...
}$$ and relations $\eta y-y \eta = 0, \eta^2 = 0$.

We now wish to lift the above diagram to this graded setting. In order to do this, we must specify how $\vartheta: B \to \tilde B$ and $\vartheta': \bigwedge(\eta) \otimes \bigwedge(\xi) \to \bigwedge(\xi) \otimes \bigwedge(\eta)$ lift to the graded setting, by specifying what they do on $\fc_\infty^0 =\widetilde \fc_\infty^0 = Z^0$. We do this by specifying that they send the idempotent
$e_l$ to $e_{p+1-l}$ (and denote this endomorphism of  $\fc_\infty^0 =\widetilde \fc_\infty^0= Z^0$ by $\sigma$) and denote the resulting algebra isomorphisms by $\sigma_\fc$ and $\sigma_Z$ respectively. More precisely $\sigma_\fc: \fc_\infty \to \widetilde\fc_{\infty}$ sends $x$ to $y$ and $\xi$ to $\eta$ and $\sigma_Z$ interchanges $\xi$ and $\eta$.

We obtain a diagram

$$\xymatrix{
\fc_\infty \ar[rr]^{Y} \ar[d]_{(\fc_\infty \otimes_{\fc^0_{\infty}} \widetilde{\fc}_\infty^{*})^{ \sigma_\fc}}
& & Z \ar[d]^{\sigma_Z \nu\langle -1 \rangle} \\
\fc_\infty \ar[rr]^{Y} & & Z  \\
}$$

where $\nu$ denotes the algebra automorphism of $Z$ sending $e_l$ to $e_{l-1}$, $\eta$ to $\eta$ and $\xi$ to $\xi$ (coming from the fact that $\bigwedge(\eta)^* \cong \bigwedge(\eta) \triangleleft -1 \triangleright$),
and where $Y$ denotes a $\fc_\infty$-$Z$-bimodule whose adjoint has homology $\oplus_{l \in \mathbb{Z}} Fe_l$.

The differential on $\fc_\infty \otimes_{\fc^0_{\infty}} \widetilde{\fc}_\infty^{*}$
is now again given by $a \otimes b \mapsto (-1)^{|a|} (ax \otimes \eta b +a \xi \otimes yb)$, where the degree of $a$ is given by the power of $x$ appearing.

We now carefully analyse the bimodule $(\fc_\infty \otimes_{\fc^0_\infty} \widetilde{\fc}_\infty^{*})^{ \sigma_\fc}$. Note that as a $\fc^0_{\infty}$-$\fc$-bimodule $(\widetilde{\fc}_\infty^{*})^{ \sigma_\fc} \cong {}^\sigma (\fc_\infty^*)$ so $(\fc_\infty \otimes_{\fc^0_\infty} \widetilde{\fc}_\infty^{*})^{ \sigma_\fc} \cong \fc_\infty \otimes_{\fc^0_\infty} {}^\sigma(\fc_\infty^*)$ and our diagram becomes

$$\xymatrix{
\fc_\infty \ar[rr]^{Y} \ar[d]_{ \fc_\infty \otimes_{\fc^0_\infty} {}^\sigma(\fc_\infty^*)}
& & Z \ar[d]^{\sigma_Z \nu\langle -1 \rangle} \\
\fc_\infty \ar[rr]^{Y} & & Z \\
}$$

with the differential on $ \fc_\infty \otimes_{\fc^0_\infty} {}^\sigma(\fc_\infty^*)$ given by $a \otimes b \mapsto (-1)^{|a|} (ax \otimes \xi b + a \xi \otimes xb)$, where the degree of $a$ is given by the power of $x$ appearing.

To truncate this to our finite-dimensional setting, it is more convenient to work with
$\fc_\infty^* \otimes_{\fc_\infty^0} {}^\sigma \fc_\infty$ than $\fc \otimes_{\fc^0_\infty} {}^\sigma(\fc_\infty^*)$.
Note that the former is right adjoint to $\fc_{\infty} \otimes_{\fc^0_\infty} {}^\sigma \fc_\infty$,
whereas the latter it left adjoint to the same bimodule.
Since we know we are dealing with equivalences, we can choose which dg bimodule we work with.
The compatibility of adjunctions with differentials articulated in the appendix on Koszul duality concerning adjunction implies that on
all these bimodules the differentials are given by the formula
$$\ttd(a \otimes b) = (-1)^{|a|_k}(ax \otimes \xi b + a \xi \otimes x b).$$
Under the adjoint equivalences between derived categories of $Z$ and $\fc_\infty$ determined by $Y$ the subcategory of modules generated by vertices
$l$ with $l \geq r$ corresponds to the subcategory of modules given by vertices with $l \geq r$.
Indeed, since under the functors between $S(x)$ and $\bigwedge(\eta)$ determined by $C_x$ simple $\bigwedge(\eta)$-modules
correspond to injective $S(x)$-modules,
it follows that under $Y$ the simple $Z$-module indexed by $l$ corresponds to the $\fc_\infty$-module with socle $l$
and composition factors $l,l+1,l+2,...$ in ascending radical degrees.

It is a general feature of the theory of quasi-hereditary algebras that cutting at an idempotent corresponding to a set of vertices
which forms an ideal in the partial order corresponds to taking a quotient of derived categories in which the simple objects
corresponding to the complementary set of vertices are sent to zero.
We can consequently cut on both sides at the idempotent $e_{\leq p}$ given by the vertices $\leq p$ and obtain a commutative diagram

$$\xymatrix{
e_{\leq p}\fc_\infty e_{\leq p}\ar[rr]^{Y} \ar[d]_{e_{\leq p}\fc_\infty^* \otimes_{\fc_\infty^0} {}^\sigma \fc_\infty e_{\leq p}}
& &e_{\leq p} Ze_{\leq p}\ar[d]^{ e_{\leq p}(Z^{ \sigma_Z \nu})  e_{\leq p}\langle -1 \rangle} \\
e_{\leq p}\fc_\infty e_{\leq p}\ar[rr]^{Y} & &e_{\leq p}Z e_{\leq p}. \\
}$$

Computing $e_{\leq p}\fc_\infty^* \otimes_{\fc_\infty^0} {}^\sigma \fc_\infty e_{\leq p}$ we see that $e_{a}\fc_\infty^* e_k \otimes_{\fc_\infty^0} e_k{}^\sigma \fc_\infty e_{b} \neq 0$,
for $a,b \leq p$ forces $k\leq p$ and also $p+1-k \leq p$, and therefore $k \geq 1$. But $k \geq 1$ implies $a \geq 1$ and $k \geq p$ implies $b \geq 1$, so we have $1 \leq a,b,k \leq p$ and $e_{\leq p}\fc_\infty^* \otimes_{\fc_\infty^0} {}^\sigma \fc_\infty e_{\leq p} \cong \fc^* \otimes_{\fc^0} {}^\sigma \fc$. Therefore our diagram becomes

$$\xymatrix{
\fc \ar[rr]^{Y} \ar[d]_{\fc^* \otimes_{\fc^0} {}^\sigma \fc }
& & e_{\leq p} Ze_{\leq p}\ar[d]^{ e_{\leq p}(Z^{ \sigma_Z \nu})  e_{\leq p}\langle -1 \rangle} \\
\fc \ar[rr]^{Y} & & e_{\leq p}Z e_{\leq p}. \\
}$$

Twisting by the anti-automorphisms $\sigma$ of $Z$ and $\fc_{\infty}$ which exchange $e_l$ and $e_{p+l-1}$
and fix $\xi$, $\eta$ and $x$, gives us a commutative diagram

$$\xymatrix{
\fc && \ar[ll]^{P}
e_{\geq 1} Ze_{\geq 1} \\
\fc \ar[u]^{\fc^\sigma \otimes_{\fc^0} \fc^* }& & \ar[ll]^{P} e_{\geq 1}Z e_{\geq 1}\ar[u]_{ e_{\geq 1}(Z^{ \nu \sigma_Z })e_{\geq 1} \langle -1 \rangle} \\
}$$

since $P$ is the opposite of $Y$.
Note that upon taking opposites, left differentials on dg bimodules become right differentials, and vice versa.

Computing  $e_{\geq 1}(Z^{ \nu\sigma_Z })  e_{\geq 1} = (e_{\geq 1}Z e_{\leq p-1})^{\nu \sigma_Z }$, which from Section \ref{GL2} we know to be isomorphic to $\bft$ as bimodule over $e_{\geq 1} Ze_{\geq 1} e_{\geq 1}/\mathrm{Ann}(e_{\geq 1}(Z^{ \nu\sigma_Z })  e_{\geq 1}) \cong \bfc$, so our commutative diagram finally yields

$$\xymatrix{
\fc && \ar[ll]^{P}
\bfc \\
\fc  \ar[u]^{\fc^\sigma \otimes_{\fc^0} \fc^* }& &\ar[ll]^{P} \bfc \ar[u]_{ \bft \langle -1 \rangle}. \\
}$$

As above in the infinite case, $\fc^\sigma \otimes_{\fc^0} \fc^* $
induces a derived self-equivalence of $\fc$ and is left adjoint to the equivalence given by
$\fc \otimes {}^\sigma \fc$, whose right adjoint is $\fc^* \otimes_{\fc^0} {}^\sigma \fc$,
so it is left to us which adjoint we use in our computations.
\endproof

\begin{remark}
As Remark \ref{koszulselfdualalgebra} provides us with a $\fc$-$\fc^0$-bimodule isomorphism $\fc^! \cong \fc^\sigma$, we have a $\fc$-$\fc$-bimodule isomorphism between the homological dual of $\bft$ given by $\fc \otimes_{\fc^0} {}^\sigma (\fc^*)$ and the Koszul complex $\fc^! \otimes_{\fc^0} \fc^*$. Note that however the differentials are \emph{not} the same:
In the Koszul complex (after applying our bimodule iso $\fc^! \cong \fc^\sigma$),
the differential is given by $$a e_{p+1-l} \otimes e_l b \mapsto axe_{p+2-l} \otimes e_{l-1}\xi b + (-1)^l a\xi e_{p+2-l} \otimes e_{l-1}x b$$
whereas our complex has a differential
$$ a e_{p+1-l} \otimes e_l b \mapsto (-1)^{|a|_k} (axe_{p+2-l} \otimes e_{l-1}\xi b + a\xi e_{p+2-l} \otimes e_{l-1}x b)$$ where again the $k$-degree of $a$ is determined by the power of $x$ appearing rather than by idempotents as in the Koszul complex.
\end{remark}

\begin{comment}
Since $\Delta$ has top concentrated in $j$-degree $0$, and $\bft \otimes \Delta \cong \nabla$ has socle concentrated in $j$-degree $1$,
$$\Hom_\bfc(\Delta, \tilde{\bft} \otimes_\bfc\Delta) \cong \Hom_\bfc(\Delta,\nabla)$$
is concentrated in $j$-degree $1$ and $k$-degree $0$, and is, as a $\fc$-module,
isomorphic to the sum $\bigoplus_{l=1}^p \fc^0e_l \langle 1 \rangle$ over the simple $\fc$-modules
(where $\fc^0$ denotes the degree zero part in the radical grading on $j$).
The algebra $\fc$ is Koszul self-dual by Remark \ref{koszulselfdualalgebra},
and this $\fc$-$\fc$ bimodule $\fc^{0 \sigma}$,
where $\sigma$ is the automorphism of $\fc^0$ swapping idempotents $e_l$ and $e_{p-l+1}$, restricts to a one sided tilting complex for $\fc$.
A more careful analysis (see Remark \ref{homologicalselfdualitypsi}) yields that
\end{comment}

It follows from Proposition \ref{longlemma} that
$\cF(\tilde{\bft})$ is quasi-isomorphic to a twisted version of the Koszul complex for $\fc$, shifted in $j$-degree by $1$, namely to
\begin{equation}\label{maltese}
\ft := \fc^\sigma \otimes_{\fc^0} \fc^* \langle 1 \rangle
\end{equation}

with differential given by internal multiplication by $x \otimes \xi + \xi \otimes x$.

\subsection{The triple.}\label{triple}

Adjunction gives us a quasi-isomorphism
$$\fc^\sigma \otimes_{\fc^0} \fc^{* \sigma} \otimes_{\fc^0} { \fc} \rightarrow \fc.$$
It follows immediately that $\cF(\tilde{\bft}^{-1})$ is quasi-isomorphic to
\begin{equation}\label{malteseinv}
\ft^{-1} := \fc^\sigma \otimes_{\fc^0} \fc \langle -1 \rangle 
\end{equation}
with differential  $\ttd$ given by internal multiplication by $x \otimes \xi + \xi \otimes x$.
The $k$-grading on $\ft^{-1}$ is identified with the $k$-grading on $\fc^\sigma \otimes_{\fc^0} \fc \langle -1 \rangle$
inherited from the $k$-grading on $\fc$, which is \emph{different} from the usual homological $h$-grading on the Koszul complex;
nevertheless the differential has $k$-degree $1$.

Putting these observations together and setting $\underline\ft = (\ft,\ft^{-1})$, we have

\begin{prop}\label{daggermaltese}
The triple $(\fc, \underline\ft)$ is quasi-isomorphic to $\bbF(\bfc, \ut)$.
\end{prop}

\section{Expressing $\mathbf{w}$ via $(\fc, (\ft, \ft^{-1}))$.} \label{computation}

We demonstrate here how the analysis of homological duality of algebraic operators made above
can be used to reduce the computation of the Weyl extension algebra $\mathbf{w}$ of the principal block of $GL_2$
to the computation of the homology of a certain tensor algebra.

\begin{prop} \label{stock}
We have $\mathbf{w_q} \cong \fO_{ F} \fO_{\mathbb{H}\bbT_{\fc}(\underline\ft)}^q (F[z]).$
\end{prop}

\proof
We have algebra isomorphisms

$$\begin{array}{lll}
\mathbf{w_q} & \cong \mathbb{H} \bbP_{F,0} \bbF \bbP_{\bfc,\ut}^q(F,(F,F)) & \mbox{Theorem \ref{GL2theorem}} \\
& \cong \mathbb{H} \bbP_{F,0} \bbP_{\bbF(\bfc,\ut)}^q \bbF(F,(F,F)) &  \mbox{Theorem \ref{chainquasi}, Proposition \ref{conds}}\\
& \cong \mathbb{H} \bbP_{F,0} \bbP_{\bbF(\bfc,\ut)}^q(F,(F,F)) & \mbox{$\bbF(F,(F,F))=(F,(F,F))$} \\
& \cong \mathbb{H} \fO_F \fO_{\bbT_{\cF(\bfc)}(\cF(\ut))}^q(F[z,z^{-1}]) &  \mbox{Lemma \ref{compare}}\\
& \cong \fO_{\mathbb{H} F} \fO_{\mathbb{H}\bbT_{\cF(\bfc)}(\cF(\ut))}^q \mathbb{H}(F[z,z^{-1}]) &  \mbox{Lemma \ref{commute}}\\
& \cong \fO_{\mathbb{H} F} \fO_{\mathbb{H}\bbT_{\cF(\bfc)}(\cF(\ut))}^q (F[z,z^{-1}]) &  \mbox{$\mathbb{H}(F[z]) = F[z]$}\\
&  \cong \fO_{ F} \fO_{\mathbb{H}\bbT_{\fc}(\underline\ft)}^q (F[z,z^{-1}])   & \mbox{Proposition \ref{daggermaltese}}
.\end{array}$$ \endproof

In the following section we study the algebra $\Upsilon:=\mathbb{H}\bbT_{\fc}(\underline\ft)$
and show that it has the description referred to in the Section \ref{monomial}, which will then complete the proof of Theorem \ref{Yoneda}.

\section{The algebra $\Upsilon$.}\label{tensoralg}

The aim of this section is to obtain a combinatorial description of 
$\mathbb{HT}_{\fc}(\underline{\ft})$, which we will denote by $\Upsilon.$ We will first show that we only need to consider a certain subspace. Then we will analyse the $\fc$-$\fc$-bimodule components of this subspace,
then how to multiply these components,
before reducing the structure to a basis given by points in a polytope.

\subsection{Truncating $\Upsilon$.}

It turns out that to describe the algebra structure on $\fO_{F} \fO_{\Upsilon}^q(F[z])$
we only need to know about multiplication on the part of $\Upsilon$ featuring $\mathbb{H}(\ft^i)$ for $i \leq 1$.
Here we justify this fact.

Note that for $i \in \ZZ$,

\begin{equation}\label{iter}\begin{split}
 \fO_{\mathbb{H}\bbT_{\fc}(\ft, \ft^{-1})}^q& \mathbb{H}(F[z,z^{-1}])^{i,k}\\& = \bigoplus_{\substack{j_1,\\ k_1 + \tilde k_1 = k}} (\mathbb{H}(\ft^{\otimes i}))^{j_1,k_1} \otimes_F (\fO_{\mathbb{H}\bbT_{\fc}(\ft, \ft^{-1})}^{q-1} \mathbb{H}(F[z,z^{-1}]))^{j_1,\tilde k_1}\\
&=\bigoplus_{\substack{j_1,\\ k_1 + \tilde k_1 = k}} (\mathbb{H}(\ft^{\otimes i}))^{j_1,k_1} \otimes_F \\
& \left( \bigoplus_{\substack{j_2,\\ k_2 + \tilde k_2 = \tilde k_1}} ( \mathbb{H}(\ft^{\otimes j_1}))^{j_2,k_2} \otimes_F (\fO_{\mathbb{H}\bbT_{\fc}(\ft, \ft^{-1})}^{q-2} \mathbb{H}(F[z,z^{-1}]))^{j_2,\tilde k_2}  )   \right) \\
& = \cdots\\
\end{split}
\end{equation}

where, for $l <0$  we interpret $\ft^{\otimes l }$ as $(\ft^{-1})^{\otimes -l}$.

A typical direct summand of this looks like

$$(\mathbb{H}(\ft^{\otimes i}))^{j_1,k_1} \otimes_F (\mathbb{H}(\ft^{\otimes j_1}))^{j_2,k_2}  \otimes_F \cdots \otimes_F (\mathbb{H}(\ft^{\otimes j_{q-1}}))^{j_q,k_q}.  $$

\begin{lem} \label{trunc}
For $i \leq 1$ and any $k \in \ZZ$, no direct summand of $\fO_{\mathbb{H}\bbT_{\fc}(\ft, \ft^{-1})}^q \mathbb{H}(F[z,z^{-1}])^{i,k }$ involves tensor factors $\mathbb{H}(\ft^{\otimes i})^{j,k} $ for $i >1$.
\end{lem}

\proof
The proof is by induction on $q$. The case $q=1$ is trivial.
Assume it is true for $\fO_{\mathbb{H}\bbT_{\fc}(\ft, \ft^{-1})}^{q-1} \mathbb{H}(F[z,z^{-1}])^{i,k }$.
The first step in \eqref{iter} together with the observation that for $i \leq 1$,
the graded piece $\mathbb{H}(\ft^{\otimes i})^{j, \diamond}$ is zero for all $j >1$ then implies the inductive step and the Lemma.
\endproof

\subsection{The homology of $\ft$.}

\begin{lem}\label{hmalt}
As a $\fc^0$-$\fc^0$-bimodule, $\bbH(\ft) \cong (\fc^0)^\sigma$.  
\end{lem}

\proof
This follows from direct computation. The basis elements spanning the homology are $e_i \otimes e_{p+1-i}^*$ for $i=1,\dots p$.
\endproof

\begin{remark}\label{truncatemaltese}
We remark that both as a left and as a right dg module, $\ft$ is quasi-isomorphic to its homology, but not as a dg bimodule. However the projection map $\fc^\sigma \otimes_{\fc^0} \fc^* \twoheadrightarrow Fe_p \otimes Fe_1^*$ is a morphism of dg bimodules, and we denote by 
$\widehat{\fc^0}$ the kernel of this projection;
it is quasi-isomorphic to $\overline{\fc^0} = \oplus_{h=1}^{p-1} (\fc^{0}e_h)^\sigma$ as a right and left dg module.
\end{remark}

\subsection{A bimodule.}\label{bimodule}

Here we define and study a certain bimodule $M$ for the algebra $\fc$ which is prominent in the homology of $\ft^{i}$.

We define $L_l \in \fc \ml$ to be the module with basis $\{ x.x^l, \xi.x^l  | 0 \leq l \leq p-1 \}$ where $e_f$ acts as the identity
on $ x.x^{p-f}$ and $\xi.x^{p-f}$, the generator $x$ acts on the second component in the obvious way,
and $\xi$ sends $x.x^l$ to $\xi.x^{l+1}$, whilst killing $\xi.x^{p-l}$.
Thus $L_l$ is a left $\fc$-module of dimension $2p$, whose Loewy structure is
{\footnotesize
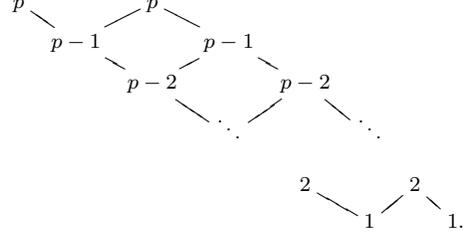
\xymatrix@C=5pt@R=4pt{&p\ar@{-}[dr]&&p\ar@{-}[dl]\ar@{-}[dr]&&&&&&&\\
&&p-1\ar@{-}[dr]&&p-1\ar@{-}[dl]\ar@{-}[dr]&&&&&&\\
&& &p-2 \ar@{-}[dr] &&p-2 \ar@{-}[dl]\ar@{-}[dr]&&&&&\\
&&&&\ddots &&\ddots &&&&\\
&&&&&2\ar@{-}[dr]&&2\ar@{-}[dl]\ar@{-}[dr]&&&\\
&&&&&&1&&1.&&}
}

We define $L_r \in \mr \fc$ to be the module with basis $\{ x^l.x, x^l.\xi | 0 \leq l \leq p-1 \} $ where $e_l$ acts as the identity
on $x^{l-1}.x$ and $x^{l-1}.\xi$, the generator $x$ acts on the first component in the obvious way and $\xi$ sends $x^l.x$ to $x^{l+1}.\xi$.
Thus $L_r$ is a right $\fc$-module of dimension $2p$, whose Loewy structure is

{\footnotesize
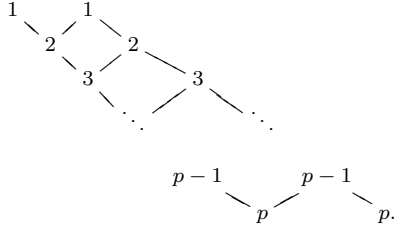
\xymatrix@C=5pt@R=3pt{&1\ar@{-}[dr]&&1\ar@{-}[dl]\ar@{-}[dr]&&&&&&&\\
&&2\ar@{-}[dr]&&2\ar@{-}[dl]\ar@{-}[dr]&&&&&&\\
&& &3 \ar@{-}[dr] &&3 \ar@{-}[dl]\ar@{-}[dr]&&&&&\\
&&&&\ddots &&\ddots &&&&\\
&&&&&p-1\ar@{-}[dr]&&p-1\ar@{-}[dl]\ar@{-}[dr]&&&\\
&&&&&&p&&p.&&}
}

The algebra $\fc$ has a $k$-grading,
with the symbol $x$ in degree $1$ and the symbol $\xi$ in degree $0$.
This gives $\fc$ the structure of a dg algebra with trivial differential.
The dg bimodules $\ft$ and $\ft^{-1}$ admit natural $k$-gradings, with $x$ in degree $1$ and $\xi$ in degree $0$.

Both $L_l$ and $L_r$ carry a natural $j$-grading,
by placing the symbols $x$ and $\xi$ appearing in the basis elements in degrees $-1$ and $1$ respectively.
We give $L_l$ and $L_r$ a $k$-grading by placing $\xi.x^l$, $x^l. \xi$ in degree $l-1$ and $x.x^l$, $x^l.x$ in degree $l$.

\medskip 

We denote by $K'_l = \fc^\sigma \otimes \fc^*$ the twisted left Koszul complex, with differential given by
$\ttd(a \otimes b) = (-1)^{|a|_k}(ax \otimes \xi b + a \xi \otimes x b)$.
We denote by $K'_r = \fc^* \otimes {}^\sigma \fc$ the twisted right Koszul complex, with differential given by
$\ttd(a \otimes b) = (-1)^{|a|_k}(ax \otimes \xi b + a \xi \otimes x b)$.

\begin{lem}\label{Ls}
\begin{enumerate}[(i)]
\item\label{27i} We have $L_l^* \cong L_r \langle p-1 \rangle[2-p]$ and $L_r^* \cong L_l \langle p-1 \rangle[2-p]$.
\item\label{27ii} As ungraded modules, every nonsplit extension of the injective $\fc^*e_h$ by the projective $\fc e_h$ is isomorphic to $L_l$;
every nonsplit extension of the injective $e_h\fc^*$ by the projective $e_h\fc$ is isomorphic to $L_r$.
\item\label{27iii} $K'_l \otimes_\fc L_l$ is quasi-isomorphic to $L_l \langle p-1 \rangle[1-p]$ and $L_r \otimes_\fc K'_r$ is quasi-isomorphic to $ L_r \langle p-1\rangle[1-p]$.
\item\label{27iv} As ungraded modules, we have $\Ext^m_\fc(L_\bullet, \fc) = \Ext^m_\fc(\fc^*, L_\bullet) =0$, for $0 \leq m \leq p-2$, and $\bullet \in \{l,r\}$.
\item\label{27v} There is a unique nonsplit extension of $L_\bullet$ by an irreducible module for $\bullet \in \{ l,r \}$, forming the middle terms of short exact sequences
$$0 \rightarrow \fc^0 e_{p-1} \rightarrow E_l \rightarrow L_l \rightarrow 0,$$
$$0 \rightarrow e_2 \fc^0 \rightarrow E_r \rightarrow L_r \rightarrow 0.$$
\end{enumerate}
\end{lem}

\proof
Both \eqref{27i} and \eqref{27ii} are proved by easy explicit calculations.
The proof of \eqref{27iii} is a little more subtle as we need to drag $j$ and $k$-gradings through Koszul duality.
The algebra $\fc$ is Koszul self-dual.
We consider the image in $D^b(\fc)$ of modules
under the endofunctor $K' \otimes_\fc -$, where $K' = K'_l$ is the twisted Koszul complex for $\fc$.
Here the $k$-grading on the dg $\fc$-$\fc$ bimodule $K'$ is inherited from the $k$-grading on $\fc$,
and \emph{not} the homological $h$-grading on the Koszul complex.
Classical Koszul duality sends projectives to simples and simples to injectives.
Likewise, here we know that $K' \otimes_\fc \fc^0e_s$ is quasi-isomorphic to $\fc^*e_{p+1-s}$
and $K' \otimes_\fc \fc e_s$ is quasi-isomorphic to $\fc^0 e_{p+1-s}$.
Applying $K' \otimes_{\fc} -$ to the exact triangle $\fc e_p \rightarrow L_l \rightarrow \fc^0e_p \rightsquigarrow$
gives us an exact triangle
$\fc^0e_1 \rightarrow K' \otimes_\fc L_l \rightarrow \fc^* e_1 \rightsquigarrow$.
The unique extension of $\fc^* e_1$ by $\fc^0e_1$ is $L_l$ and thus $K' \otimes_{\fc} L_l$ is quasi-isomorphic to $L_l$.
To see the gradings, note that a projective resolution of $L_l$ is
$$P_1 \langle -(p-2) \rangle [-1] \oplus P_1 \langle -p \rangle [0] \to \cdots \hspace{5cm}$$
$$\hspace{2cm} \cdots\to P_{p-1} \langle 0\rangle [-1] \oplus P_{p-1} \langle -2\rangle[0] \to
P_p \langle 1\rangle [-1] \oplus P_p \langle -1\rangle [0].$$
This has a filtration with sections $P_{p+k} \langle k+1 \rangle[k-1] \oplus P_{p+k} \langle k-1 \rangle[k]$, for $k=0,-1,...,1-p$.
Tensoring with $K'$, which is quasi-isomorphic to $\fc^{0 \sigma}$, we obtain $L_l \langle p-1 \rangle[1-p]$ as required.
The proof for $L_r$ is similar.

Claim \eqref{27iv} follows from our projective resolution of $L_l$ above because applying $\Hom_\fc(-,\fc)$ to that linear resolution
gives a linear resolution of $L_r$, whose homology is concentrated in a single degree.

Claim \eqref{27v} again follows from explicit computation.
\endproof

We have an endomorphism $t$ of $L_l$ of $jk$-degree $(2,-1)$ which sends $x.x^l$ to $\xi.x^l$, and $\xi.x^l$ to zero.
Writing $\Lambda$ for the exterior algebra on $t$,
we thus give $L_l$ the structure of a $\fc$-$\Lambda$-bimodule.
We have an endomorphism $t$ of $L_r$ of degree $2$ which sends $x^l.x$ to $x^l.\xi$, and $x^l.\xi$ to zero.
We thus give $L_r$ the structure of a $\Lambda$-$\fc$ bimodule.

We define $M$ to be the $\fc$-$\fc$-bimodule $L_l \otimes_\Lambda L_r$.
We define $\tau$ to be the involution on $\fc$ that sends $x$ to $x$ and $\xi$ to $- \xi$.
The bimodule $M$ has some intriguing properties:

\begin{prop} \label{Ms}
\begin{enumerate}[(i)]
\item\label{28i} We have $M^* \cong M \langle 2p \rangle[-2p+3]$.

\item\label{28ii}  We have $_{\fc}M \cong \bigoplus_{h=0}^{p-1} L_l \langle -1-h \rangle[h]$ and $M_{\fc} \cong \bigoplus_{h=0}^{p-1} L_r \langle -1-h\rangle[h]$.

\item\label{28iii} We have a short exact sequence of bimodules
$$0 \rightarrow \fc \langle -p-1 \rangle[p-1] \rightarrow M \rightarrow \fc^{*} \langle -p+1 \rangle [p-2]\rightarrow 0.$$

\item\label{28iv}  We have a quasi-isomorphism between $\ft^i \otimes_\fc M$ and
$M^{\tau^i} \langle ip \rangle[i(1-p)]$, for $i \in \mathbb{Z}$. Similarly we have a quasi-isomorphism between $M \otimes_\fc \ft^i$ and $M^{\tau^i} \langle ip \rangle[i(1-p)]$, for $i \in \mathbb{Z}$.

\item\label{28v}  We have an exact triangle in the derived category of differential graded $\fc$-$\fc$-bimodules
$$\ft^{-1} \rightarrow M^\tau \langle 1 \rangle[0] \rightarrow %\fc^\sigma \otimes_{\fc^0} \fc^* \langle 1 \rangle
\ft[-1] \rightsquigarrow.$$

\item\label{28vi}  We have an exact triangle in the derived category of differential graded $\fc$-$\fc$-bimodules
$$\ft^{-2} \rightarrow M \langle -p+1 \rangle[p-1] \rightarrow \fc \langle 0 \rangle[-1] \rightsquigarrow.$$
\end{enumerate}
Assume $p \geq 3$.
\begin{enumerate}[(i)]\setcounter{enumi}{6}
\item\label{28vii} We have $\Hom_{\fc}(M,M) \cong M \langle p+1 \rangle[1-p]$ and $M \otimes_{\fc} M \cong M \langle 1-p \rangle[p-1]$ as $jk$-graded $\fc$-$\fc$-bimodules.
\end{enumerate}
\end{prop}
\proof
\eqref{28i} We have isomorphisms of bimodules
\begin{equation*}
\begin{split}
M^* &= \Hom(L_l \otimes_\Lambda L_r, F) \\
& \cong \Hom_\Lambda( L_l, \Hom(L_r,F)) \\
& \cong \Hom_\Lambda(L_l, L_l)  \\
& \cong \Hom_\Lambda(L_l, L_l \otimes_\Lambda \Lambda))  \\
& \cong L_l \otimes_\Lambda \Hom_\Lambda(L_l, \Lambda)  \\
& \cong L_r \otimes_\Lambda L_r \\
& \cong M
\end{split}
\end{equation*}
We are using here that $L_l$ is projective as a right $\Lambda$-module, some adjunctions,
and the fact that $\Lambda$ is a symmetric algebra.
The gradings match up as described: $M$ is concentrated in $j$-degrees $0$ down to $-2p$,
thus $M^*$ is concentrated in degrees $2p$ down to $0$; $M$ is concentrated in $k$-degrees $-1$ up to $2p-2$,
thus $M^*$ is concentrated in degrees $2-2p$ up to $1$.

\eqref{28ii} This follows directly from the definition.

\eqref{28iii} To embed $\fc$ in $M$, we send $e_h \in \fc$ to $x.x^{p-h} \otimes x^{h-1}.x \in M$, a map of $(j,k)$-degree
$(-p-1,p-1)$. The quotient map corresponding to this embedding is then dual to the embedding.

\eqref{28iv} As one-sided modules, this follows from Lemma \ref{Ls} \eqref{27iii}.
We know that $\ft$ is $K_l$ accompanied by a shift in $j$-degree by $1$.
It follows that $\ft \otimes_\fc L_l = L_l \langle p \rangle[1-p]$
and consequently that $\ft \otimes_\fc M = M \langle p \rangle[1-p]$.
The left-sided statement follows from the fact that $K' \otimes_\fc M \cong M ^\tau$ as $\fc$-$\fc$-bimodules,
the right-sided statement similarly.
We do not establish the twist by $\tau$ here;
a detailed analysis confirming the twist appears in Corollary \ref{twistbytau} and Lemma \ref{notwistbytau}.

\eqref{28v} The triangle is obtained by tensoring the exact sequence of part \eqref{28iii} on the left with
$\ft^{-1}\langle 1\rangle =\fc^\sigma \otimes_{\fc^0} \fc$, shifting in $j$-degree by $p$ and using part \eqref{27iv}.

\eqref{28vi} The triangle is obtained by tensoring the triangle of part \eqref{28v} on the left with
$\ft^{-1} = \fc^\sigma \otimes \fc \langle -1 \rangle$, again using part \eqref{27iv}.

\eqref{28vii} First note that applying $\Hom_\fc(-,M)$ to the short exact sequence of \eqref{28iii} gives us a long exact sequence
$$... \leftarrow \Ext^1(\fc^*,M) \leftarrow \Hom(\fc \langle -p-1 \rangle[p-1],M) \leftarrow \Hom(M,M) \leftarrow \Hom(\fc^*,M) \leftarrow 0.$$
By the preceding lemma $\Hom_\fc(\fc^*,M) = \Ext^1_\fc(\fc^*,M) = 0$, and so $$\Hom_\fc(M,M) \cong M \langle p+1 \rangle[1-p].$$
Since the map
$\Hom_\fc(M,M) \rightarrow \Hom_\fc(\fc,M)$ is a bimodule homomorphism,
we have $\Hom_\fc(M,M) \cong M \langle p+1 \rangle[1-p]$ as $jk$-graded $\fc$-$\fc$-bimodules.
By adjunction and part \eqref{28i}  and  the previous sentence, we have
\begin{equation*}\begin{split}\Hom_F(M \otimes_\fc M, F) & \cong \Hom_\fc(M,M^*)\\
&\cong\Hom_\fc(M,M \langle 2p \rangle[3-2p]) \\
&\cong M \langle 3p+1 \rangle[4-3p]
\end{split}
\end{equation*}
so duality gives us $M \otimes_\fc M \cong M \langle -p-1 \rangle[p-1]$ by part \eqref{28i}.
\endproof

\begin{lem}
We have a bimodule $L_l^\theta  \otimes_\Lambda L_r$ where $\theta$ is the automorphism of $\Lambda$ sending $t$ to $-t$.
As a $\fc$-$\fc$-bimodule, $L_l^\theta  \otimes L_r$ is isomorphic to $M^\tau$.
\end{lem}
\proof
On both $L_l^\theta \otimes_\Lambda L_r$ and $M^\tau$,
internal multiplication by $x \otimes \xi + \xi \otimes x$ is zero. Since the spaces have the same dimension
and are quotients of $L_l \otimes_F L_r$, they are isomorphic.
\endproof

It is convenient to truncate $M$ a little, as both $M$ and its truncated version will appear in our explicit computations of $\mathbb{HT}_{\fc}(\underline \ft)$.
We have a $j$-graded bimodule homomorphism
$M \rightarrow (\fc^0 e_p)^\sigma$ which sends $x.1 \otimes 1.\xi \in M$ to $e_p$,
and sends all other basis elements of $M$ to zero; we define $\overline{M}$ to be the kernel of this homomorphism.
Recall the dg bimodule $\widehat{\fc^0}$ from Remark \ref{truncatemaltese}.

\begin{prop}\label{Ms2}
\begin{enumerate}[(i)]
\item\label{29i} We have an exact triangle in the derived category of $\fc$-$\fc$-bimodules,
$$\ft^{-1} \rightarrow \overline{M}^\tau \langle 1 \rangle \rightarrow \widehat{\fc^0} \langle 1 \rangle [-1] \rightsquigarrow.$$
\item\label{29ii} The homology of the left $\fc$-module $\ft^{-1}e_h$ is isomorphic to
$\fc e_p$ if $h=1$, to $E_l$ if $h=2$, and to $L_l \oplus \fc^0e_h$ if $h>2$.
The homology of the right $\fc$-module $e_h \ft^{-1}$ is isomorphic to
$e_1\fc$ if $h=p$, to $E_r$ if $h=p-1$, and to $L_r \oplus e_i\fc^0$ if $h<p-1$.
\end{enumerate}
Assume $p \geq 3$.
\begin{enumerate}[(i)]\setcounter{enumi}{2}
\item\label{29iii} We have $$\Hom_{\fc}(\overline{M},M) \cong M \langle p+1\rangle [1-p]$$ and $$\overline{M}\otimes_{\fc} M \cong M \cong
M \otimes_{\fc} \overline{M} \cong \overline{M}\otimes_{\fc}\overline{M} \cong M  \langle -p-1\rangle [p-1]$$
as $\fc$-$\fc$-bimodules.
\end{enumerate}
\end{prop}
\proof

\eqref{29i}
By Proposition \ref{Ms}\eqref{28v}, we have a triangle $\ft^{-1} \rightarrow M^\tau \langle 1 \rangle[0] \rightarrow %\fc^\sigma \otimes_{\fc^0} \fc^* \langle 1 \rangle
\ft[-1] \rightsquigarrow.$ We obtain the desired triangle by cancelling the one-dimensional cokernels of $\overline{M} \to M$ and $\widehat{\fc^0} \to \ft$ with the induced non-zero map between them, that is zero in the derived category.

\eqref{29ii} Using the quasi-isomorphism between $\ft^{-1}$ and $\overline{M}^\tau \rightarrow \widehat{\fc^0}[-1]$ from Proposition \ref{Ms}\eqref{28v} as well as the one-sided quasi-isomorphisms $\widehat{\fc^0} \cong \overline{\fc^0}$, this follows by direct computation.

\eqref{29iii} We have an exact sequence of bimodules
$$0 \rightarrow \overline{M} \rightarrow M \rightarrow (\fc^0 e_p)^\sigma \rightarrow 0.$$
Applying $\Hom(-,M)$ in the category of left modules gives us a long exact sequence
$$\leftarrow \Ext^1(\fc^0 e_p,M) \leftarrow \Hom(\overline{M},M) \leftarrow \Hom(M,M) \leftarrow \Hom(\fc^0 e_p,M).$$
Since $\Hom(\fc^0 e_p,M)$ and $\Ext^1(\fc^0 e_p,M)$ are both zero, we find $$\Hom(\overline{M},M) \cong \Hom(M,M) \cong M.$$
As ungraded modules, $$\Hom_{F}(M \otimes_{\fc} \overline{M} , F) \cong \Hom_{\fc}(\overline{M} ,M^*) \cong \Hom_{\fc}(\overline{M} ,M) \cong M,  $$ and dualising (using again $M^* \cong M$ from Proposition \ref{Ms}\eqref{28i}), we obtain $M \otimes_{\fc} \overline{M} \cong  M$.
Working instead on the right, we obtain $\overline{M} \otimes_\fc M \cong M$.
By duality, we have an exact sequence of bimodules
$$0 \leftarrow \overline{M}^* \leftarrow M \leftarrow  (\fc^0 e_1)^\sigma \leftarrow 0.$$
Applying $\Hom(\overline{M},-)$ in the category of left modules gives us a long exact sequence
$$\Ext^1(\overline{M},(\fc^0 e_1)^\sigma) \leftarrow \Hom(\overline{M}, \overline{M}^*) \leftarrow \Hom(\overline{M},M) \leftarrow \Hom(\overline{M},(\fc^0 e_1)^\sigma) \leftarrow 0$$
Vanishing of the first and last terms of this sequence gives us an isomorphism
$\Hom(\overline{M}, \overline{M}^*) \cong \Hom(\overline{M}, M)$, so $\Hom(\overline{M}, \overline{M}^*) \cong M\langle -p\rangle$.
Adjunction and the isomorphism $M^*\cong M$ again give us $\overline{M} \otimes_\fc \overline{M} \cong M$. All gradings are checked using the gradings prescribed in Proposition \ref{Ms}.
\endproof

\begin{lem}\label{mapsgamma}
We have natural homomorphisms of $\fc$-$\fc$-bimodules $\alpha$, $\beta$, and $\gamma$, that make the triangles
$$\xymatrix@=14pt{
M \otimes_\fc \fc^{0 \sigma} \ar[r]^\alpha \ar[d]^\gamma & M e_p \langle p-1 \rangle[1-p] \ar[dl]^\beta & \fc^{0 \sigma} \otimes_\fc M \ar[r]^\alpha \ar[d]^\gamma & e_1 M  \langle p-1 \rangle [1-p] \ar[dl]^\beta \\
M \langle p-1 \rangle[1-p] &  & M \langle p-1 \rangle[1-p] & \\
}$$
$$\xymatrix{
\overline{M} \otimes_\fc \fc^{0 \sigma} \ar[r]^\alpha \ar[d]^\gamma & \fc e_p \langle -2 \rangle[0] \ar[dl]^\beta & \fc^{0 \sigma} \otimes_\fc \overline{M} \ar[r]^\alpha \ar[d]^\gamma & e_1 \fc  \langle -2 \rangle[0] \ar[dl]^\beta \\
\fc \langle -2 \rangle[0] &  & \fc \langle -2 \rangle[0] & \\
}$$
commute.
\end{lem}
\proof
The arrows $\beta$ are the obvious embeddings. The arrow $\alpha$ in the top left hand diagram sends an element
$a \otimes b \otimes e_l \in L_l \otimes L_r \otimes \fc^{0 \sigma}$ to $a \otimes b e_p$ for $l=1$
and to zero for $l>1$;
The arrow $\alpha$ in the top right hand diagram sends an element
$e_l \otimes b \otimes c \in \fc^{0 \sigma} \otimes L_l \otimes L_r$ to $e_1 b \otimes c$ for $l=1$ and to zero for $l>1$.
These maps restrict to the arrows $\alpha$ in the bottom diagram.
The arrows $\gamma$ are merely compositions of $\alpha$ and $\beta$.
\endproof

\subsection{The $\fc$-$\fc$-bimodule structure of $\mathbb{T}_\fc(\ft^{-1})$.}

We now obtain the following description of $\mathbb{T}_\fc(\ft^{-1})$. 

\begin{prop}\label{dgminus} As a differential $jk$-graded $\fc$-$\fc$-bimodule, $\mathbb{T}_\fc(\ft^{-1})$ has a filtration with sections
$$\xymatrix@C=5pt@R=4pt{
          &           &                                               &             & \fc &                                  &        &        & \\
          &           &                             &   \overline{M}^\tau \langle 1 \rangle \ar[rr] &      & \widehat{\fc^0}\langle 1 \rangle &    &        &  \\
          &           & M\langle -p+1 \rangle \ar[rr] &                               & \fc &          &        &        & \\
          & M^\tau \ar[rr]\langle -2p+1 \rangle &     &  \overline{M}^\tau \langle 1 \rangle \ar[rr] &   & \widehat{\fc^0}\langle1\rangle  & & & \\
M \ar[rr] \langle -3p+1 \rangle &    & M \ar[rr] \langle -p+1 \rangle &       & \fc &          &        &        & \\
          &           &           &             &......&          &        &        &     \\
}$$
in the $j$-grading, and
$$\xymatrix@C=5pt@R=4pt{
          &           &                                               &             & \fc &                                  &        &        & \\
          &           &                             &   \overline{M}^\tau \ar[rr] &      & \widehat{\fc^0}[-1] &    &        &  \\
          &           & M[p-1] \ar[rr] &                               & \fc [-1] &          &        &        & \\
          & M^\tau \ar[rr][2p-2] &     &  \overline{M}^\tau [-1] \ar[rr] &   & \widehat{\fc^0}[-2]  & & & \\
M \ar[rr][3p-3] &    & M \ar[rr] [p-2] &       & \fc[-2] &          &        &        & \\
          &           &           &             &......&          &        &        &     \\
}$$
in the $k$-grading.
\end{prop}
\proof
We have separated $j$-gradings and $k$-gradings for typesetting reasons: the diagrams get too complicated otherwise.
In the statement of the lemma, and in similar statements in our script, we take a sequence of arrows
$$X_0 \rightarrow X_1 \rightarrow ...  \rightarrow X_l$$
to represent an object of the derived category with a filtration whose sections are isomorphic to $X_l[l]$;
if we replace $X_l$ by an isomorphic projective resolution,
then the object is isomorphic to a complex formed by taking the total complex of a complex of such resolutions.
The rows of the diagram correspond to $\ft^{i}$, for $i \in \mathbb{Z}_{\leq 0}$.
Cases $i=-1,-2$ are described in Proposition \ref{Ms2}\eqref{29i} and \ref{Ms}\eqref{28vi}, and we obtain the lower rows by successively applying $\ft^{-1} \otimes -$ to the terms in the row above using Proposition \ref{Ms}\eqref{28iv}.
\endproof

\subsection{Explicit formulas for the homology of $\ft^{i}$ for $i \leq -1$.}\label{malteser}

We define $\Upsilon^-$ to be the homology of $\mathbb{T}_\fc(\ft^{-1})$.

Let $i \leq -1$. We are interested in describing explicit formulas for a basis of $\Upsilon^{i}$, as a subquotient of
$$\ft^{i} \cong \fc^\sigma \otimes_{\fc^0} \fc^\sigma \otimes_{\fc^0} ... \otimes \fc^\sigma \otimes_{\fc^0} \fc \langle i \rangle$$
To simplify notation, in the remains of this section we omit shifts in $jk$-degree.

\medskip 

{\bf The case $i=-1$.} The differential on $\ft^{-1} = \fc^\sigma \otimes \fc \langle -1 \rangle$ is given by
$\ttd(a \otimes b) = (-1)^{|a|_k}(ax \otimes \xi b + a \xi \otimes x b).$
We record that $x$ and $\xi$ super-commute in homology:
\begin{lem}\label{xxicommute}
In $\mathbb{H}(\fc^\sigma \otimes \fc)$, we have $x^{m} e_h\otimes e_{h} \xi x^{l} \equiv - x^{m-1}\xi e_h\otimes e_{h} x^{l+1}$.
\end{lem}
\proof
This follows immediately from the fact that the image of $x^{m-1}e_{h-1}\otimes e_{h-1} x^{l}$
under the differential is $(-1)^{m-1}(x^{m} e_h \otimes e_{h} \xi x^{l} + x^{m-1}\xi e_h \otimes e_{h} x^{ l+1})$.
\endproof

We deduce from this the twist required on our copy of $\overline{M}$ in homology:

\begin{cor} \label{twistbytau}
The component of $\mathbb{H}(\fc^\sigma \otimes_{\fc^0}  \fc)$
generated by $e_p \otimes e_1$ is isomorphic to $\overline{M}^\tau$ as a $\fc$-$\fc$-bimodule.
\end{cor}
\proof
The element $e_p \otimes e_p$ of $\fc^\sigma \otimes_{\fc^0} \fc$ maps to zero under the differential.
Since the bimodule $\overline{M}$ has simple top $p \otimes_F 1^{op}\langle -1 \rangle$
and we have a unique composition factor of $\fc^\sigma \otimes \fc$ isomorphic to $p \otimes 1^{op}$,
we conclude the element $e_p \otimes e_1$ of $\fc^\sigma \otimes \fc$
generates the factor $\overline{M}^\tau$ of $\mathbb{H}(\ft^{-1})$,
whose basis elements written into pictures of the composition structures for $\overline{M}^\tau$ look like

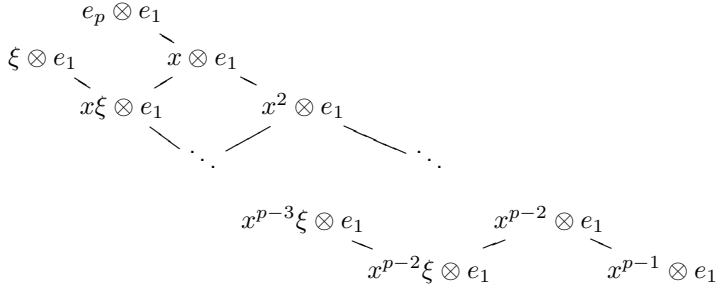
\xymatrix@C=-5pt@R=3pt{&e_p \otimes e_1\ar@{-}[dl]\ar@{-}[dr]&&&&&\\
\xi \otimes e_1 \ar@{-}[dr]&&x \otimes e_1\ar@{-}[dl]\ar@{-}[dr]&&&&\\
 &x\xi \otimes e_1  \ar@{-}[dr] &&x^2 \otimes e_1\ar@{-}[dl]\ar@{-}[dr]&&&\\
&&\ddots &&\ddots &&\\
&&& x^{p-3}\xi \otimes e_1 \ar@{-}[dr]&&x^{p-2} \otimes e_1\ar@{-}[dl]\ar@{-}[dr]&\\
&&&&x^{p-2}\xi \otimes e_1&&x^{p-1}\otimes e_1 }

for $\mathbb{H}(\fc^\sigma \otimes_{\fc^0}  \fc)e_1$, and

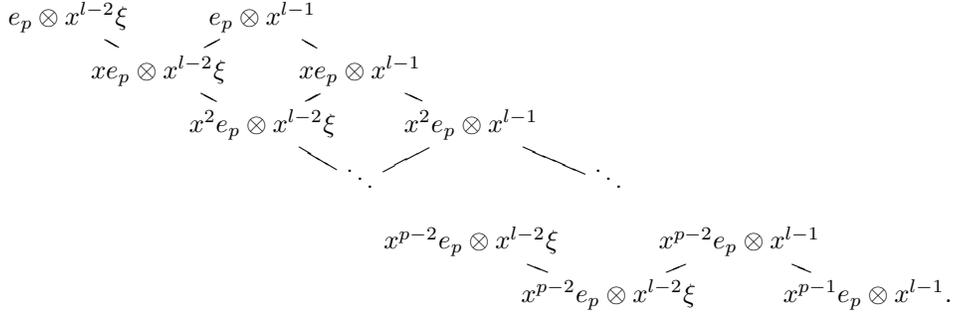
\xymatrix@C=-20pt@R=3pt{&e_p \otimes x^{l-2}\xi\ar@{-}[dr]&&e_p \otimes x^{l-1} \ar@{-}[dl]\ar@{-}[dr]&&&&&\\
&& x e_p \otimes x^{l-2}\xi\ar@{-}[dr]&& x e_p \otimes x^{l-1}\ar@{-}[dl]\ar@{-}[dr]&&&&\\
&& &x^2 e_p \otimes x^{l-2}\xi \ar@{-}[dr] &&x^2 e_p \otimes x^{l-1}\ar@{-}[dl]\ar@{-}[dr]&&& \\
&&&&\ddots &&\ddots &&&&\\
&&&&&x^{p-2} e_p \otimes x^{l-2}\xi  \ar@{-}[dr]&&x^{p-2} e_p \otimes x^{l-1} \ar@{-}[dl]\ar@{-}[dr]&\\
&&&&&& x^{p-2} e_p \otimes x^{l-2}\xi  && x^{p-1} e_p \otimes x^{l-1}.}

for $\mathbb{H}(\fc^\sigma\otimes_{\fc^0}  \fc)e_l$ for $l \geq 2$,
where we have used Lemma \ref{xxicommute} to move all the $\xi$ to the rear.

To see the twist by $\tau$,
note that internal multiplication by $x \otimes \xi$ can be identified with the negative of internal multiplication by $\xi \otimes x$,
by Lemma \ref{xxicommute}.
\endproof

It is easy to see that $\xi \otimes \xi \in \fc^\sigma \otimes \fc$
lies in the kernel of the differential, but not the image, hence contributes to homology;
this factor of homology has zero intersection with the factor of $(\fc e_p)^\sigma \otimes e_1 \fc$ described above,
since the components $\xi e_h \otimes e_{p+1-h} \xi$ of $\xi \otimes \xi$ in homology all belong to a subquotient of
$(\fc e_h)^\sigma \otimes e_{p+1-h} \fc$ for $1 \leq h \leq p-1$.
These components $\xi e_h \otimes e_{p+1-h} \xi$ of $\xi \otimes \xi$ for $1 \leq h \leq p-1$, shifted by $-1$ in $j$-degree,
thus give us the factor $\overline{\fc^0}$ of homology.

\medskip {\bf The case $i=-2$.}
Following the super sign convention,
the left differential on the dg $\fc$-$\fc$-bimodule $\fc^\sigma \otimes_{\fc^0} \fc^\sigma \otimes_{\fc^0} \fc$ is given by
$$\ttd(a \otimes b \otimes c) =
(-1)^{|a|_k} (ax \otimes \xi b \otimes c + a\xi \otimes x b \otimes c) +
(-1)^{|a|_k + |b|_k} (a \otimes bx \otimes \xi c + a \otimes b\xi \otimes xc).$$
for $a,b \in \fc^\sigma$, $c \in \fc$.

\begin{lem}\label{welement}
$w =(\xi \otimes \xi \otimes 1 - 1 \otimes \xi \otimes \xi) \in
\fc^\sigma \otimes \fc^\sigma \otimes \fc$
lies in the kernel of the differential.
\end{lem}
\proof
This is a straightforward computation. There is no nontrivial super-commutation involved since $\xi$ has $k$-degree zero.
\endproof

The two dimensional subspace
$e_p \otimes (e_1 \fc e_p)^\sigma \otimes e_1$ of $\fc^\sigma \otimes \fc^\sigma \otimes \fc$ which has basis given by
$\{ e_p \otimes x^{p-1} \otimes e_1, e_p \otimes x^{p-2}\xi \otimes e_1 \}$ maps to zero under the differential,
and represents the two composition factors of $\fc^\sigma \otimes \fc^\sigma \otimes \fc$ isomorphic to $p \otimes 1^{op}$ (where again $i$ denotes the simple at vertex $i$).
Since as an ungraded bimodule, $M$ has a two dimensional top isomorphic to $(p \otimes 1^{op})^{\oplus 2}$,
we conclude that the  factor $M$ of homology is generated as a bimodule by
$\{ e_p \otimes x^{p-1} \otimes e_1, e_p \otimes x^{p-2}\xi \otimes e_1 \}$.

We know that $w =(\xi \otimes \xi \otimes 1 - 1 \otimes \xi \otimes \xi) \in
\fc^\sigma \otimes_{\fc^0} \fc^\sigma \otimes_{\fc^0} \fc$
lies in the kernel of the differential but not the image, hence contributes to homology;
this factor of homology has zero intersection with the factor of $(\fc e_p)^\sigma \otimes (e_1 \fc e_p)^\sigma \otimes e_1 \fc$ described above,
since the components $\fc \otimes \xi \otimes \fc$ and $\fc^\sigma \otimes \langle x^{p-1}, x^{p-2}\xi \rangle \otimes \fc$
have different degrees in the middle tensor.
The elements $e_lw e_l$ must be generators of the factor $_\fc\fc_\fc$, for $1 \leq l \leq p$
because $e_l Fw e_l \cong  l \otimes l^{op}$ as a $\fc^0$-$\fc^0$-bimodule,
and all composition factors of the regular bimodule for $\fc$ outside the top are isomorphic to $l \otimes m^{op}$ for $l \neq m$.

\begin{lem}\label{twoximove}
We have
$\xi \otimes \xi \otimes x^{p-1} = (-1)^{p} x^{p-1} \otimes \xi \otimes \xi$
in $\mathbb{H}(\fc^\sigma \otimes_{\fc^0}  \fc^\sigma \otimes_{\fc^0} \fc)$.
\end{lem}
\proof We get the first element by multiplying $w =(\xi \otimes \xi \otimes 1 - 1 \otimes \xi \otimes \xi)$
on the right by $x^{p-1}$, and the second by multiplying $w$
on the left by $(-1)^{p-1}x^{p-1}$. We can show explicitly that $xw = -wx$ in homology, since
$$wx = \xi \otimes \xi \otimes x - 1 \otimes \xi \otimes \xi x, \quad -xw =  x \otimes \xi \otimes \xi - x\xi \otimes \xi \otimes 1,$$
and both are equal to $\xi \otimes \xi \otimes x + \xi \otimes x^* \otimes \xi + x \otimes \xi \otimes \xi$
in homology thanks to the following computations of images under the differential acting on the left:
$$\ttd: \xi \otimes 1 \otimes 1 \mapsto (\xi \otimes \xi \otimes x + \xi \otimes x \otimes \xi + \xi x \otimes \xi \otimes 1),$$
$$\ttd: 1 \otimes 1 \otimes \xi \mapsto (\xi \otimes x \otimes \xi + x \otimes \xi \otimes \xi + 1 \otimes \xi \otimes x \xi).$$
The result follows.
\endproof

\begin{lem}
The component of $\mathbb{H}(\fc^\sigma \otimes_{\fc^0}  \fc^\sigma \otimes_{\fc^0} \fc)$
generated by $w$ is quasi-isomorphic to $\fc$ as a $\fc$-$\fc$-bimodule.
\end{lem}
\proof
We have $xw = - wx$ and $\xi w = - w \xi$. So the component of homology is given by $\fc^\zeta$,
where $\zeta$ is the automorphism that sends $\xi$ to $-\xi$ and $x$ to $-x$. This automorphism is inner,
being given by conjugation by $\sum (-1)^le_l$. Consequently $\fc^\zeta \cong \fc$ as $\fc$-$\fc$-bimodules.
\endproof

\begin{lem} \label{notwistbytau}
The component of $\mathbb{H}(\fc^\sigma \otimes_{\fc^0}  \fc^\sigma \otimes_{\fc^0} \fc)$
generated by $e_p \otimes x^{p-1} \otimes e_1$ and $e_p \otimes x^{p-2} \xi \otimes e_1$ is isomorphic to $M$
as a $\fc$-$\fc$-bimodule.
\end{lem}
\proof
We have
$$\ttd(x e_p \otimes e_1 x^{p-2} \otimes e_2) = (-1)^{p-1}(x \otimes x^{p-1} \otimes \xi + x \otimes x^{p-2} \xi \otimes x),$$
$$\ttd(e_{p-1} \otimes e_2 x^{p-2} \otimes e_1 x) = ( x \otimes \xi x^{p-2} \otimes x + \xi \otimes x^{p-1} \otimes x).$$
Thus
$x \otimes x^{p-1} \otimes \xi = - x \otimes x^{p-2}\xi \otimes x = \xi \otimes x^{p-1} \otimes x$,
so internal multiplication  by $x \otimes \xi - \xi \otimes x$ in our copy of $M$ is zero,
which means the bimodule $M$ is untwisted.
\endproof

\medskip {\bf The case $i \leq -3$.}

We now assume $i \leq -3$, and write a basis element in $\mathbb{H}((\fc^\sigma \otimes_{\fc^0}  \fc)^{\otimes -i})$ as an $-i+1$-fold tensor.

\begin{thm}\label{xandys}
\begin{enumerate}[(i)]
 \item For $i$ even, we set
\begin{equation*}\begin{split}\x_{f,-i} &= e_p \otimes (\xi \otimes \xi)^{\otimes f} \otimes (x^{p-1})^{\otimes -i-1-2f}\otimes e_1 \\&
\equiv \pm e_p \otimes (x^{p-1})^{\otimes  -i-1-2f} \otimes (\xi \otimes \xi)^{\otimes f}\otimes e_1 \\&\not \equiv 0  \in \mathbb{H}((\fc^\sigma\otimes_{\fc^0}  \fc)^{\otimes -i}) \end{split}
\end{equation*}
and
\begin{equation*}\begin{split}\y_{f,-i} &= e_p \otimes (\xi \otimes \xi)^{\otimes f} \otimes (x^{p-1})^{\otimes -i-2-2f} \otimes x^{p-2} \xi \otimes e_1\\&\equiv \pm e_p \otimes  x^{p-2} \xi\otimes (x^{p-1})^{\otimes  -i-2-2f} \otimes (\xi \otimes \xi)^{\otimes f}\otimes e_1\\&\not \equiv 0 \in \mathbb{H}((\fc^\sigma\otimes_{\fc^0}  \fc)^{\otimes -i}) \end{split}
\end{equation*}
for $0 \leq f \leq \frac{-i-2}{2}$. %i-1 (odd) nontrivial entries and at most -i-2 (even) of them \xi

\item For $i$ odd, we set
\begin{equation*}\begin{split}
     \x_{f,-i} &= e_p \otimes (\xi \otimes \xi)^{\otimes f} \otimes (x^{p-1})^{\otimes -i-1-2f} \otimes e_1\\ &\equiv \pm e_p \otimes (x^{p-1})^{\otimes -i-1-2f} \otimes (\xi \otimes \xi)^{\otimes f} \\&\not \equiv 0     \in \mathbb{H}((\fc^\sigma\otimes_{\fc^0}  \fc)^{\otimes -i})
\end{split}
\end{equation*}

for $0 \leq f \leq \frac{-i-1}{2}$,and
\begin{equation*}\begin{split}\y_{f,-i}& = e_p \otimes (\xi \otimes \xi)^{\otimes f} \otimes (x^{p-1})^{\otimes -i-2 -2f} \otimes x^{p-2} \xi\otimes e_1\\& \equiv \pm e_p \otimes  x^{p-2} \xi\otimes (x^{p-1})^{\otimes -i-2-2f} \otimes (\xi \otimes \xi)^{\otimes f}\\& \not \equiv 0 \in \mathbb{H}((\fc^\sigma\otimes_{\fc^0}  \fc)^{\otimes -i}) \end{split}
\end{equation*}
for $0 \leq f \leq \frac{-i-3}{2}$.
\end{enumerate}
For  $f < \frac{-i-2}{2}$, the elements % i.e. in all cases for i even and all except f = (-i-1)/2 for i odd
$\x_{f,-i}$ and $\y_{f,-i}$  generate the possibly twisted copy of $M$ in the corresponding homological degree as a bimodule.
For $f = \frac{-i-1}{2}$ (and $i$ odd), $\x_{f,-i}$ generates a possibly twisted copy of $\overline{M}$.
\end{thm}

\proof In all cases the equivalence in homology of the two given representatives follows from Lemma \ref{twoximove}.
%The given elements are certainly killed by the differential, so it suffice to prove they are not in its image. This follows from a similar computation as that in the proof of Lemma \ref{twoximove}.
We proceed by induction on $i$. Cases $i=-1$ and $i=-2$ have been examined above.
So assume the statement is true for $-i+1$,
i.e. $\x_{f,-i+1}$ and $\y_{f,-i+1}$ are nonzero in $\mathbb{H}((\fc^\sigma\otimes_{\fc^0}  \fc)^{\otimes -i+1})$ for $f \leq \frac{-i-3}{2}$
and generate copies of $M$ as a bimodule.
In particular, this means that $\x_{f,-i+1} x^{p-1}$ and $\y_{f,-i+1}x^{p-1}$% = e_p \otimes (\xi \otimes \xi)^{\otimes f} \otimes (x^{p-1})^{\otimes i-3 -2f} \otimes x^{p-2} \xi\otimes x^{p-1} \equiv $
are nonzero in $\mathbb{H}((\fc^\sigma\otimes_{\fc^0}  \fc)^{\otimes -i+1})$. But since $\mathbb{H}((\fc^\sigma\otimes_{\fc^0}  \fc)^{\otimes -i})e_1  = \mathbb{H}((\fc^\sigma\otimes_{\fc^0}  \fc)^{\otimes -i+1})e_p \otimes e_1 $, this means that
$$\x_{f,-i+1} x^{p-1} \otimes e_1 = \x_{f,-i}$$ and
\begin{equation*}\begin{split}
\y_{f,-i+1}x^{p-1} \otimes e_1  &= e_p \otimes (\xi \otimes \xi)^{\otimes f} \otimes (x^{p-1})^{\otimes -i-3 -2f} \otimes x^{p-2} \xi\otimes x^{p-1} \otimes e_1 \\ &\equiv  e_p \otimes (\xi \otimes \xi)^{\otimes f} \otimes (x^{p-1})^{\otimes -i-2 -2f} \otimes x^{p-2} \xi\otimes e_1 \\&= \y_{f,-i} \end{split}
\end{equation*}
are non-zero in $\mathbb{H}((\fc^\sigma\otimes_{\fc^0}  \fc)^{\otimes -i})$.
Since $e_pMe_1$ is two dimensional, we see that they must be the generators of the corresponding copy of $M$.
It remains to consider the cases $i$ even, $f =  \frac{-i-2}{2}$ and $i$ odd, $f =  \frac{-i-1}{2}$.
For the first case, we know by the inductive assumption that $$x_{f,-i+1} = e_p \otimes (\xi \otimes \xi)^{\otimes f} \otimes e_1$$
is non-zero in homology and generates a copy of $\overline{M}$.
In particular this means that $x_{f,-i+1}x^{p-1}$ and $x_{f,-i+1}x^{p-2}\xi$ are non-zero and, arguing as above, wee see
that $x_{f,-i+1}x^{p-1}\otimes e_1 = \x_{f,-i}$ and $x_{f,-i+1}x^{p-2}\xi \otimes e_1 = \y_{f,-i}$ are nonzero in
$\mathbb{H}((\fc^\sigma\otimes_{c^0} \fc)^{\otimes -i})$ and because $e_pMe_1$ is two dimensional, generate a copy of $M$.
Now consider the case where $i$ is odd and $f =  \frac{-i-1}{2}$,
i.e. consider the element $\x_{\frac{-i-1}{2}, -i} = e_p \otimes (\xi \otimes \xi)^{\frac{-i-1}{2}} \otimes e_1$.
This is obviously in the kernel, but not in the image of the differential and hence non-zero in
$\mathbb{H}((\fc^\sigma\otimes_{\fc^0}  \fc)^{\otimes -i})$.
Since $e_p \overline{M} e_1$ is one dimensional,
it must be the generator of the corresponding copy of $\mathbb{H}(\fc^\sigma\otimes_{\fc^0}  \fc)$.
\endproof

The remaining homology is given as follows:

\begin{lem}\label{lowest}
\begin{enumerate}[(i)]
\item For $i$ even, the elements $e_l w^{\frac{-i}{2}} e_l $  (for $1 \leq l \leq p$) generate, as a $\fc$-$\fc$-bimodule the factor isomorphic to $\fc$ in $\mathbb{H}((\fc^\sigma \otimes_{\fc^0}  \fc)^{\otimes -i})$.
\item For $i$ odd, the elements $e_l (\xi \otimes \xi)^{\otimes \frac{1-i}{2}} e_l$ for $1 \leq l \leq p-1$ generate the bimodule isomorphic to $\overline{\fc^{0\sigma}}$ in $\mathbb{H}((\fc^\sigma\otimes_{\fc^0}  \fc)^{\otimes -i})$.
 \end{enumerate}
\end{lem}

\proof
For $i$ even, the elements $e_l w^{\frac{-i}{2}} e_l$ for $1 \leq l \leq p$ are certainly in the kernel of the differential,
and the elements are not in the image of differential, for lack of tensor factors $x$. They must be generators,
because $e_l (Fw^{\frac{-i}{2}}) e_l \cong  l \otimes l^{op}$ as an ungraded $\fc^0$-$\fc^0$-bimodule,
and all composition factors of the regular bimodule for $\fc$ outside the top are isomorphic to $l \otimes m^{op}$ for $l \neq m$.

For $i$ odd, the given elements have $\xi$ in every tensor factor, hence are non-zero in
$\mathbb{H}((\fc^\sigma\otimes_{\fc^0}  \fc)^{\otimes -i})$.
They generate the desired semi-simple bimodule since that is all the homology not accounted for by factors isomorphic to $M$ and $\overline{M}$.
\endproof

\subsection{The $\fc$-$\fc$-bimodule structure of $\Upsilon^{\leq 1}$.} \label{explicit}

Suppose $p \geq 3$.
Due to our analysis in Lemma \ref{trunc} we do not care about the entire algebra structure of $\Upsilon$:
all we want to know about the product is what it looks like when composed with projection onto the subspace
$\ft \mathbb{T}_\fc(\ft^{-1})$, which we denote by $\Upsilon^{\leq 1}$.

\begin{lem} As $j$-graded $\fc$-$\fc$-bimodules, $\Upsilon^-$ can be identified with the homology of
$$\xymatrix@C=1pt@R=2pt{
    &        &        &             & \fc &          \\
    &        &        & \overline{M}^\tau \ar[rr]^0 \langle 1 \rangle &      &  \overline{\fc^{0\sigma}} \langle 1 \rangle \\
    &        &  M  \langle -p+1 \rangle  \ar[rr]^0 &      & \fc &      \\
    &   M^\tau \langle -2p+1 \rangle \ar[rr]^0 & & \overline{M}^\tau \langle 1 \rangle \ar[rr]^0  &      &  \overline{\fc^{0\sigma}} \langle 1 \rangle \\
M  \langle -3p+1 \rangle \ar[rr]^0 &  &   M \langle -p+1 \rangle  \ar[rr]^0 &       & \fc &         \\
    &        &        &             &......&          \\
}$$
As $k$-graded $\fc$-$\fc$-bimodules, $\Upsilon^-$ can be identified with the homology of
$$\xymatrix@C=1pt@R=2pt{
    &        &        &             & \fc &          \\
    &        &        & \overline{M}^\tau \ar^0[rr][0] &      &  \overline{\fc^{0\sigma}}[-1] \\
    &        &  M [p-1]  \ar[rr]^0 &      & \fc[-1] &      \\
    &   M^\tau [2p-2] \ar[rr]^0 &  & \overline{M}^\tau [1] \ar[rr]^0  &      &  \overline{\fc^{0\sigma}}[-2] \\
M [3p-3] \ar[rr]^0 & &   M [p-2]  \ar[rr]^0 &      & \fc[-2] &         \\
    &        &        &             &......&          \\
}$$
% As an ungraded $\fc$-$\fc$-bimodule $\Upsilon^+$ can be identified with the homology of
% $$\xymatrix@C=1pt@R=2pt{
% &             &......&          &        &        &     \\
% &             &\fc^* \ar[rr]^0 &    & M \ar[rr]^0 & & M \\
% &\fc^{0\sigma} \ar[rr] &  &   M^\tau  \ar[rr]^0 & & M^\tau  & \\
% &             &\fc^* \ar[rr]^0 &   &   M   &        &     \\
% &\fc^{0\sigma} \ar[rr] &      &    M^\tau   &        &        &     \\
% &             &\fc^* &          &        &        &     \\
% &   \fc^{0\sigma}   &      &          &        &        &     \\
% & & \fc & & & & \\
% }$$

\end{lem}
\proof
As in Proposition \ref{dgminus} rows correspond to $\ft^{i}$, for $i \in \mathbb{Z}_{\leq 0}$. The Lemma is obtained from Proposition \ref{dgminus} by applying homology, using Remark \ref{truncatemaltese}, and noting that there is no nonzero map between $\overline{M}^\tau$ and $\overline{\fc^{0\sigma}}$, since as left $\fc$-module, all simples in the top of $\overline{M}^\tau$ are isomorphic to the simple at the vertex $p$, and $\overline{\fc^{0 \sigma}}$ is preciesly a direct sum of the other simples.
It remains to check that $\fc$ splits off from $M$ in homology in $\ft^{-2}$. To see this, recall that it is generated as a bimodule by $\xi \otimes \xi \otimes 1 - 1 \otimes \xi \otimes \xi$,
whereas the copy of $M$ is generated as a bimodule by
$e_p \otimes x^{p-1} \otimes e_1$ and $e_p \otimes x^{p-2} \xi \otimes e_1$.
Since $p \geq 3$ the middle term in the tensor generating $\fc$ has different radical degree
from the middle term in the tensors generating $M$.
As images of monomials in $x$ and $\xi$ under the differential on $\ft^{-2}$ have middle terms with identical radical degrees,
whilst multiplying such a monomial on the left or right by a monomial in $\fc$ does not alter the radical degree of the middle term,
there is a splitting as required.
\endproof

Rewriting the above expressions we see that
as $j$-graded $\fc^0$-$\fc^0$-bimodules, $\Upsilon^{\leq 1}$ can be identified with
$$\xymatrix@C=1pt@R=2pt{
    &        &        &   & & \fc^{0 \sigma}  \langle 1 \rangle      \\
    &        &        &             & \fc &          \\
    &        &        & \overline{M}^\tau \langle 1 \rangle  & \oplus    &  \overline{\fc^{0\sigma}} \langle 1 \rangle \\
    &        &  M  \langle -p+1 \rangle  &  \oplus    & \fc &      \\
    &   M^\tau \langle -2p+1 \rangle & \oplus & \overline{M}^\tau \langle 1 \rangle  & \oplus  &  \overline{\fc^{0\sigma}} \langle 1 \rangle \\
M  \langle -3p+1 \rangle & \oplus &   M \langle -p+1 \rangle & \oplus      & \fc &         \\
    &        &        &             &......&          \\
}$$
whilst as $k$-graded $\fc^0$-$\fc^0$-bimodules, $\Upsilon^{\leq 1}$ can be identified with
$$\xymatrix@C=1pt@R=2pt{
    &        &        &   & & \fc^{0 \sigma}      \\
    &        &        &             & \fc &          \\
    &        &        & \overline{M}^\tau  & \oplus    &  \overline{\fc^{0\sigma}} \\
    &        &  M  [p-1]  &  \oplus    & \fc &      \\
    &   M^\tau [2p-2] & \oplus & \overline{M}^\tau  & \oplus  &  \overline{\fc^{0\sigma}} \\
M  [3p-3] & \oplus &   M [p-1]  & \oplus        & \fc &         \\
    &        &        &             &......&          \\
}$$

\begin{thm} \label{dgalgebraandhomology}
The algebra $\Upsilon^-$ is isomorphic to the algebra given by
$$\mathbb{T}_{\fc}(\overline{M} \langle 1 \rangle \ttu \oplus \overline{\fc^{0 \sigma}} \langle 1 \rangle \ttv) \otimes_F F[\ttw ]$$
modulo the relations $\ttv^2 = 0$, $\ttu\ttv= -\ttv\ttu =\ttw $, the relations $x\ttu = -\ttu x$, $x\ttv = -\ttv x$, $\xi \ttu = \ttu \xi$, $\xi \ttv = \ttv \xi$ for $x, \xi \in \fc$,
and relations that ensure products on the generators $\overline{M}$ and $\overline{\fc^{0 \sigma}}$ are given by the maps
$$\gamma: \overline{M}  \otimes \overline{\fc^{0\sigma}} \rightarrow \fc \ttw ,
\quad \gamma: \overline{\fc^{0 \sigma}}  \otimes \overline{M}  \rightarrow \fc \ttw $$
defined in Lemma \ref{mapsgamma}; here $\ttu$, $\ttv$, and $\ttw $ are formal variables.
\end{thm}
\proof
Again, to simplify notation, we omit shifts in $j$-degree.
It is clear from the defining relations, and properties of the bimodules,
that the algebra defined by generators and relations is a quotient of the space
$$\xymatrix@C=5pt@R=4pt{
     &        &        &                        & \fc     &                       &        &        &     \\
     &        &        & \overline{M}^\tau \ttu  &  \oplus        &  \overline{\fc^{0\sigma}} \ttv  &        &        &     \\
     &        &  M\ttu^2[p-1]  & \oplus                 & \fc \ttw    &                       &        &        &     \\
     &   M^\tau \ttu^3[2p-2] & \oplus & \overline{M}^\tau \ttu \ttw   &   \oplus       &  \overline{\fc^{0\sigma}} \ttv \ttw  &        &        &     \\
M\ttu^4[3p-3] & \oplus &  M\ttu^2\ttw [p-1] & \oplus                 & \fc \ttw^2 &                       &        &        &     \\
    &         &        &                        &  ......  &                       &        &        &     \\
}$$
as a $\fc$-$\fc$-bimodule. Thanks to the relations between $\ttu$, $x$, and $\xi$,
we have $M \ttu^i \cong M^{\tau^i}$ as $\fc$-$\fc$-bimodules.
We now observe the existence of an algebra homomorphism from the algebra defined above to $\Upsilon^-$ which identifies
$\overline{M} \ttu \oplus \overline{\fc^{0\sigma}} \ttv $ with $\ft^{-1}$,
identifies $\ttv $ with $(\xi \otimes \xi)$,
identifies $M \ttu^2$ with the component $M$ of $H(\ft^{-2})$,
and identifies $\ttw$ with the element$w=\xi \otimes \xi \otimes 1 - 1 \otimes \xi \otimes \xi \in \ft^{-2}$ from Lemma \ref{welement}.
Note the relation $\ttu\ttv  = -\ttv \ttu$ comes from the fact that the sign preceding $1 \otimes \xi \otimes \xi$ in $\xi \otimes \xi \otimes 1 - 1 \otimes \xi \otimes \xi$ is the negative of the sign preceding $\xi \otimes \xi \otimes 1$.
The preceding lemma implies that we have an algebra isomorphism from the homology of this algebra to $\Upsilon^-$.
\endproof

\medskip 

We can extend the above to include multiplication by the subspace $\mathbb{H}(\ft)$,
which we identify with $\fc^{0 \sigma} \langle 1 \rangle[0] \ttz$,
where $\ttz$ is a formal central variable. The products
$$M \ttu^a \ttw^b \otimes \fc^{0 \sigma} \ttz \rightarrow M \ttu^{a-1} \ttw^b, \quad \fc^{0 \sigma} \ttz \otimes M \ttu^a \ttw^b \rightarrow M \ttu^{a-1} \ttw^b,$$
$$M \ttu \ttw^b \otimes \fc^{0 \sigma} \ttz \rightarrow \fc \ttw^b, \quad \fc^{0 \sigma} \ttz \otimes M \ttu \ttw^b \rightarrow \fc \ttw^b,$$
are given by maps $\gamma$; the products
$$\fc \ttw^b \otimes \fc^{0 \sigma} \ttz \rightarrow \overline{\fc^{0\sigma}} \ttv\ttw^{b-1}, \quad \fc^{0 \sigma} \ttz \otimes \fc \ttw^b \rightarrow \overline{\fc^{0 \sigma}}\ttv\ttw^{b-1}$$
are given by multiplication followed by projection; and the product of $\overline{\fc^{0 \sigma}} \ttz$ and $\fc^{0 \sigma} \ttv$ is zero.

\medskip 

In this section, we have described the product on $\Upsilon^{\leq 1}$ in its entirety: the space is isomorphic to
$$\xymatrix@C=1pt@R=2pt{
    &        &        &   & & \fc^{0 \sigma}        \\
    &        &        &             & \fc &          \\
    &        &        & \overline{M}^\tau   & \oplus    &  \overline{\fc^{0\sigma}} \\
    &        &  M    &  \oplus    & \fc &      \\
    &   M^\tau  & \oplus & \overline{M}^\tau   & \oplus  &  \overline{\fc^{0\sigma}}  \\
M  & \oplus &   M  & \oplus      & \fc &         \\
    &        &        &             &......&          \\
}$$
and all products on these components are described by the natural maps
$M \otimes M \rightarrow M$, $M \otimes \fc \rightarrow M$, $\fc \otimes M \rightarrow M$,
and $\pm \gamma$;
we use here the natural bimodule isomorphism $^\tau M^\tau \cong M$
that sends $x^l \otimes x^m \in {^\tau L_l} \otimes L_r^\tau$ to $(-1)^{l+m} x^l \otimes x^m$.
We describe a polytopal version of this in section \ref{polybasis}.

\begin{example}\label{exagain}
Let $p=3$. 
Then $\fc_3$ is the subalgebra of the algebra computed in Section \ref{example} given by the first three idempotents.
In this example we illustrate that computing $\fO_{ F} \fO_{\Upsilon}^2 (F[z,z^{-1}]) $ does indeed give the algebra defined in Section \ref{example} by computing the last (and largest) projective.
In the first step, we just obtain $\fO_{\Upsilon} (F[z,z^{-1}])  = \Upsilon$. As we apply $\fO_{ F}$ in the end we need to concentrate only on the part of $i$-dgree $0$, i.e. $\fc$. Now to obtain the $9$-th projective, we identify $9 =(3-1)3+3$ with its $3$-adic expansion $(3,3)$. Our algorithm then tells us that we need look at the third projective in $\fc$, which has $jk$-graded Loewy series
{\footnotesize 
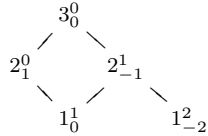
$$\xymatrix@=5pt{ & 3^0_0\ar@{-}[dl]\ar@{-}[dr] & & \\2^0_1 \ar@{-}[dr]&& 2^1_{-1}\ar@{-}[dl]\ar@{-}[dr] & \\&1^1_0 &&1^2_{-2} }$$}
so by tensoring the quotient $3^0_0$ with the third projective of $\fc$, we should get a quotient 
{\footnotesize 
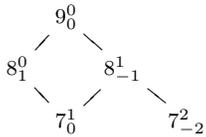
$$\xymatrix@=5pt{
& 9^0_0\ar@{-}[dl]\ar@{-}[dr] & & \\8^0_1 \ar@{-}[dr]&& 8^1_{-1}\ar@{-}[dl]\ar@{-}[dr] & \\&7^1_0 &&7^2_{-2}.
}$$}
The simple subquotient $2^0_1$ should be tensored with $\bbH(\ft)e_3 = 1^0_1$, to give a simple subquotient $(2,1)^0_2$ which we identify with $4^0_2$.
The simple $2^1_{-1}$ should be tensored with $\bbH(\ft^{-1})e_3$ which is the direct sum of a simple $1^0_1$ (from $\overline{\fc^{0\sigma}}\langle 1\rangle$) and the indecomposable left summand of $\overline{M}\tau e_3 \langle 1\rangle= M\tau e_3\langle 1\rangle$ (the truncation is only seen on $Me_1$) which has Loewy filtration  
{\footnotesize
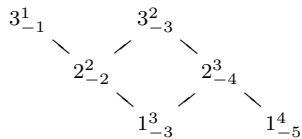
$$\xymatrix@=5pt{ 3^1_{-1}\ar@{-}[dr] && 3^2_{-3}\ar@{-}[dl]\ar@{-}[dr] & & \\&2^2_{-2} \ar@{-}[dr]&& 2^3_{-4}\ar@{-}[dl]\ar@{-}[dr] & \\&&1^3_{-3} &&1^4_{-5} }$$}
to obtain a subquotient (again identifying $(2,1),(2,2),(2,3)$ with $4,5,6$ respectively)
{\footnotesize
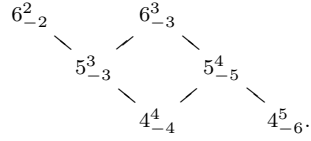
$$\xymatrix@=5pt{ 6^2_{-2}\ar@{-}[dr] && 6^3_{-3}\ar@{-}[dl]\ar@{-}[dr] & & \\&5^3_{-3} \ar@{-}[dr]&& 5^4_{-5}\ar@{-}[dl]\ar@{-}[dr] & \\&&4^4_{-4} &&4^5_{-6}. }$$}
The simple $1^1_0$ should be tensored with $\fc e_3$ to obtain 
{\footnotesize 
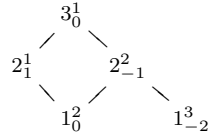
$$\xymatrix@=5pt{ & 3^1_0\ar@{-}[dl]\ar@{-}[dr] & & \\2^1_1 \ar@{-}[dr]&& 2^2_{-1}\ar@{-}[dl]\ar@{-}[dr] & \\&1^2_0 &&1^3_{-2} }$$}
and finally the simple $1^2_{-2}$ should be tensored with  $\bbH(\ft^{-2})e_3$, which is isomorphic to the direct sum of $\fc e_3$ giving a submodule
{\footnotesize 
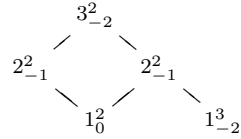
$$\xymatrix@=5pt{ & 3^2_{-2}\ar@{-}[dl]\ar@{-}[dr] & & \\2^2_{-1} \ar@{-}[dr]&& 2^2_{-1}\ar@{-}[dl]\ar@{-}[dr] & \\&1^2_0 &&1^3_{-2} }$$}

and $Me_3\langle -2\rangle[2]$ which has Loewy filtration
{\footnotesize
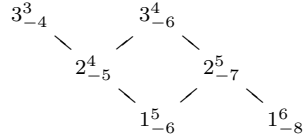
$$\xymatrix@=5pt{ 3^3_{-4}\ar@{-}[dr] && 3^4_{-6}\ar@{-}[dl]\ar@{-}[dr] & & \\&2^4_{-5} \ar@{-}[dr]&& 2^5_{-7}\ar@{-}[dl]\ar@{-}[dr] & \\&&1^5_{-6} &&1^6_{-8} }$$}
and hence, when tensored with $1^2_{-2}$, produces a submodule
{\footnotesize
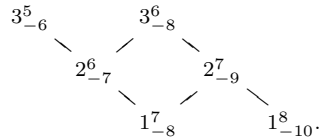
$$\xymatrix@=5pt{ 3^5_{-6}\ar@{-}[dr] && 3^6_{-8}\ar@{-}[dl]\ar@{-}[dr] & & \\&2^6_{-7} \ar@{-}[dr]&& 2^7_{-9}\ar@{-}[dl]\ar@{-}[dr] & \\&&1^7_{-8} &&1^8_{-10}. }$$}
All of these are visible in the projective given in Section \ref{example} and exhaust the module. %To obtain the explicit extensions in Section \ref{example} we have cheated and done some additional computations in the derived category of $\bfc\ml$: while, of course, the algebra structure in our general presentation is completely explicit, we do not know how to derive a general formula for quiver and relations for the extension algebra of Weyl modules.

\end{example}

\subsection{The case $p=2.$}\label{p2}

In the preceding section we assumed $p \geq 3$. Here we explain how to adapt the results to the case $p=2$.
This case is exceptional because when $p=2$ we can extend $L_l$ by a projective $\fc$-module. We do not need to worry about signs
because $1 = -1$ modulo $2$.

Let $S^h(x, \xi)$ be the collection of polynomials in $x$ and $\xi$ of degree $h$, a vector space of dimension $h+1$.
We have a natural product map $$S^{h_1}(x, \xi) \otimes S^{h_2}(x,\xi) \rightarrow S^{h_1+h_2}(x, \xi).$$

For $i < 0$, we define $V_{-i}$ to be the $\fc$-$\fc$-bimodule
$$\xymatrix@C=-20pt@R=3pt{
      & S^{-i-1}(x, \xi) & \\
S^{-i}(x,\xi)_l & \oplus & S^{-i}(x,\xi)_r \\
 & S^{-i+1}(x,\xi)&
}$$
where $x$ and $\xi$ act by left multiplication sending $S^{-i-1}(x,\xi)$ to $S^{-i}(x,\xi)_l$,
sending $S^{-i}(x,\xi)_r$ to $S^{-i+1}(x,\xi)$, sending
$S^{-i}(x,\xi)_l$ and $S^{-i+1}(x,\xi)$ to zero;
where $x$ and $\xi$ act by right multiplication sending $S^{-i-1}(x,\xi)$ to $S^{-i}(x,\xi)_r$,
sending $S^{-i}(x,\xi)_l$ to $S^{-i+1}(x,\xi)$, sending
$S^{-i}(x,\xi)_r$ and $S^{-i+1}(x,\xi)$ to zero.
We have
$$S^{-i-1}(x, \xi) = e_2 V_{-i} e_1, \quad S^{-i}(x,\xi)_l = e_1 V_{-i} e_1,$$
$$S^{-i}(x,\xi)_r = e_2 V_{-i} e_2 , \quad S^{-i+1}(x,\xi) = e_1 V_{-i} e_2$$

\begin{lem} Let $p=2$.
The homology of $\ft^{i}$ is isomorphic to $V_{-i}$, for $i < 0$.
\end{lem}
\proof
The top has a basis of elements of the form
$1 \otimes a_1 \otimes ... \otimes a_{-i-1} \otimes 1$, where $a_1,...,a_{-i-1}$ is a list of letters consisting of
$l$ $x$s and $i-1-l$ $\xi$s. The basis element is independent of the order of the elements in the list.
\endproof

The product on $\Upsilon^-$ is given by the product
on the algebra $S(x, \xi)$ whenever degrees and idempotents match up appropriately,
and zero otherwise; for example the product of $S^{-i}_l$ and $S^{-i}_r$ is zero, whilst the product
$S^{-i}_l \otimes S^{-i}_l \rightarrow S^{-2i}_l$ is given by the natural product on $S(x,\xi)$.

\medskip 

The action of $\Upsilon^1 = \langle e_1 \otimes e_2, e_2 \otimes e_1 \rangle$
on $\Upsilon^i$ for $i \leq 1$ is given here by the maps
$$\Upsilon^1 \otimes \Upsilon^i \rightarrow \Upsilon^{i+1}, \quad \Upsilon^i \otimes \Upsilon^1 \rightarrow \Upsilon^{i+1}$$
which are the identity on a component $S^h(x, \xi) \subset \Upsilon^i, \Upsilon^{i+1}$
whenever idempotents match up appropriately, and zero otherwise.

\subsection{Polytopal basis for $\Upsilon^{\leq 1}$. \label{polybasis}}
{\bf For a small correction to this basis, see Section \ref{indexing}. This mistake does not affect the cardinality or the $ijk$-degree, but does need to be corrected for multiplication to be correct as stated.}

We first define \emph{monomials} in $\Upsilon$.
These form a collection of elements each of which is determined by its $jk$-degree and, upon multiplying on the left and right by idempotents $e_s$, gives an element of a basis of $\Upsilon$, or zero.

We define a monomial in $\Upsilon^0$ to be an element of $\Upsilon^0 = \fc$ with representative $x^l\xi$, for some $l$.
We define a monomial in $\Upsilon^i$, for $i \leq -1$, to be an element of $\Upsilon^i = \mathbb{H}(\ft^{\otimes i})$
with representative $x^l\xi w^{\frac{-i}{2}}$ or $x^l w^{\frac{-i}{2}}$for $i$ even and $e_l\xi^{\otimes 1-i}e_l$ for $i$ odd (see Lemma \ref{lowest}) or
$x^m \x_{f,-i} x^l$, $x^m \y_{f,-i} x^l$ (see Theorem \ref{xandys}).
We say there is a unique monomial in $\Upsilon^1$, namely $1$.

\begin{thm}
$$\dim e_s \Upsilon^{ijk} e_t \leq 1$$ for all $i \in \mathbb{Z}_{\leq 1}$, $j,k \in \mathbb{Z}$, $1 \leq s,t \leq p$.
\end{thm}
\proof
Suppose $p \geq 3$.
$\Upsilon^i$ is the homology of $\ft^i$, for $i \in \mathbb{Z}_{\leq 0}$, which is the homology of
$$M  \langle (i+1)p+1 \rangle [(i+1)(p-1)] \oplus M  \langle (i+3)p+1 \rangle [(i+3)(p-1) -1] \oplus $$
$$...\oplus  M \langle -p+1 \rangle[p+\frac{i}{2}]   \oplus    \fc[\frac{i}{2}]$$
for $i$ even, and the homology of
$$M^\tau  \langle (i+1)p+1 \rangle [(i+1)(p-1)] \oplus M^\tau  \langle (i+3)p+1 \rangle [(i-1)(p-1)] \oplus $$
$$...\oplus M^\tau \langle 1 \rangle[\frac{i+1}{2}] \rightarrow \fc^{0\sigma} \langle 1 \rangle[\frac{i-1}{2}]$$
for $i$ odd. In either case, the number of $\xi$s appearing in any monomial in $\Upsilon$ determines which
bimodule factor $M$, $\fc$, $\fc^{0 \sigma}$ that monomial appears in.
The number of $\xi$s appearing in a monomial can be computed from the $jk$-degree of that monomial,
since $\xi \in \fc$ has $jk$-degree $(1,0)$ whilst $x$ has $jk$-degree $(-1,1)$.
Therefore the $ijk$-degree of a monomial determines which component $\Upsilon^i$ and
which bimodule factor $M[l]$, $\fc[0]$, $\fc^{0 \sigma}$ of that component that monomial appears in.
To establish that $\dim e_s \Upsilon^{ijk} e_t \leq 1$,
it is enough for us to observe that a monomial element of $e_sMe_t$ or $e_s \fc e_t$
is determined by its $jk$-degree, which is the case since it is determined by the number of $x$s and $\xi$s appearing.

The case $p=2$ is similar.
\endproof

\begin{lem}
The polytope $\mathcal{P}_\fc$ corresponding to $\fc$ is the set of elements $(s,j,k,t) \in \mathbb{Z}^4$ such that
$1 \leq s \leq t \leq p$, $0 \leq j+k \leq 1$, $t-s = j+2k$ and that $k=j=0$ if $s=t$.

The polytope $\mathcal{P}_0$ corresponding to $\overline{\fc^{0 \sigma}}$ is the set of elements $(s,j,k,t)$ in $\mathbb{Z}^4$ such that $t+s=p+1$,
$1 \leq s, t \leq p$, minus the element $(p,0,0,1)$.
and $j=k=0$.

The polytope $\mathcal{P}_M$ corresponding to $M$ is the set of elements $(s,j,k,t)$ in $\mathbb{Z}^4$ such that
$1 \leq s,t \leq p$, $j+2k+2=t-1-s+p$, $0 \leq j+k+2 \leq 1$.

The polytope $\mathcal{P}_{\overline{M}}$ corresponding to $\overline{M}$ is the set of elements $\mathcal{P}_{M}$ minus the element $(p,0,-1,1)$.
\end{lem}
\proof
We pick out the number of $\xi$s in an element of $M$  with $j+k+2$, and the number of $x$s with $k+2$;
the total number of $x$s and $\xi$s is thus $j+2k+4$;
the restrictions are that the number of $\xi$s is $0$ or $1$, the number of $x$s and $\xi$s is at least $2$, and
the number of $x$s and $\xi$s is $2+(t-1)-(s-p)$; the element $(p,0,-1,1)$ corresponds to the element $x \otimes \xi \in M$.
This gives us the description of $\mathcal{P}_M$.

To obtain our description of $\mathcal{P}_\fc$,
we similarly pick out the number of $\xi$s with $j+k$ and the total number of $x$s and $\xi$s with $j+2k$.

The description of $\mathcal{P}_0$ comes since nonzero elements of $e_s \fc^{0 \sigma} e_t$ satisfy the constraint $t+s = p+1$.
\endproof
\begin{example}\label{Mbasis}
The following is a diagram of the polytope for $M$ in case $p=3$ (we depict its structure as a left module):
{\footnotesize
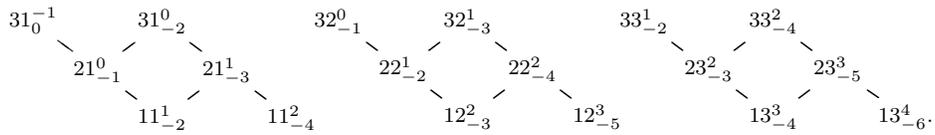
$$
\xymatrix@R=3pt@C=1pt{
31^{-1}_0 \ar@{-}[dr]&& 31_{-2}^0 \ar@{-}[dl]\ar@{-}[dr] && \\
&21_{-1}^0 \ar@{-}[dr]&& 21_{-3}^1 \ar@{-}[dl]\ar@{-}[dr]& \\
&&11_{-2}^1 && 11_{-4}^2 
}\hspace*{-2mm}\xymatrix@R=3pt@C=1pt{
        32_{-1}^0 \ar@{-}[dr]&&32_{-3}^1 \ar@{-}[dl]\ar@{-}[dr]&& \\
&22_{-2}^1 \ar@{-}[dr]&& 22_{-4}^2 \ar@{-}[dl]\ar@{-}[dr]&\\
 &&12_{-3}^2 &&12_{-5}^3 
}\hspace*{-2mm}
\xymatrix@R=3pt@C=1pt{
33_{-2}^1 \ar@{-}[dr]&&33_{-4}^2 \ar@{-}[dl]\ar@{-}[dr]& &\\
&23_{-3}^2 \ar@{-}[dr]&&23_{-5}^3 \ar@{-}[dl]\ar@{-}[dr]&\\
&&13_{-4}^3 &&13_{-6}^4.
}$$}
In the diagram an element $(s,j,k,t)$ is written $st^k_j$. Similarly a diagram of the polytope for $\fc$ in case $p=3$ is given by
{\footnotesize
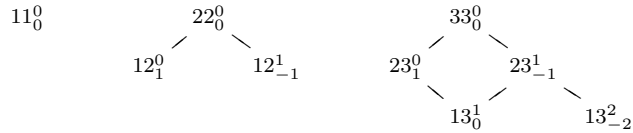
$$
\xymatrix@C=5pt@R=4pt{ 11_{0}^0 }\hspace*{10mm}
\xymatrix@C=5pt@R=4pt{ & 22_{0}^0\ar@{-}[dl]\ar@{-}[dr]  &\\ 12_{1}^0 && 12_{-1}^1}\hspace*{10mm}
\xymatrix@C=5pt@R=4pt{ &33_{0}^0 \ar@{-}[dl]\ar@{-}[dr] && \\23_{1}^0 \ar@{-}[dr]&& 23_{-1}^1 \ar@{-}[dl]\ar@{-}[dr]\\ & 13_{0}^1 &&13_{-2}^2}.
$$}

\end{example}

We introduce integers $a$ and $b$ indexing the powers of $\ttu$ and $\ttv $ appearing in a homogeneous element of $\Upsilon^{\leq 0}$, see Theorem \ref{dgalgebraandhomology}.

By Theorem \ref{dgalgebraandhomology} we have a basis for $\Upsilon^{\leq 0}$ indexed by the subset

\begin{equation*}
\begin{split}
\mathcal{P}_{\leq 0}:=&\{ v = (s,j_0,k_0,a,b,t) \in \mathbb{Z}^6 | (s,j_0,k_0,t) \in \mathcal{P}_\fc, a, b \geq 0, a=b   \}\\ \cup
&\{ v = (s,j_0,k_0,a,b,t) \in \mathbb{Z}^6 | (s,j_0,k_0,t) \in \mathcal{P}_{\overline{\fc^0}}, a, b \geq 0, a= b-1   \}\\ \cup
&\{ v = (s,j_0,k_0,a,b,t) \in \mathbb{Z}^6 | (s,j_0,k_0,t) \in \mathcal{P}_{\overline{M}}, a, b \geq 0, a=b+1   \}\\ \cup
&\{ v = (s,j_0,k_0,a,b,t) \in \mathbb{Z}^6 | (s,j_0,k_0,t) \in \mathcal{P}_{M}, a, b \geq 0, a > b+1   \}.
\end{split}
\end{equation*}

The $ijk$-degree of such an element is given by the formulas

$i = -a-b$;

$j=j_0- (a-b-1)p+1$ for $a \geq b+1$, $j=j_0$ for $a=b$, $j=j_0+1$ for $a=b-1$;

$k=k_0+ (a-b-1)(p-1)$ for $a \geq b+1$, $k=k_0$ for $a \leq b$.

We define $\mathcal{P}_{\Upsilon^{\leq 0}}$ to be the corresponding set of elements $(s,i,j,k,a,b,t)$ in $\mathbb{Z}^7$.
We define
$\mathcal{P}_{\Upsilon^{\leq 1}}$ to be $\mathcal{P}_{\Upsilon^{\leq 0}} \cup \{ (s,1,1,0,0,0,p+1-s) \in \mathbb{Z}^7 | 1 \leq s \leq p\}$.

Theorem \ref{dgalgebraandhomology} leads us to the following combinatorial description of $\Upsilon^{\leq 1}$:

\begin{thm}\label{basisupsilon} $\Upsilon^{\leq 1}$ has basis $\{ m_v \}_{v \in \mathcal{P}_{\Upsilon^{\leq 1}} }$ with product given by
$$m_u m_{u'} =
\left\{ \begin{array}{ll}
(-1)^{a j'_0 + b j'_0 + b a'}m_{v} & \textrm{if }
v_1 = u_1, u_7 = u_1', u_7'=v_7,
v_l = u_l+u'_l \\ & \textrm{for } 2 \leq l \leq 5   \textrm{ and } v \in \mathcal{P}_{\Upsilon^{\leq 1}}. \\
0 & \textrm{otherwise.} \end{array} \right.$$
\end{thm}

This is precisely the algebra described in Section \ref{monomial}, and using Proposition \ref{stock}  as well as the definition of $\Upsilon = \bbH\bbT_{\fc}(\underline{\ft})$ and Lemma \ref{trunc} this completes the proof of Theorem \ref{Yoneda}.

\begin{example}\label{yuck}
Suppose that in the example computed in Section \ref{example}, we would like to see the basis element that gives the non-zero $\Ext^8(\Delta(1),\Delta(9))$ in the polytopal basis. In the construction this came from the basis element $13_{-2}^2$ in $\fc$ with the degree of both $\ttu, \ttv $ equal to zero in the description of Theorem  \ref{dgalgebraandhomology}, tensored with the basis element $13_{-6}^4$ of a copy of $M$, which belonging to $\ft^{-2}$ had a factor $\ttu^2$. If we express this in terms of the polytopal basis we obtain $a=b=0$ for the first, with $i=0$, $j=-2, k=2$, hence the polytopal element $(1,0,-2,2,0,0,3)$ for the first, and $a=2,b=0$, therefore $i=-2, j=-6-3+1 = -8, k=4+2=6$, and hence the polytopal element $(1,-2,-8,6,2,0,3)$ for the second.
Hence in the polytopal basis of $\fO_F \fO_\Upsilon^2(F[z,z^{-1}])$, this corresponds to the element
$((1,0,0,0,0,0,1),(1,0,-2,2,0,0,3),(1,-2,-8,6,2,0,3),z^8)$. Note that the $k$-degree of this element, as the sum over the $k$-degrees of the constituents, is given by the sum over the fourth entries and is indeed $8$. 

\end{example}

\section*{Appendix 1: Signs.} \label{signsappendix}

{\bf Super sign convention.} Here we record some aspects of the super sign convention that are of relevance for us.
A differential graded vector space is a $\mathbb{Z}$-graded vector space $V = \oplus_k V^k$ with a graded endomorphism $\ttd$ of degree $1$.
We write $|v|$ for the degree of a homogeneous element of $V$.
We assume $\ttd$ can act both on the left and the right of $V$, with the convention $\ttd(v) = (-1)^{|v|}(v)\ttd$.
A differential graded algebra is a $\mathbb{Z}$-graded algebra $A = \oplus_k A^k$ with a differential $d$ such that
$$\ttd(ab) = \ttd(a).b + (-1)^{|a|}a.\ttd(b),$$
or equivalently
$$(ab)\ttd = a.(b)\ttd + (-1)^{|b|}(a)\ttd.b.$$
If $A$ is a differential graded algebra then a differential graded left $A$-module is a graded left $A$-module $M$ with differential $\ttd$ such that
$$\ttd(a.m) = \ttd(a).m + (-1)^{|a|}a.\ttd(m);$$
a differential graded right $A$-module is a graded right $A$-module $M$ with differential $\ttd$ such that
$$\ttd(m.a) = \ttd(m).a + (-1)^{|m|}m.\ttd(a).$$
If $A$ and $B$ are dg algebras then a dg $A$-$B$-bimodule is a graded $A$-$B$-bimodule with a differential
which is both a left dg $A$-module and a right dg $B$-module.
If $_AM$ is a left dg $A$-module,
then $\End_A(M)$ is a differential graded algebra which acts on the right of $M$, giving $M$ the structure of an $A$-$\End_A(M)$-bimodule,
the differential on $\End_A(M)$ being given by $m.(\phi)\ttd = ((m)\phi)\ttd - (-1)^{|\phi|} ((m)\ttd)\phi$.
If $M_B$ is a right dg $A$-module,
then $\End_A(M)$ is a differential graded algebra which acts on the left of $M$, giving $M$ the structure of an $\End_B(M)$-$B$-bimodule,
the differential on $\End_B(M)$ being given by $\ttd(\phi).m =  \ttd(\phi.m) - (-1)^{|\phi|} \phi.\ttd(m)$.

If $_AM_B$ and $_BN_C$ are dg bimodules where $A$, $B$, and $C$ are dg algebras, then $M \otimes_B N$ is a dg $A$-$C$-bimodule
with differential $$\ttd(m \otimes n) = \ttd(m) \otimes n + (-1)^{|m|} m \otimes \ttd(n).$$

If $_AM_B$ and $_AN_C$ are dg bimodules where $A$, $B$, and $C$ are dg algebras, then $Hom_A(M, N)$ is a dg $B$-$C$-bimodule
with differential $$\ttd(\phi(m)) = \ttd(\phi)(m) + (-1)^{|\phi|} \phi(\ttd(m)).$$

\section*{Appendix 2: Koszul duality.}\label{kosz}

Here we give an account of Koszul duality for Koszul algebras,
synthesising the work of Beilinson, Ginzburg, and Soergel, and Keller \cite{BGS}, \cite{Ke}.

\medskip \emph{In what follows we will exceptionally denote homological degree by $h$ rather than $k$,
because the homological $h$-grading on the Koszul algebras described here is different from the homological
$k$-grading for the Koszul algebra $\fc$ used in the rest of the paper.}

\medskip 

{\bf Vector space duals.} Let $A^0$ be a direct product of finitely many fields, thought of as an algebra.
We write $M^*$ for the
$F$-linear graded dual of a graded vector space over $F$; the dual of a component in degree $j$ lies in degree $-j$;
if $M$ is an $A$-$B$-bimodule, then $M^*$ is a $B$-$A$-bimodule.
Given a right $A^0$-module $M$ with dual $M^*$, we have an isomorphism $(Me)^* \cong e(M^*)$ for each primitive idempotent $e$ in $A^0$;
we write $\eta_M$ for the sum $\eta: A^0 \rightarrow M \otimes_{A^0} M^*$ of units $F \cong A^0e \rightarrow Me \otimes_F eM^*$.
If $M$ is an $A^0$-$A^0$-bimodule, then $\eta$ is a homomorphism of $A^0$-$A^0$-bimodules.
We have a fixed isomorphism $A^0 \cong A^{0*}$ which sends $1 \in F$ to its dual in $F^*$.

\medskip 

{\bf Quadratic duals.} Let $A = \mathbb{T}_{A^0}(A^1)/R$, with $R \subset A^1 \otimes_{A^0} A^1$
be a quadratic algebra whose degree zero part is $A^0$ and whose degree one part $A^1$ is finite dimensional.
Let $A^! = \mathbb{T}_{A^{0}}(A^{!-1})/R^!$ be its quadratic dual, where the $A^0$-$A^0$ bimodules $A^1$ and $A^{!-1}$,
and the short exact sequences of $A^0$-$A^0$-bimodules
$$0 \rightarrow R \rightarrow A^1 \otimes_{A^0} A^1 \rightarrow A^2 \rightarrow 0$$
$$0 \leftarrow A^{!-2} \leftarrow A^{!-1} \otimes_{A^0} A^{!-1} \leftarrow R^{!} \leftarrow 0,$$
are duals of each other.
The grading that is implicit here is the \emph{radical grading}, or $r$-grading:
we insist $A$ is generated in $r$-degrees $0$ and $1$, and $A^!$ is generated in $r$-degrees $0$ and $-1$.

\medskip 

{\bf Differential bimodules.} The composition
$$\xymatrix{
A^0 \ar[rrr]^(.3){\eta_{A^{1} \otimes_{A^0} A^1}} & &
& (A^{1} \otimes_{A^0} A^1) \otimes_{A^0} (A^{!-1} \otimes_{A^0} A^{!-1}) \ar[r] & A^2 \otimes_{A^0} A^{!2}
}$$
is equal to zero, because the first map can be written $\sum_i (b_i \otimes b^*_i) + \sum_j (b^{!*}_i \otimes b^!_i)$ where
$\{b_i \}$ is a basis for $R$ with dual basis $\{ b_i^* \}$ and
$\{ b^!_i\}$ is a dual basis for $R^!$ with dual basis $\{ b_i^{!*}\}$.
Consequently the space $A \otimes_{A^0} A^{!}$ is a differential bimodule,
with differential given by the composition map
$$A \otimes_{A^0} A^!  \overset{\sim}{\longrightarrow}  A \otimes_{A^0} A^0 \otimes_{A^0} A^! \overset{1 \otimes \eta \otimes 1}{\longrightarrow} A \otimes_{A^0} A^1 \otimes_{A^0} A^{!-1} \otimes_{A^0} A^! \overset{}{\longrightarrow}A \otimes_{A^0} A^!$$
We denote this differential bimodule $C$.
This differential acts naturally on the \emph{inside} of $C$ which is a little awkward notationally:
we adopt the convention that this differential applied to $a \otimes \alpha$ is written
$\begin{array}{ccc}
& \ttd & \\
(a & \otimes & \alpha)
\end{array}$.

Suppose that $A$ is an $rh$-graded algebra with an $A^0$-$A^0$-bimodule decomposition
$A^1 = A^{10} \oplus A^{11}$; thus the $r$-degree $1$ part $A^1$ of $A$
decomposes as a direct sum of a $h$-degree $0$ part and a $h$-degree $1$ part. We have $A^{!-1} = A^{1*}$,
and write $A^{!-11} = A^{10*}$, $A^{!-10} = A^{11*}$.
Thus $A^{!-1} = A^{!-10} \oplus A^{!-11}$, and $A^!$ are $h$-graded algebras,
in such a way that the differential on $C$ has $h$-degree one.
If we want to write down the corresponding map on the left or right we apply the super sign convention:
$$\ttd(a \otimes \alpha) = (-1)^{|a|_h}
\begin{array}{ccc}
& \ttd & \\
(a & \otimes & \alpha)
\end{array}$$
$$(a \otimes \alpha)\ttd = (-1)^{|\alpha|_h}
\begin{array}{ccc}
& \ttd & \\
(a & \otimes & \alpha)
\end{array}$$
There are only two ways in which we can obtain a left and right differential this way:
either $A$ is concentrated in $h$-degree $0$, or $A^!$ is concentrated in $h$-degree $0$.
From now on we assume that one of these is the case.
In this way we give $C$ the structure of a differential $h$-graded $A$-$A^!$-bimodule.
We denote by $C^! = A^! \otimes_{A^0} A$ the corresponding differential bigraded $A^!$-$A$-bimodule.
The \emph{left Koszul complex} is the differential $rh$-bigraded $A$-$A^!$-bimodule
$$K_l = C \otimes_{A^!} A^{!*} \cong A \otimes_{A^0} A^{!*}.$$

\medskip 

{\bf Adjunction.} Given a pair of $rh$-graded modules $M$ and $N$ we define $\Hom(M,N)$ to be the sum of the spaces of
$rh$-graded homomorphisms from $M$ to $N$ shifted in degree by $(r,h)$.
We have
$$\Hom_A(A \otimes_{A^0} A^{!*}, M) \cong \Hom_{A^0}(A^{!*}, M) \cong
A^! \otimes_{A^0} M \cong A^! \otimes_{A^0} A \otimes_A M,$$
by adjunction and the fact that $A^0$ is semisimple.
Consequently, there is an equivalence of functors between categories of differential $rh$-bigraded modules
$$\Hom_A(K_l,-) \cong C^! \otimes_{A} -: A \bigra_{rh} \rightarrow A^! \bigra_{rh},$$
and therefore an adjunction $(K_l \otimes_{A^!}-, C^! \otimes_A -)$;
we have a homomorphism of differential $rh$-bigraded bimodules $\rho_l: K_l \otimes_{A^!} C^! \rightarrow A$
corresponding to the counit of this adjunction given by the composition of natural maps
$$A \otimes_{A^0} A^{!*} \otimes_{A^0} A \rightarrow A \otimes_{A^0} A^{0*} \otimes_{A^0} A \rightarrow A \otimes_{A^0} A \rightarrow A;$$
we have a homomorphism of differential $rh$-bigraded bimodules $\phi_l: A^! \rightarrow C^! \otimes_{A} K_l$ corresponding to the unit of this adjunction given by the composition of natural maps
$$A^! \rightarrow \Hom_A(K_l, K_l) \cong \Hom_{A^!}(A^!, C^! \otimes_A K_l)\cong C^! \otimes_{A} K_l.$$
Let us explain why the differentials on $C^!$ and $K_l$ match under the natural isomorphisms.
Here is a diagram depicting the counit of the adjunction
$$\xymatrix{
A \ar@{-} '[ddrrr] '[rrrrrr] & \otimes_{A^0} &  A^{!*} \ar@{-} '[dr] '[rr] &  \otimes & A^! & \otimes_{A^0} &  A \\
  &               &         &     F     &     &               &   \\
  &               &         &     A     &     &               &
}$$
The map sends $a \otimes \alpha \otimes \alpha' \otimes a'$ to $a \langle \alpha, \alpha' \rangle b$.
Under the differential on $K_l$ we obtain the map sending $a \otimes \alpha \otimes \alpha' \otimes a'$ to
$\sum_x ax \langle x^* \alpha, \alpha' \rangle b$;
under the differential on $C^!$ we obtain the map sending $a \otimes \alpha \otimes \alpha' \otimes a'$ to
$\sum_x a \langle  \alpha, \alpha'x^* \rangle x b$;
these maps are identical by the super sign convention since whenever $\langle \alpha, \alpha' \rangle$
is nonzero for homogeneous $\alpha$ and $\alpha'$, we have $|\alpha|_h + |\alpha'|_h = 0$, and $|x|_h|x^*|_h$ is always zero.

The algebra $A^0$ admits the structure of a differential bigraded $A$-$A^!$-bimodule:
we have zero differential, and all elements of strictly positive or strictly negative degrees act as zero.
There is a natural homomorphism of differential bigraded $A$-modules $\pi_l: K_l \rightarrow A^0$
given by $\pi_l = \rho_l \otimes_A 1_{A^0}$.
There is a natural homomorphism of differential bigraded $A$-modules $\iota_l: A^0 \rightarrow K_l$
given by $1_{A^0} \otimes_{A^!} \phi_l$.

\medskip 

{\bf Koszul algebras.}
The algebra $A$ is said to be \emph{Koszul} if $\pi_l$ is a quasi-isomorphism.
This is equivalent to $\iota_l$ being a quasi-isomorphism, or $\rho_l$ being a quasi-isomorphism, or $\phi_l$ being a quasi-isomorphism.
Since the Koszul complexes for $A$ and $A^!$ are duals of each other, $A$ is Koszul if and only if $A^!$ is Koszul.

We have a category $A \bigra_{rh}$ whose objects are differential $rh$-bigraded modules.
Localising the quasi-isomorphisms gives us a triangulated category $D(A \bigra_{rh})$.
If $A$ is Koszul then we have adjoint equivalences of derived categories of differential bigraded modules
$$\xymatrix{ D(A \bigra_{rh}) \ar@/^/[r]^{\Hom_A(K_l,-)} & \ar@/^/[l]^{K_l \otimes_{A^!}-} D(A^! \bigra_{rh}) }.$$

We now restrict to the case in which $A$ is concentrated in homological degree $0$.

We give $A$ the structure of a differential $rh$-bigraded algebra as follows:
$A$ is concentrated entirely in $h$-degree $0$, whilst $A^!$ is concentrated in positive degrees;
the $h$-grading on $A^!$ is the negative of the $r$-grading.

Writing $C^{r,h} = A^{r+h} \otimes_{A^0} A^{!-h}$, we give $C$ the structure of a differential bigraded $A$-$A^!$-bimodule;
the differential has $(r,h)$-degree $(0,1)$.
Writing $C^{! r,h} = A^{! -h} \otimes_{A^0} A^{r+h}$, we give $C^!$ the structure of a differential bigraded $A^!$-$A$-bimodule;
the differential has $(r,h)$-degree $(0,1)$.

The \emph{left Koszul complex} is the differential bigraded $A$-$A^!$-bimodule
$K_l = C \otimes_{A^!} A^{!*} \cong A \otimes_{A^0} A^{!*}$,
where the bigrading is given by
$$K_l^{r,h} = A^{r+h} \otimes_{A^0} A^{!* -h} = A^{r+h} \otimes_{A^0} A^{!h *};$$
we denote by $K_l^! = C^! \otimes_{A} A^*$ the left Koszul complex for $A^!$; the differential has $(r,h)$-degree $(0,1)$.

There is a natural homomorphism of differential bigraded $A$-modules $\pi_l: K_l \rightarrow A^0$
obtained by tensoring $1_C$ on the right with the dual of the embedding $A^0 \rightarrow A^!$
and on the left the homomorphism $A \rightarrow A^0$ which sends all elements of positive degree to zero.
The algebra $A$ is Koszul if $\pi_l$ is a quasi-isomorphism.
In that case, the structures described above collapse favourably:

\begin{thm} \label{Koszulgeneralities} (Beilinson, Ginzburg and Soergel \cite{BGS}, Keller \cite{Ke})
Suppose $A$ is Koszul. Then
the map $$A^! \rightarrow \Ext^\bullet_A(A^0,A^0)$$ induced by the action of $A^!$ on $K_l$ is an isomorphism;
if $A^!$ is finite dimensional, then we have adjoint equivalences of bounded derived categories
$$\xymatrix{ D^b(A \gr_r) \ar@/^/[r]^{\Hom_A(K_l,-)} & \ar@/^/[l]^{K_l \otimes_{A^!}-} D^b(A^! \gr_r)},$$
where $A \gr_r$ is the category of $r$-graded modules for $A$.
We have adjoint equivalences of derived dg categories
$$\xymatrix{ D_{dg}(A) \ar@/^/[r]^{\Hom_A(K_l,-)} & \ar@/^/[l]^{K_l \otimes_{A^!}-} D_{dg}(A^!) }.$$
where $D_{dg}(A)$ is the derived category of differential $h$-graded modules for $A$.
\end{thm}
\begin{remark}
Obviously
$$CC^!CC^!...CC^! \cong A \otimes_{A^0} A^! \otimes_{A^0} A \otimes_{A^0} ... \otimes_{A^0} A,$$
$$CC^!CC^!...CC^!C \cong A \otimes_{A^0} A^! \otimes_{A^0} A \otimes_{A^0} ... \otimes_{A^0} A^!,$$
where there is one more tensor factor on the right hand side of the isomorphisms than on the left hand side.
If $A$ is Koszul, then we have quasi-isomorphisms
$$K_lK_l^!K_lK_l^!...K_lK_l^! \leftarrow A^* \otimes_{A^0} A^{!*} \otimes_{A^0} ... \otimes_{A^0} A^*,$$
$$K_lK_l^!K_lK_l^!...K_lK_l^!K_l \leftarrow  A^* \otimes_{A^0} A^{!*} \otimes_{A^0} ... \otimes_{A^0} A^{!*},$$
where there is one fewer tensor factor on the right hand side of the quasi-isomorphism than on the left hand side
(here we write $MN$ for a tensor product $M \otimes_A N$).
The complex $CC^!CC^!...CC^!$ is quasi-isomorphic to the dual of $K_lK_l^!K_lK_l^!...K_lK_l^!$,
where there are two more tensor factors in the second term;
the complex $CC^!CC^!...CC^!C$ is quasi-isomorphic to the dual of $K_l^!K_lK_l^!...K_lK_l^!K_lK_l^!$,
where there are two more tensor factors in the second term.
\end{remark}
\vspace*{-5mm}
\begin{remark}
There is a right Koszul complex $K_r = A^* \otimes_{A} C$ which is obtained by tensoring $C$ on the left with $A^*$.
We have an adjunction $(C^! \otimes_{A} -, K_r \otimes_{A^!} -)$. We have a homomorphism of dg bimodules
$\phi_r: C^! \otimes_A K_r \rightarrow A^!$ corresponding to the counit of this adjunction.
We have a dg homomorphism
$\pi_r: A^0 \rightarrow K_r$ obtained by taking the dual of $\pi_l^!$, which is a quasi-isomorphism if and only if $A$ is Koszul.
We have an an analogue of Theorem \ref{Koszulgeneralities} for the right Koszul complex.
\end{remark}

\section{Appendix 3: Corrections to Section 5}\label{errata}

This section should replace Section 5.

\subsection{The $2$-category $\mathcal{T}$}
Let $\mathcal{T}$ denote the
collection of  pairs $(A, M)$ where $A$ is a differential $k$-graded algebra and $M$
is a  differential $k$-graded $A$-$A$-bimodule.

The collection $\mathcal{T}$ in fact forms the set of objects of a $2$-category:  $1$-morphisms  between two objects $(A,M)$ and $(B,N)$ are given by a pair $(S, \phi_S)$, consisting of a differential (bi-)graded $A$-$B$-bimodule ${}_AS_B$ and a quasi-isomorphism
$$\phi_S: S \otimes_BN\to M\otimes_AS; $$

$2$-morphisms from $(S, \phi_S)$ to $(T, \phi_T)$ are given by homomorphisms of differential (bi-)graded $A$-$B$-bimodules $f: S \to T$ such that the diagram
$$\xymatrix{
S \otimes_BN\ar^{\phi_S}[r]\ar^{ f\otimes id }[d]&M \otimes_A S \ar^{id \otimes f}[d] \\
T \otimes_BN\ar^{\phi_T}[r]                      &M \otimes_A ST                     
}$$
commutes.

\begin{defn}
We define a {\bf Rickard object} of  $\mathcal{T}$ to be an object
$(A,M)$ of $\mathcal{T}$, where ${}_AM_A \in A\perf \cap \rperf A$, the natural morphism of dg algebras $A\to \End_A(M)$ is a quasi-isomorphism, there is a quasi-isomorphism $A \to \bbH A$, and $\bbH A$ is a finite-dimensional algebra of finite global dimension. 
\end{defn}

For a Rickard object $(A,M)$ of $\mathcal{T}$, write $M^{-1}$ for $\Hom_A(M,A)$.

\begin{defn}
We define a {\bf $j$-graded
object} of $\mathcal{T}$ to be an object $(a,m)$ of $\mathcal{T}$,
where $a= \bigoplus a^{jk}$ is a differential bigraded algebra, and $m
= \bigoplus m^{jk}$ a differential bigraded $a$-$a$-bimodule, and $a^{j\bullet} = m^{j\bullet} = 0$ for $j<0$. 
\end{defn}

\subsection{The operator $\mathbb{O}$.} Let $(\bfa,\bfm) $ be a $j$-graded object of $\mathcal{T}$. We define
$$\mathbb{O}_{\bfa,\bfm} \circlearrowright \mathcal{T}$$
to be the operator given by
$$\mathbb{O}_{\bfa,\bfm}(A,M) = (\bfa(A,M), \bfm(A,M)),$$
where
$$\alpha(A,M) =(\alpha^0 \otimes A)\oplus (\bigoplus_{j > 0} \alpha^j \otimes M^{\otimes_A j})$$ for $\alpha \in \{\bfa,\bfm\}.$ 
The algebra structure on $\bigoplus \bfa^{jk} \otimes_F M^{\otimes_A j}$ is the restriction of the algebra structure on the tensor
product of algebras $\bfa \otimes \bbT_{A}(M)$, where $ \bbT_{A}(M)$ denotes the tensor algebra of $M$ over $A$.
The $k$-grading and differential on the complex $\bigoplus \bfa^{jk} \otimes M^{\otimes_A j}$ are defined to be the total $k$-grading and total differential on the tensor product of complexes. The bimodule structure, grading and differential on $\bigoplus \bfm^{jk} \otimes M^{\otimes_A j}$ are defined likewise. 

We remark that this extends to a $2$-endofunctor of $\mathcal{T}$ (cf. \cite[Lemma 9]{MT2}).
\begin{lem} \label{tensor}\cite[Lemma 14]{MT3}
Let $\bfa, \bfb, \bfc$ be a differential bigraded algebras, ${}_\bfb \x_\bfa$ and ${}_\bfa \y_\bfc$
differential bigraded modules, all concentrated in nonnegative $j$-degrees. Let  $(A, M)$ be an object of $\mathcal{T}$.
Then
$$\x(A,M) \otimes_{\bfa(A,M)} \y(A,M) \cong (\x \otimes_\bfa \y)(A,M)$$ as differential bigraded $\bfb(A,M)$-$\bfc(A,M)$-bimodules.
\end{lem}

Given a differential bigraded $\bfa$-module $\x$, with components in
positive and negative $j$-degrees, we define $\x(A,M)$ to be the
$\bfa(A,M)$-module given by
$$\x(A,M) = (\bigoplus_{j<0} x^{j\bullet} \otimes (M^{-1})^{\otimes_A -j}) \oplus (x^{0\bullet} \otimes A) \oplus
(\bigoplus_{j>0} x^{j\bullet} \otimes M^{\otimes_A j}),$$
where $M^{-1}:=\Hom_A(M,A)$.

\begin{lem} \label{homo}(cf.\cite[Lemma 15]{MT3})
Let $\bfc$ be a differential bigraded algebra, $\x$ and $\y$ are differential bigraded
$\bfc$-modules, all concentrated in nonnegative $j$-degrees, and let $(A,M)$ be a Rickard
object of $\mathcal{T}$. Then we have a quasi-isomorphism of differential bigraded $(\bfc^0\otimes A)\otimes (\bfc^0\otimes A)^{op}$-modules
$$\Hom_\bfc(\x,\y)(A,M) \rightarrow \Hom_{\bfc(A,M)}(\x(A,M), \y(A,M)).$$
\end{lem}

\proof
We have previously stated this only as a quasi-isomorphism of differential bigraded vector spaces. However, the quasi-isomorphism we constructed in \cite[Proof of
Theorem 13]{MT2} is in fact a quasi-isomorphism of $(\bfc^0\otimes A)\otimes (\bfc^0\otimes A)^{op}$-modules.
\endproof

\subsection{ The operator $\fO$.} We now recall the definition of the operator $\fO$ from \cite{MT3}. Let $\Gamma = \bigoplus \Gamma^{ijk}$ be a
differential trigraded algebra. We have an operator
$$\fO_\Gamma \circlearrowright \{ \Sigma | \mbox{ $\Sigma = \bigoplus \Sigma^{jk}$ a differential bigraded algebra } \}$$
given by
\begin{equation}\label{fO}\fO_\Gamma(\Sigma)^{ik} = \bigoplus_{j, k_1+k_2=k} \Gamma^{ijk_1} \otimes \Sigma^{jk_2}.\end{equation}
The algebra structure and differential are obtained by restricting
the algebra structure and differential from $\Gamma \otimes \Sigma$.
If we forget the differential and the $k$-grading, the operator
$\fO_\Gamma$ is identical to the operator $\fO_\Gamma$ defined in the
introduction.

\subsection{Comparing $\mathbb{O}$ and $\fO$.}
Throughout this section, let $(\bfa,\bfm)$ be a $j$-graded object of $\mathcal{T}$ and $(A,M)\in \mathcal{T}$.
Note that the algebra $\bbT_\bfa(\bfm)(A,M)$, formed with respect to the $j$-grading
on $\bbT_\bfa(\bfm)$, is a differential bigraded algebra, with
$$\bbT_\bfa(m)(A,M)^{ik} = \oplus_{j} \bbT_\bfa(\bfm)^{ijk} \otimes M^{\otimes_A
j}.$$ The algebra $\bbT_{\bfa(A,M)}(\bfm(A,M))$ is a differential bigraded
algebra, with $$\bbT_{\bfa(A,M)}(\bfm(A,M))^{ik} =
(\bfm(A,M)^{\otimes_{\bfa(A,M)} i})^k.$$

We write $X^{i \diamond \bullet} \cong Y^{i \diamond \bullet}$ to signify that $X^{i jk} \cong Y^{ijk}$ for all $j,k$. 

We have the following lemmas.

\begin{lem}\label{properties}\cite[Lemma 19]{MT3}
\begin{enumerate}[(i)]
\item\label{8i} We have an isomorphism of objects of $\mathcal{T}$
$$\bbO_{\bfa,\bfm}(A,M) \cong
(\fO_{\bbT_\bfa(\bfm)}(\bbT_A(M))^{0 \diamond \bullet},
\fO_{\bbT_\bfa(\bfm)}(\bbT_A(M))^{1 \diamond \bullet}),$$ where the $k$-grading
on the components of $\bbO_{\bfa,m}(A,M)$ can be identified with the
$k$-grading on $\fO_{\bbT_\bfa(\bfm)}(\bbT_A(M))$.
\item\label{8ii} We have an isomorphism of differential bigraded algebras
$$\fO_{\bbT_\bfa(\bfm)}(\bbT_A(M)) \cong \bbT_\bfa(\bfm)(A,M).$$
\item\label{8iii} We have an isomorphism of differential bigraded algebras
$$\bbT_\bfa(\bfm)(A,M) \cong \bbT_{\bfa(A,M)}(\bfm(A,M)).$$
\end{enumerate}
\end{lem}

Suppose we are given $(\bfa_i,\bfm_i)$ for $1 \leq i \leq n$, and $(A,M)$.
Let us define $(A_i,M_i)$ recursively via $(A_i,M_i) =
\bbO_{\bfa_i,\bfm_i}(A_{i-1},M_{i-1})$ and $(A_0,M_0) = (A,M)$.

\begin{lem} \label{comparetwo}\cite[Lemma 20]{MT3}
\begin{enumerate}[(i)]
 \item \label{9i} We have an algebra isomorphism
$$\bbT_{A_n}(M_n) \cong \fO_{\bbT_{\bfa_n}(\bfm_n)}
...\fO_{\bbT_{\bfa_1}(\bfm_1)}(\bbT_A(M)).$$
\item\label{9ii} We have an isomorphism
of objects of $\mathcal{T}$
\begin{equation*}\begin{split}\bbO_{\bfa_n,\bfm_n} &\cdots\bbO_{\bfa_1,\bfm_1}(A,M) \\
&\cong (\fO_{\bbT_{\bfa_1}(\bfm_1)}...\fO_{\bbT_{\bfa_n}(\bfm_n)}(\bbT_A(M))^{0\diamond \bullet},
\fO_{\bbT_{\bfa_1}(\bfm_1)}...\fO_{\bbT_{\bfa_n}(\bfm_n)}(\bbT_A(M))^{1\diamond \bullet}).
\end{split}\end{equation*}
\end{enumerate}
\end{lem}

\begin{cor}\label{compareagain}\cite[Corollary 21]{MT3}
Let $\bfa$ be a dg algebra and $\bfm$ dg $\bfa$-$\bfa$-bimodule. Then we have an
isomorphism of algebras
$$\bbO_{F,0} \bbO^n_{\bfa,\bfm}(F,F) \cong \fO_F\fO_{\bbT_\bfa(\bfm)}^n(F[z,z^{-1}]).$$
\end{cor}

\subsection{Keller duality.}
From now on, we always assume that our Rickard object $(A,M)$ has an additional vector space $d$-grading, such that the $d$-degree $0$ part of $A$ generates the derived category $D(A)$. Let $P_A$ be a projective resolution of the $d$-degree $0$ part of $A$ and let $\cF (A)$ denote the dg algebra
$\cF(A) = \End_{A}(P_A)$. There are mutually inverse equivalences
$$\xymatrix{ D_{dg}(A) \ar@/^/[r]^{\Hom_A(P_A,-)} & \ar@/^/[l]^{P_A \otimes_{\cF (A)}-} D_{dg}(\cF (A))},$$
by a theorem of Keller \cite[Theorem 3.10]{Ke}. Since $P_A$ is
projective as an $A$-module, we have a natural isomorphism of
functors
$$\Hom_A(P,-) \cong \Hom_A(P_A,A) \otimes_A -.$$
Given a differential graded $A$-$A$-bimodule $N$, denote by $\cF(M)$ the dg $\cF(A)$-$\cF(A)$-bimodule
$$\cF(N) = \Hom_A(P_A,A) \otimes_A N \otimes_A P_A \cong \Hom_A(P_A, N\otimes_AP_A).$$

\begin{lem}\label{invout}\begin{enumerate} We have isomorphisms
\item $\cF(M^{-1}) \cong \cF(M)^{-1}$;
\item $\cF(M^{\otimes_A i}) \cong \cF(M)^{\otimes_{\cF(A) i}}$ for $i\in \ZZ$.
\end{enumerate}
\end{lem}

\proof For simplicity, write $P=P_A$
On the on hand, we have 

\begin{equation*}
\begin{split}
\cF(M^{-1}) &\cong \Hom_A(P,\Hom_A(M,A)\otimes_AP)\\&\cong \Hom_A(P,\Hom_A(M,P)) \cong \Hom_A(M\otimes_A P,P)
\end{split}
\end{equation*}
by definition, $A$-projectivity of $M$ and adjunction.
On the other hand, we have 
\begin{equation*}
\begin{split}\cF(M)^{-1} &= \Hom_{\cF(A)}(\Hom_A(P,A)\otimes_A M\otimes_AP, \Hom_A(P,P)) \\&\cong \Hom_A(P\otimes_{\cF(A)} \Hom_A(P,A)\otimes_A M\otimes_AP, P).\end{split}
\end{equation*}
Since $P$ is a projective generator for $A$, we have $P\otimes_{\cF(A)} \Hom_A(P,A) = P\otimes_{\End_A(P)} \Hom_A(P,A) \cong A$ and (i) follows.
(ii) is proved similarly.
\endproof

\begin{defn}\label{dagger}
Suppose that $(a,\um) = (a,(m, m^{-1}))$ is a $j$-graded Rickard object of $\mathcal{U}$, such that $a$ is concentrated in non-negative $d$-degrees and that a projective resolution $P$ of the degree zero part of $a$ is differential $djk$-trigraded.
Then $\cF(a)$ inherits a differential $djk$-trigrading.
In this case we call $(a,\um)$ a \emph{dagger object} of $\mathcal{U}$.
\end{defn}

\subsection{The $d$-grading}

While the $d$-grading as defined in the article is a vector space grading with the property that taking the $d$-degree zero part of $\bbH\bbO^q_{\bfc,\bft}(F,F)$ gives the direct sum of standard modules for a block of polynomial representations of $G$ with $p^q$ simple modules (and more generally, the graded pices provide a standard filtration of the corresponding algebra), the grading defined on $\bft$ is not a module grading over the $d$-graded algebra $\bfc$, and as such the $d$-grading on $\bbH\bbO^q_{\bfc,\bft}(F,F)$ is not an algebra grading. An algebra grading is not needed anywhere in the article, since the iterative construction only carries around the $d$-grading as an adornment, and can and should be phrased with the $d$-grading being a vector space grading only (apart from in Definition \ref{dagger}), with the only assumption being that the differential is of $d$-degree $0$.

\subsection{Bonding.}

%\begin{defn}
%We say a Rickard object
%$(A,M)$ is \emph{bonded} if the quasi-isomorphism $A\to \End_A(M)$ has an $A$-$A$-bimodule splitting.
%\end{defn}

\begin{lem}\label{hbond}
The natural maps, called bonding maps,  $M\otimes_A\Hom_A(M,A) \to A$ (given by evaluation) and $\Hom_A(M,A)\otimes_A M \to \End_A(M) \overset{qim}{\leftarrow} A$ induce a structure of associative algebra on $$\bbH\bbT_A(\uM) = \bbH\left((\bigoplus_{i \geq 1} M^{\otimes_A i}) \oplus A \oplus (\bigoplus_{i \leq -1} (M'^{-1})^{\otimes_A -i})\right).$$
\end{lem}

\proof
Setting $E:=\End_A(M)$, we have a commutative diagram
$$\xymatrix{
M \otimes_A M^{-1}\otimes_AM\ar_\psi[d]\ar^\phi[rr] && M \otimes_A E\ar^{\theta}@{->>}[d] && M \otimes_A A\ar@{=}[dd]  \ar_{\eta}[ll]    \\
A\otimes_AM\ar@{=}[d]                              && M\otimes_E E \ar@{=}[d]     &&\\
M\ar@{=}[rr]&&M\ar@{=}[rr]&&M
}$$
where $\phi, \psi$ are the bonding maps, $\eta:m\otimes 1\mapsto m\otimes 1$ and $\theta: m\otimes e \mapsto m\otimes e$ are the natural maps induced by the quasi-isomorphism $A\to E$, and the equalities denote the canonical isomorphisms.
Indeed in the square on the left $x\otimes f\otimes y$ maps to $f(x)y$ either way, and in the right square, $m\otimes 1$ gets sent to $m$ either way. In homology, $\eta$ and $\theta$ become isomorphisms, and $\bbH(\eta)^{-1}\circ \phi$ produces a map that makes multiplication associative.

Similarly, the commutative diagram
$$\xymatrix{
M^{-1} \otimes_A M\otimes_AM^{-1}\ar_\phi[d]\ar^\psi[rr] && E \otimes_A M^{-1}\ar^{\theta}@{->>}[d] && A \otimes_A M^{-1}\ar@{=}[dd]  \ar_{\eta}[ll]    \\
M^{-1}\otimes_AA\ar@{=}[d]                              && E\otimes_E M^{-1} \ar@{=}[d]     &&\\
M^{-1}\ar@{=}[rr]&&M^{-1}\ar@{=}[rr]&&M^{-1}
}$$
gives the desired associativity in homology.

\endproof

%\begin{lem}\label{Kbond}
%If $(A,M)$ is a  bonded dagger object, then so is $(\cF(A),\cF(M))$.
%\end{lem}
%
%\proof
%We need to check that the quasi-isomorphism $\cF(A)\to \End_{\cF(A)}(\cF(M))$ has a $\cF(A)$-$\cF(A)$-bimodule splitting.
%We have 
%\begin{equation*}
%\begin{split}
%\Hom_{\cF(A)}(\cF(M),\cF(M)) &= \Hom_{\cF(A)}(\Hom_{A}(P_A,M\otimes_AP_A),\Hom_{A}(P_A,M\otimes_AP_A))\\
%&\cong  \Hom_{A}(P_A\otimes_A\Hom_{A}(P_A,A)\otimes_AM\otimes_AP_A,M\otimes_AP_A)\\
%&\cong \Hom_{A}(M\otimes_AP_A,M\otimes_AP_A)\\
%&\cong \Hom_{A}(P_A,\Hom_A(M,M)\otimes_AP_A)\\
%&\rightarrow^{qim}  \Hom_{A}(P_A,A\otimes_AP_A) = \cF(A)
%\end{split}\end{equation*}
%where we have used the bonding on $(A,M)$. We claim that this is the desired bimodule splitting.
%\endproof
%Similarly, if $\bfa$ is a not necessarily positively $j$-graded, and $\bfm$ is a differential bigraded $\bfa$-$\bfa$-bimodule, we obtain the structure of a $jk$-graded algebra on $\bbH(\bfa(A,M))$.

%\begin{lem}\label{itbond}
%Let $(\bfa,\bfm)$ be a $j$-graded bonded Rickard object of $\mathcal{T}$ and $(A,M)$ a bonded Rickard object in $\mathcal{T}$. Then $(\bfa(A,M),\bfm(A,M))$ is a bonded Rickard object in $\mathcal{T}$.
%\end{lem}
%\proof
%We know it is a Rickard object by ??? and check that the bonding survives.
%We have 
%\begin{equation*}
%\begin{split}
%
%
%\end{split}\end{equation*}
%\endproof
%

\subsection{Iteration}

For a Rickard object $(A,M)$, define
$$\bbF\bbT_A(M):= \bigoplus_{i\in\ZZ} \cF(M^{\otimes_A i}).$$

Note that this is a priori only a $\cF(A)$-$\cF(A)$ bimodule but, after taking homology we obtain

\begin{equation}\label{use}
\begin{split}
\bbH\bbF\bbT_A(M)&= \bbH\bigoplus_{i\in\ZZ} \cF(M^{\otimes_A i})\\
&\cong \bbH\bigoplus_{i\in\ZZ} \cF(M)^{\otimes_{\cF(A) i}}\\
&\cong \bbH\bbT_{\cF(A)}(\cF(\uM))
\end{split}\end{equation}
which is an associative algebra thanks to Lemmas \ref{hbond} and \ref{invout}, with the algebra structure induced by the bonding maps on $\bbH\bbT_A(M)$.

\begin{thm}\label{iteration}
Assume the following:
\begin{enumerate}[(i)]
\item\label{cond2} Let $(\bfa,\bfm)$ be a dagger object of $\mathcal{T}$ such that $\bfa$ is concentrated in non-negative $j$-degrees.
\item\label{cond1} Let $(A,M) = (A,M)$ be a Rickard object of $\mathcal{T}$, such that
\begin{enumerate}[(a)]
\item\label{cond11} both $A$ and $M$ are concentrated in non-negative $d$-degrees;
\item\label{cond12} $M^{\otimes j} \otimes_A A^{0 \bullet} \cong (M^{\otimes j})^{0 \bullet} $ for all $j$ such that $a^{0 j \bullet}\neq 0$, where $()^{0 \bullet}$ refers to taking the $d$-degree $0$ part.
\end{enumerate}
\item\label{cond3} $(\bfa(A,M))^{0 \bullet}$ generates the derived category $D_{dg}(\bfa(A,M))$.
\end{enumerate}

Then $$\bbH\bbF\bbT_{\bfa(A,M)}(\bfm(A,M))\cong \fO_{\bbH\bbT_{\cF(\bfa)}(\cF(\um))}\bbH\bbT_{\cF(A)}( \cF(\uM))$$
\end{thm}

\proof
Since $$\bbH\bbF\bbT_{\bfa(A,M)}(\bfm(A,M))=\bbH\bigoplus_{i\in\ZZ} \cF(\bfm(A,M)^{\otimes_{\bfa(A,M)} i}),$$ we do this for each $i$ separately.
We can choose a projective resolution of the $d$-degree $0$ part of $\bfa(A,M)$ as $P_\bfa(A,M)\otimes_AP_A$ (see paper).

\begin{equation*}
\begin{split}
\cF(&\bfm(A,M)^{\otimes_{\bfa(A,M)} i})\\&=\Hom_{\bfa(A,M)}(P_\bfa(A,M)\otimes_AP_A,\bfm(A,M)^{\otimes_{\bfa(A,M) i}}\otimes_{\bfa(A,M)}P_\bfa(A,M)\otimes_AP_A )
\\
&=\Hom_A(P_A, \Hom_{\bfa(A,M)}(P_\bfa(A,M) , \bfm(A,M)^{\otimes_{\bfa(A,M) i}}\otimes_{\bfa(A,M)}P_\bfa(A,M)\otimes_AP_A  ))\\
&=\Hom_A(P_A, \Hom_{\bfa(A,M)}(P_\bfa(A,M) , \bfm(A,M)^{\otimes_{\bfa(A,M) i}}\otimes_{\bfa(A,M)}P_\bfa(A,M))\otimes_AP_A  )
\end{split}\end{equation*}

We now consider $$\Hom_{\bfa(A,M)}(P_\bfa(A,M) , \bfm(A,M)^{\otimes_{\bfa(A,M) i}}\otimes_{\bfa(A,M)}P_\bfa(A,M)).$$

If $i\geq 0$, then $\bfm(A,M)^{\otimes_{\bfa(A,M) i}}\otimes_{\bfa(A,M)}P_\bfa(A,M) \cong (\bfm^{\otimes_\bfa i}\otimes_\bfa P_\bfa)(A,M)$ by Lemma \ref{tensor}, so by Lemma \ref{homo}, we obtain a quasi-isomorphism
\begin{equation*}
\begin{split}
\Hom_{\bfa(A,M)}&(P_\bfa(A,M) , \bfm(A,M)^{\otimes_{\bfa(A,M) }i}\otimes_{\bfa(A,M)}P_\bfa(A,M))\\
&\leftarrow \Hom_\bfa(P_\bfa, \bfm^{\otimes_\bfa i}\otimes_\bfa P_\bfa) (A,M),
\end{split}\end{equation*}

and hence

\begin{equation*}
\begin{split}
\bbH\Hom_A&(P_A, \Hom_{\bfa(A,M)}(P_\bfa(A,M) , \bfm(A,M)^{\otimes_{\bfa(A,M) }i}\otimes_{\bfa(A,M)}P_\bfa(A,M))\otimes_AP_A  )\\
&\cong \bbH\Hom_A(P_A,  \Hom_\bfa(P_\bfa, \bfm^{\otimes_\bfa i}\otimes_\bfa P_\bfa) (A,M)\otimes_A P_A)\\
&\cong \bbH\Hom_A(P_A,  \cF(\bfm^{\otimes_\bfa i})(A,M)\otimes_A P_A)\\
&\cong \bbH\bigoplus_{j\in \ZZ} \Hom_A(P_A,\cF(\bfm^{\otimes_\bfa i})^j\otimes M^{\otimes_A j}\otimes_AP_A)\\
&\cong \bbH\bigoplus_{j\in \ZZ}\cF(\bfm^{\otimes_\bfa i})^j\otimes \Hom_A(P_A, M^{\otimes_A j}\otimes_AP_A)\\
&\cong\bbH\bigoplus_{j\in \ZZ}\cF(\bfm^{\otimes_\bfa i})^j\otimes \cF(M^{\otimes_Aj})\\
&\cong \bigoplus_{j\in \ZZ}\bbH\cF(\bfm^{\otimes_\bfa i})^j\otimes \bbH(\cF(M)^{\otimes_\cF(A)j})\\
%&\cong \cF(\bfm^{\otimes_\bfa i})(\cF(A),\cF(M)).
&= \fO_{\bbH\bbT_{\cF(\bfa)}(\cF(m), \cF(m^{-1}))}\bbH\bbT_{\cF(A)}(\cF(A), \cF(\uM))^{i\bullet\diamond}
\end{split}\end{equation*}

For $i<0$

\begin{equation*}
\begin{split}
\Hom_{\bfa(A,M)}&(P_\bfa(A,M) , \bfm(A,M)^{\otimes_{\bfa(A,M) }i}\otimes_{\bfa(A,M)}P_\bfa(A,M))\\
&=\Hom_{\bfa(A,M)}(P_\bfa(A,M) , \Hom_{\bfa(A,M)}(\bfm(A,M)^{\otimes_{\bfa(A,M)} -i}, \bfa(A,M))\otimes_{\bfa(A,M)}P_\bfa(A,M))\\
&\cong \Hom_{\bfa(A,M)}(\bfm(A,M)^{\otimes_{\bfa(A,M)}-i}\otimes_{\bfa(A,M)}P_\bfa(A,M), P_\bfa(A,M))\\
&\cong \Hom_{\bfa(A,M)}( (\bfm^{\otimes_\bfa -i} \otimes_\bfa P_\bfa)(A,M), P_\bfa(A,M))\\
&\leftarrow^{qim} \Hom_\bfa(\bfm^{\otimes_\bfa -i} \otimes_\bfa P_\bfa, P_\bfa) (A,M),
\end{split}\end{equation*}

so we obtain 

\begin{equation*}
\begin{split}
\bbH\Hom_A&(P_A, \Hom_{\bfa(A,M)}(P_\bfa(A,M) , \bfm(A,M)^{\otimes_{\bfa(A,M) }i}\otimes_{\bfa(A,M)}P_\bfa(A,M))\otimes_AP_A  )\\
&\cong \bbH\Hom_A(P_A, \Hom_\bfa(\bfm^{\otimes_\bfa -i} \otimes_\bfa P_\bfa, P_\bfa) (A,M)\otimes_AP_A)\\
&\cong \bbH\Hom_A(P_A,  \cF(\bfm^{\otimes_\bfa i})(A,M)\otimes_A P_A)
\end{split}\end{equation*}

from where the proof continues as in the $i>0$ case.

The fact that this is an isomorphism of algebras follows from the algebra structure on both sides being naturally induced by the bonding maps on $(\bfm,\bfm^{-1})$ and $(M,M^{-1})$.  Indeed, the algebra structure on $\bbH\bbF\bbT_{\bfa(A,M)}(\bfm(A,M))$ is induced, thanks to Lemma \ref{invout}, by that on $\bbH\bbT_{\bfa(A,M)}(\bfm(A,M))$ which is induced by viewing it as a subalgebra inside $\bbH\bbT_\bfa(\bfm) \otimes\bbH\bbT_A(M)$. On the other hand, $\fO_{\bbH\bbT_{\cF(\bfa)}(\cF(\um))}\bbH\bbT_{\cF(A)}(\cF(A), \cF(\uM))$ inherits its algebra structure by viewing it inside $ \bbH\bbT_{\cF(\bfa)}(\cF(\um))\otimes\bbH\bbT_{\cF(A)}(\cF(\uM))$ and the algebra structure of the latter is, again by Lemma \ref{invout}, induced by that on $\bbH\bbT_\bfa(\bfm) \otimes\bbH\bbT_A(M)$.
\endproof

This replaces Theorem 18 in the paper, and only this statement is applied in the proof of Proposition 28, which in general form can be phrased as follows.

\begin{cor}\label{red}
Let $(\bfa,\bfm)$ is a dagger object in $\mathcal{T}$ and $(A,M) = \bbO^q_{\bfa,\bfm}(F,F)$. Setting $(A_{k}, M_{k}) =  \bbO^{k}_{\bfa,\bfm}(F,F)$ for $k=1,\dots q$, assume that the hypotheses of Theorem \ref{iteration} are satisfied  for  $(\bfa,\bfm)$  and all $(A_{k}, M_{k})$.
Then $$\bbH\bbF\bbT_A(M) \cong \fO^q_{\bbH\bbT_{\cF(\bfa)}(\cF(\um))}F[z,z^{-1}].$$ Moreover, this yields an isomorphism $$\bbH\cF(A)\cong\fO_F\fO^q_{\bbH\bbT_{\cF(\bfa)}(\cF(\um))}F[z,z^{-1}].$$
\end{cor}

\proof
We prove this by induction: in case $q=0$ both sides equal $F[z,z^{-1}]$.
%Setting $(A_{k}, M_{k}) =  \bbO^{k}_{\bfa,\bfm}(F,F)$ for $k=1,\dots q$, and noting that by Lemma \ref{itbond} all these are bonded Rickard objects,  
We then obtain
\begin{equation*}\begin{split}
\bbH\bbF\bbT_A(M) &\cong \bbH\bbF\bbT_{\bfa(A_{q-1},M_{q-1})}(\bfm(A_{q-1},M_{q-1}))\\
&\cong \fO_{\bbH\bbT_{\cF(\bfa)}(\cF(\um))}\bbH\bbT_{\cF(A_{q-1})}(\cF(A_{q-1}), \cF(\uM_{q-1}))\quad \hbox{by Theorem \ref{iteration}}\\
&\cong  \fO_{\bbH\bbT_{\cF(\bfa)}(\cF(\um))}\bbH\bbF\bbT_{A_{q-1}}(M_{q-1}) \quad \hbox{by \eqref{use}}\\
&\cong  \fO_{\bbH\bbT_{\cF(\bfa)}(\cF(\um))} \fO^{q-1}_{\bbH\bbT_{\cF(\bfa)}(\cF(\um))}F[z,z^{-1}] \quad \hbox{by the inductive assumption}\\
&=\fO^q_{\bbH\bbT_{\cF(\bfa)}(\cF(\um))}F[z,z^{-1}].
\end{split}\end{equation*}

The second statement follows by applying the operator $\fO_F$, and noticing that when applying it to $\bbH\bbF\bbT_A(M)$, we obtain precisely $\bbH\cF (A)$.
\endproof

Proposition 28 in the paper is the direct application of this corollary to the algebras relevant for $GL_2$, where, setting $\bfa=\bfc$ and $\bfm=\bft$, the algebra $\bbH\cF (A)$ is the algebra that is there called $\mathbf{w}_q$.

\subsection{Indexing parameters $a$ and $b$}\label{indexing}

In Collary 39, we missed an extension as $\fc$-$\fc$-bimodules. of $\overline{M}^\tau$ by $\overline{\fc^{0\sigma}}$. The generator $e_p\otimes e_1$ in $\overline{M}^\tau$ extends the element $\xi e_p\otimes e_1 \xi \in \overline{\fc^{0\sigma}}$, while only the other basis elements of $\overline{\fc^{0\sigma}}$, namely $\xi e_h\otimes e_{p+1-h} \xi$ for $2\leq h \leq p-1$, split off as $\fc$-$\fc$-bimodules in homology. As a consequence, in the polytopal basis, any basis element $(p-1,0,0,a,b,2) \in \mathcal{P}_{\leq 0}$ (in the portion coming from $\mathcal{P}_{\overline{\fc^{0\sigma}}}$) should satisfy $a=b+1$ while the basis elements $(p-s,0,0,a,b,s+1) \in \mathcal{P}_{\leq 0}$ for $s=2, \dots,p-1$ are correctly given as satisfying $a=b-1$. This does not affect the number of elements in the basis, nor the $ijk$-degree of the elements, but is relevant for the multiplication as stated in Theorem 53 to be correct.

\vfill

\section*{Index of Notation}
\renewcommand{\arraystretch}{1.3}
\begin{tabular}{p{2cm}p{6cm}p{3.5cm}}
{\bf Symbol}&{\bf  Explanation} &{\bf  Location}\\ &&\\
$\mathbf{W}$ & Yoneda extension algebras of all Weyl modules for $GL_2$ &Introduction\\
$\mathbf{w}$ & Yoneda extension algebras of Weyl modules in fixed block for $GL_2$ &Introduction\\ 
$i,j,k,d$&-gradings & Section \ref{grs}\\
$\langle 1\rangle$ & shift in $j$-degree & Section \ref{grs}\\ 
$[1]$ &shift in $k$-degree & Section \ref{grs}\\ 
$\bbH$ &homology functor & Section \ref{grs}\\ 
$\underline{-}$ & refers to a bonded pair & Section \ref{bonded}\\ 
$\mathbb{T}_{a}(\um)$ & tensor algebra & Section \ref{bonded}\\ 
$\mathcal{T}, \mathcal{U}$& ambient $2$-categories & Section \ref{2cats}\\ 
$\bbP$& algebraic operator & Section \ref{bbP}\\
$\fO$ & algebraic operator & Section \ref{SfO}\\
$\mathbb{D}$& $2$-functor forgetting the $d$-grading & Section \ref{comparison}\\
$\cF,\bbF$ & homological duality operators & Section \ref{bbF}\\
$Z$ & zigzag algebra & Section \ref{zigzag}\\
$\bfc$ & algebra encoding $GL_2$ & Section \ref{zigzag}\\
$\bft$ & tilting module for $\bfc$ & Section \ref{zigzag}\\
$\ut=(\tilde \bft,\tilde \bft^{-1}) $& projective bimodule resolutions of $\bft$ and its adjoint& Section \ref{zigzag}\\
$\fc$ & $\cF(\bfc)$& Section \ref{psi}\\
$\ttd$ & differential & \ref{tthroughK}\\
$\ft$ & $\cF(\tilde\bft)$ &\eqref{maltese}, end of Section \ref{tthroughK}\\
$\ft^{-1}$ & $\cF(\tilde\bft^{-1})$ &\eqref{malteseinv}, Section \ref{triple}\\
$\underline{\ft}$ & $(\ft, \ft^{-1})$  & Section \ref{triple}\\
$\Upsilon$ & $\mathbb{HT}_{\fc}(\underline{\ft})$ & Section \ref{tensoralg}\\
$\ttu,\ttv,\ttw, \ttz$ & formal variables in $\Upsilon$ & Section \ref{tensoralg}\\
\end{tabular}

\normalfont

{\sc Vanessa Miemietz, Will Turner}\\
School of Mathematics, University of East Anglia, Norwich, NR4 7TJ, UK, \\{\tt v.miemietz@uea.ac.uk}\\
Department of Mathematics, University of Aberdeen, Fraser Noble Building, King's College, Aberdeen AB24 3UE, UK, {\tt w.turner@abdn.ac.uk}.

\begin{thebibliography}{00}
\normalfont

\bibitem{BGS} A. Beilinson, V. Ginzburg, W. Soergel, \emph{Koszul duality patterns in representation theory}, J. Amer. Math. Soc.  9  (1996),  no. 2, 473--527.

\bibitem{CPS}  E. Cline, B. Parshall and L. Scott, \emph{Finite-dimensional algebras and highest weight categories},
J. Reine Angew. Math. 391 (1988), 85--99.

%\bibitem{Erd} K. Erdmann, \emph{Schur algebras of finite type}, Quart. J. Math. Oxford Ser. (2) 44 (1993), no. 173, 17--41. 

\bibitem{EH} K. Erdmann, A. Henke, \emph{On Ringel duality for Schur algebras},
Math. Proc. Cambridge Philos. Soc. 132 (2002), no. 1, 97--116. 

\bibitem{Green} J. A. Green, \emph{Polynomial representations of ${\rm
GL}_{n}$.} Lecture Notes in Mathematics, 830. Springer, Berlin,
1980.

\bibitem{Ke}
B. Keller, \emph{On differential graded categories}, International Congress of
  Mathematicians. Vol. II, Eur. Math. Soc., Z\"urich, 2006, pp.~151--190.

\bibitem{KS} M. Khovanov, P. Seidel, \emph{Quivers, Floer cohomology, and braid group actions}
J. Amer. Math. Soc. 15 (2002) 203--271.

\bibitem{Mad} D. Madsen, \emph{On a common generalization of Koszul duality and tilting equivalence}, preprint, arxiv:1007.3282v1, to appear in Adv. Math.

\bibitem{MT1} V. Miemietz, W. Turner, \emph{Rational representations of $GL_2$}, Glasgow J. Math. 53 (2011), no.2, 257--275.

\bibitem{MT2} V. Miemietz, W. Turner, \emph{Homotopy, Homology and $GL_2$}, Proc. London Math. Soc.(3) 100 (2010), no.2, 585--606.

\bibitem{MT3} V. Miemietz, W. Turner, \emph{Koszul dual $2$-functors and extension algebras of simple modules for $GL_2$}, preprint, arXiv:1106.5411v2.

\bibitem{AP} A. Parker, \emph{Higher extensions between modules for $SL_2$},  Adv. Math.  209  (2007),  no. 1, 381--405.

\bibitem{Rickard}  J. Rickard, \emph{Derived equivalences as derived functors.} J. London Math. Soc. (2) 43 (1991), no. 1, 37--48.

\bibitem{RR} R. Rouquier, \emph{Derived equivalences and finite
dimensional algebras}, Proceedings of the International Congress of
Mathematicians (Madrid, 2006 ), vol II, pp. 191-221, EMS Publishing
House, 2006.






\end{thebibliography}
\end{document}